\documentclass{amsart}

\usepackage{latexsym,amsmath,amstext,mathrsfs,amsthm,amscd,amsopn,verbatim,amsrefs}

\usepackage{graphicx}
\usepackage[dvips]{color}
\usepackage{epsfig}
\usepackage{psfrag}
\usepackage[all]{xy}
\usepackage{amscd}
\usepackage{amssymb}

\newtheorem{Thm}{Theorem}[section]
\newtheorem{Prop}[Thm]{Proposition}
\newtheorem{Lem}[Thm]{Lemma}
      
\newtheorem{Conj}[Thm]{Conjecture}

\theoremstyle{remark}
\newtheorem{Rem}[Thm]{Remark}

\theoremstyle{definition}

\newtheorem{Def}[Thm]{Definition}
\newtheorem{Exa}[Thm]{Example}

\newcommand{\norm}[1]{|\!|#1|\!|}

\begin{document}

\title[Bimodules and branes in deformation quantization]{Bimodules and branes in deformation quantization}
\author{Damien Calaque}
\email{calaque@math.univ-lyon1.fr}
\address{Institut Camille Jordan, Universit\'e Claude Bernard Lyon 1, 43
boulevard du 11 novembre 1918, F-69622 Villeurbanne Cedex France}
\author{Giovanni Felder}
\email{giovanni.felder@math.ethz.ch}
\address{Department of Mathematics, ETH Zurich, R\"amistrasse 101
8092 Zurich, Switzerland}
\author{Andrea Ferrario}
\email{andrea.ferrario@math.ethz.ch}
\address{Department of Mathematics, ETH Zurich, R\"amistrasse 101,
8092 Zurich, Switzerland}
\author{Carlo A. Rossi}
\email{crossi@math.ist.utl.pt}
\address{Centro de An\'{a}lise Matem\'{a}tica, Geometria e Sistemas Din\^{a}micos,
Departamento de Matem\'atica, Instituto Superior T\'ecnico, Av. Rovisco Pais, 1049-001
Lisboa, Portugal}

\keywords{Deformation quantization; coisotropic submanifolds; Koszul algebras; Koszul duality; $L_\infty$-algebras and morphisms; $A_\infty$-bimodules}
\thanks{We are grateful to Alberto S. Cattaneo, David Kazhdan, Bernhard Keller, Thomas Tradler and
Thomas Willwacher for useful comments, discussions and
suggestions. This research is partly supported by the french national
agency through the ANR project G\'eSAQ (project number JC$08\underline{~}320699$, by SNF Grant 200020-122126 and by the Funda\c{c}\~{a}o para a Ci\^{e}ncia e a Tecnologia (FCT /
Portugal).}

\maketitle

\begin{abstract}
We prove a version of Kontsevich's formality theorem for two
subspaces (branes) of a vector space $X$. The result implies in
particular that the Kontsevich deformation quantizations of
$\mathrm{S}(X^*)$ and $\wedge(X)$ associated with a quadratic Poisson
structure are Koszul dual. This answers an open question in
Shoikhet's recent paper on Koszul duality in deformation
quantization.
\end{abstract}

\section{Introduction}\label{s-0}
Kontsevich's proof of his formality theorem \cite{K} is based on the
Feynman diagram expansion of a topological quantum field theory. In
\cite{CF-br} a program to extend Kontsevich's construction by
including branes (i.e., submanifolds defining boundary conditions
for the quantum fields) is sketched. The case of one brane leads to
the relative formality theorem \cite{CF} for the Hochschild cochains
of the sections of the exterior algebra of the normal bundle of a
submanifold and is related to quantization of Hamiltonian reduction
of coisotropic submanifolds in Poisson manifolds. Here we consider
the case of two branes in the simplest situation where the branes
are linear subspaces $U$, $V$ of a real (or complex) vector space
$X$. The new feature is that one should associate to the
intersection $U\cap V$ an $A_\infty$-bimodule over the algebras
associated with $U$ and $V$. The formality theorem we prove holds
for the Hochschild cochains of an $A_\infty$-category corresponding
to this bimodule. It is interesting that even if $U=\{0\}$ and
$V=X=\mathbb R$ the $A_\infty$-bimodule is one-dimensional but has
infinitely many non-trivial structure maps.

Our discussion is inspired by the recent paper of B. Shoikhet
\cite{Sh} who proved a similar formality theorem in the framework of
Tamarkin's approach based on Drinfeld associators. Our result
implies that Shoikhet's theorem on Koszul duality in deformation
quantization also holds for the explicit Kontsevich quantization.
Next we review the question of Koszul duality in Kontsevich's
deformation quantization, explain how it fits in the setting of
formality theorems and state our results.

\subsection{Koszul duality}\label{ss-0-1}
Let $X$ be a real or complex finite dimensional vector space. Then
it is well known that the algebra $B=\mathrm S(X^*)$ of polynomial
functions on $X$ is a quadratic Koszul algebra and it is Koszul dual
to the exterior algebra $A=\wedge(X) $. In \cite{Sh} Shoikhet studied
the question of quantization of Koszul duality. He asked whether the
Kontsevich deformation quantization of $A$ and $B$ corresponding to
a quadratic Poisson bracket leads to Koszul dual formal associative
deformations of $A$ and $B$. Recall that a quadratic Poisson
structure on a finite dimensional vector space $X$ is by definition
a Poisson bracket on $B=\mathrm{S}(X^*)$ with the property that the
bracket of any two linear functions is a homogeneous quadratic
polynomial. A quadratic Poisson structure on $X$ also defines by
duality a (graded) Poisson bracket on $A=\wedge(X)$. If $x^1,\dots,
x^n$ are linear coordinates on $X$ and $\theta_1,\dots,\theta_n$ is
the dual basis of $X$, the brackets of generators have the form
\[
\{x^i,x^j\}_B=\sum_{k,l}C^{i,j}_{k,l}x^kx^l,\qquad
\{\theta_k,\theta_l\}_A=\sum_{i,j}C^{i,j}_{k,l}\theta_i\wedge\theta_j.
\]
Kontsevich gave a universal formula for an associative star-product
$f\star g=fg+\hbar B_1(f,g)+\hbar^2 B_2(f,g)+\cdots$ on
$\mathrm{S}(X^*)[\![\hbar]\!]$ such that $B_1$ is any given Poisson
bracket. Universal means that $B_j(f,g)$ is a differential
polynomial in $f,g$ and the components of the Poisson bivector field
with universal coefficients. Kontsevich's result also applies to
super manifolds such as the odd vector space $W=X^*[1]$, in which
case $\mathrm{S}(W^*)=\mathrm{S}(X[-1])=\wedge(X)$. Moreover, if the
Poisson bracket is quadratic, then the deformed algebras
$A[\![\hbar]\!]$, $B[\![\hbar]\!]$ are quadratic, namely they are generated
by $\theta^i$, resp.~$x_i$ with quadratic defining relations.
Shoikhet proves that Tamarkin's \cite{Tamarkin} universal
deformation quantization corresponding to any Drinfeld associator
leads to Koszul dual quantizations. Here we show that the same is
true for the original Kontsevich deformation quantization.

\subsection{Branes and bimodules}\label{ss-0-2}
In Kontsevich's approach, the associative deformations of $A$ and
$B$ are given by explicit formul\ae\ involving integrals over
configuration spaces labelled by Feynman diagrams of a topological
quantum field theory. We approach the question of Koszul duality
from the quantum field theory point of view, following a variant of
a suggestion of Shoikhet (see \cite{Sh}, 0.7). The setting is the
theory of quantization of coisotropic branes in a Poisson manifold
\cite{CF-br}. In this setting, quantum field theory predicts the
existence of an $A_\infty$-category whose set of objects $S$ is any
given collection of submanifolds (``branes'') of a Poisson manifold.
If $S$ consists of one object one obtains the $A_\infty$-algebra
related to Hamiltonian reduction \cite{CF}. Here we consider the
next simplest case of two objects that are subspaces $U,V$ of a
finite dimensional vector space $X$. In this case the
$A_\infty$-category structure is given by two (possibly curved)
$A_\infty$-algebras $A$, $B$ and an $A_\infty$-$A$-$B$-bimodule $K$ over $\mathbb R[\![\hbar]\!]$ or $\mathbb C[\![\hbar]\!]$; the $A_\infty$-algebras 
represent the spaces of endomorphism of the two objects $U$, $V$ respectively, while the $A_\infty$-$A$-$B$-bimodule represents the space of morphism from $U$ to $V$.
More precisely, we have $A=\Gamma(U,\wedge (\mathrm NU))[\![\hbar]\!]=\mathrm{S}(U^*)\otimes\wedge(X/U)[\![\hbar]\!]$,
$B=\Gamma(V,\wedge (\mathrm NV))[\![\hbar]\!]=\mathrm{S}(V^*)\otimes\wedge(X/V)[\![\hbar]\!]$, the
sections of the exterior algebras of the normal bundles, and 
\begin{equation}\label{e-Maan}
K=\Gamma(U\cap V,\wedge (\mathrm TX/(\mathrm TU+\mathrm
TV)))[\![\hbar]\!]=\mathrm{S}(U\cap V)\otimes \wedge (X/(U+V))[\![\hbar]\!].
\end{equation}
The structure maps of these algebras and bimodule are compositions
of morphisms in the $A_\infty$-category and are described by sums
over graphs with weights given by integrals of differential forms
over configuration spaces on the upper half-plane. The differential
forms are products of pull-backs of propagators, which are one-forms
on the configuration space of two points in the upper half-plane.
Additionally to the Kontsevich propagator \cite{K}, which vanishes
when the first point approaches the real axis, there are three
further propagators with brane boundary conditions \cite{CF-br,
CattaneoTorossian}. The four propagators obey the four possible
boundary conditions of vanishing if the first or second point
approaches the positive or negative real axes. In the physical model
these are the Dirichlet boundary conditions for coordinate functions
of maps from the upper half plane to $X$ such that the positive real
axis is mapped to a coordinate plane $U$ and the negative real axis
to a coordinate plane $V$.

The new feature here is that even for zero Poisson structure the
$A_\infty$-bimodule has non-trivial structure maps. Let us describe
the result first in the simplest case $U=\{0\}, V=X$ so that
$A=\wedge(X)$, $B=\mathrm S(X^*)$, $K=\mathbb R$ (here it is not
necessary to tensor by $\mathbb R[\![\hbar]\!]$ since the structure maps
are independent of $\hbar$).

\begin{Prop}\label{p-Lata}
Let $A$ be the graded associative algebra $A=\wedge(X)=\mathrm S(X[-1])$ with
generators of degree $1$ and $B=\mathrm S(X^*)$ concentrated in degree $0$. 
View $A$ and $B$ as $A_\infty$-algebras with Taylor components products $\mathrm d^j=0$ except for
$j=2$. Then there exists an $A_\infty$-$A$-$B$-bimodule $K$ whose
structure maps
 \[
\mathrm d_K^{j,k}\colon A[1]^{\otimes j}\otimes K[1]\otimes B[1]^{\otimes k}\to
K[1],
\]
obey $\mathrm d_K^{1,1}(v\otimes k\otimes u)=\langle u,v\rangle k$ for $k\in
K$, $v\in X\subset\wedge(X)$ and $u\in X^*\subset \mathrm{S}(X^*)$
and $\langle\ ,\ \rangle$ is the canonical pairing. In the general
case of subspaces $U,V\subset X$, where $A$ is generated by
$W_A=U^*\oplus (X/U)[1])$ and $V$ by $W_B=V^*\oplus (X/V)[1]$,
$\mathrm d_K^{1,1}(v\otimes k\otimes u)=\langle v,u\rangle k$ for
$v\in(V/(U\cap V))^*\oplus U/(U\cap V)[1]\subset W_B$ and $u\in
(U/(U\cap V))^*\oplus V/(U\cap V)[1] \subset W_A$.
\end{Prop}

The remaining $\mathrm d_K^{i,j}$ are given by explicit finite dimensional
integrals corresponding to the graphs depicted in Fig.~5, see
Section \ref{s-5}. There should exist a more direct description of
this basic object.

\begin{Exa} If $X$ is one-dimensional, $A=\mathbb R[\theta]$, $B=\mathbb
R[x]$ with $\theta^2=0$ the non-trivial structure maps of $K$
on monomials are
\[
\mathrm d_K^{j,1}(\underset{j}{\underbrace{\theta\otimes \cdots\otimes \theta}}\otimes 1\otimes x^j)=1;
\]
in this case, they can be computed inductively from the $A_\infty$-$A$-$B$-bimodule relations, using that $\mathrm d_K^{1,1}$ is simply the duality pairing between the generators $\theta=\partial_x$ and $x$.
\end{Exa}

\begin{Conj} The bimodule of Prop.~\ref{p-Lata} is
$A_\infty$-quasi-isomorphic to the Koszul free resolution $\wedge
(X^*)\otimes \mathrm{S}(X^*)$ of the right $\mathrm{S}(X^*)$-module
$\mathbb R$, where $\wedge(X)$ acts from the left by contraction.
\end{Conj}

\subsection{Formality theorem}\label{ss-0-3}
Our main result is a formality theorem for the differential graded
Lie algebra of Hochschild cochains of the $A_\infty$-category
associated with the $A_\infty$-$A$-$B$-bimodule $K$ (for zero
Poisson structure). Let thus as above $U$, $V$ be vector subspaces
of $X$, the objects of the category, and $A=\Gamma(U,\wedge (\mathrm
NU))=\texttt{Hom}(U,U)$, $B=\Gamma(V,\wedge (\mathrm
NV))=\texttt{Hom}(V,V)$, $K=\Gamma(U\cap V,\wedge (\mathrm
TX/(\mathrm TU+\mathrm TV)))=\texttt{Hom}(V,U)$,
$\texttt{Hom}(U,V)=0$. The nonzero composition maps in this
$A_\infty$-category are the products on $A$ and $B$ and the
$A_\infty$-bimodule maps $\mathrm d_K^{k,l}\colon A[1]^{\otimes k}\otimes
K\otimes B[1]^{\otimes l}\to K[1]$, $k,l\geq 0,i+j\geq1$. Let us call
this category $\texttt{Cat}_\infty(A,B,K)$. As for any $A_\infty$-category,
its shifted Hochschild cochain complex
$\mathrm C^{\bullet+1}(\texttt{Cat}_\infty(A,B,K))$ is a graded Lie algebra with
respect to the (obvious extension of the ) Gerstenhaber bracket.
Moreover there are natural projections to the differential graded
Lie algebras $\mathrm C^{\bullet+1}(A,A)$, $\mathrm C^{\bullet+1}(B,B)$ of
Hochschild cochains of $A$ and $B$. By Kontsevich's formality
theorem, these differential graded Lie algebras are
$L_\infty$-quasi-isomorphic to their cohomologies, that are both
isomorphic to the Schouten Lie algebra
$T^{\bullet+1}_{\mathrm{poly}}(X)=\mathrm S(X^*)\otimes
\wedge^{\bullet+1} X$ of poly-vector fields on $X$.

Thus we have a diagram of $L_\infty$-quasi-isomorphisms
\begin{equation}\label{e-Saeda}
\xymatrix{& \mathrm C^{\bullet+1}(A,A) & \\
T^{\bullet+1}_{\mathrm{poly}}(X)\ar[ur]^{\mathcal U_A}\ar[dr]_{\mathcal U_B} & & \ar[ul]_{\mathrm p_A}\ar[dl]^{\mathrm p_B}\mathrm C^{\bullet+1}(\texttt{Cat}_\infty(A,B,K))\\
& \mathrm C^{\bullet+1}(B,B) &}
\end{equation}


\begin{Thm}\label{t-Malati}
There is an $L_\infty$-quasi-isomorphism
$T^{\bullet+1}_{\mathrm{poly}}(X) \to
\mathrm C^{\bullet+1}(\texttt{Cat}_\infty(A,B,K))$ completing \eqref{e-Saeda} to a
commutative diagram of $L_\infty$-morphisms.
\end{Thm}

The coefficients of the $L_\infty$-morphisms are given by integrals
over configuration spaces of points in the upper half plane of
differential forms similar to Kontsevich's but with different
(brane) boundary conditions. This ``formality theorem for pairs of
branes'' is an $A_\infty$ analogue of Shoikhet's formality theorem
\cite{Sh}, who considered the case $U=\{0\}$, $V=X$ and $K$ replaced
by the Koszul complex and used Tamarkin's $L_\infty$-morphism
instead of Kontsevich's. Theorem \ref{t-Malati} follows from Theorem
\ref{t-form-cat} which is formulated and proved in Section
\ref{s-6}.

\subsection{Maurer--Cartan elements}\label{ss-0-4}
An $L_\infty$-quasi-isomorphism $\mathfrak g_1^\bullet\to\mathfrak
g_2^{\bullet}$ induces an isomorphism between the sets $MC(\mathfrak
g_i)=\{x\in \hbar\mathfrak g_i^1[\![\hbar]\!], \,
dx+\frac12[x,x]=0\}/\exp(\hbar \mathfrak g_i^0[\![\hbar]\!])$ of equivalence
classes of Maurer--Cartan elements (shortly, MCE), see \cite{K}. MCEs 
in $T_{\mathrm{poly}}(X)$ are formal Poisson structures on
$X$. They are mapped to MCEs in $\mathrm C^\bullet(A,A)$
and $\mathrm C^\bullet(B,B)$, which are $A_\infty$-deformations of the
product in $A$ and $B$. The previous theorem implies that the image of a
Poisson structure in $X$ in $\mathrm C^\bullet(A,B,K)$ is an
$A_\infty$-bimodule structure on $K[\![\hbar]\!]$ over the $A_\infty$-algebras $A[\![\hbar]\!]$, $B[\![\hbar]\!]$.

\subsection{Keller's condition}\label{ss-0-5}
The key property of the bimodule $K$, which is preserved under
deformation and implies the Koszul duality and the fact that the
projections $\mathrm p_A, \mathrm p_B$ are quasi-isomorphisms, is that it obeys an
$A_\infty$-version of {\em Keller's condition} \cite{Keller}. Before
formulating it, we introduce some necessary notions, see Section
\ref{s-3} for more details.

Recall that an $A_\infty$-algebra over a commutative unital ring $R$
is a $\mathbb Z$-graded free $R$-module $A$ with a codifferential
$\mathrm d_A$ on the counital tensor $R$-coalgebra $\mathrm{T}(A[1])$. The
{\em DG category of right $A_\infty$-modules} over an
$A_\infty$-algebra $A$ has as objects pairs $(M,\mathrm d_M)$ where $M$ is a
$\mathbb Z$-graded free $R$-module and $\mathrm d_M$ is a codifferential on
the cofree right $\mathrm{T}(A[1])$-comodule
$FM=M[1]\otimes_R\mathrm{T}(A[1])$. The complex of morphisms
$\underline{\mathrm{Hom}}_{-A}(M,N)$ is the graded $R$-module whose
degree $j$ subspace consists of homomorphisms $FM\to FN$ of
comodules of degree $j$, with differential $\phi\mapsto
\mathrm d_N\circ\phi-\phi\circ \mathrm d_M$. In particular
$\underline{\mathrm{End}}_{-A}(M)=\underline{\mathrm{Hom}}_{-A}(M,M)$,
for any module $M$, is a differential graded algebra. If $A$ is an
ordinary associative algebra and $M$, $N$ are ordinary modules, the
cohomology of $\underline{\mathrm{Hom}}_{-A}(M,N)$ is the direct sum
of the Ext-groups $\mathrm{Ext}^i_{-A}(M,N)$. The DG category of left
$A$-modules is defined analogously; its morphism spaces are denoted
$\underline{\mathrm{Hom}}_{A-}(M,N)$. If $A$ and $B$ are
$A_\infty$-algebras, an $A_\infty$-$A$-$B$-bimodule structure on $K$ is the same as a
codifferential on the cofree
$\mathrm{T}(A[1])-\mathrm{T}(B[1])$-comodule $T(A[1])\otimes K[1]\otimes
T(B[1])$ namely a codifferential compatible with coproducts and
codifferentials $\mathrm d_A$, $\mathrm d_B$.

The {\em curvature} of an $A_\infty$-algebra $(A,\mathrm d_A)$ is the
component $\mathrm F_A\in A^2$ in $A[1]=T^1(A[1])$ of $\mathrm d_A(1)$ where $1\in
R=T^0(A[1])$. If $\mathrm F_A$ vanishes then $\mathrm d_A(1)=0$ and $A$ is called
{\em flat}. If $A$ and $B$ are flat then an
$A_\infty$-$A$-$B$-bimodule is in particular an $A_\infty$ left
$A$-module and an $A_\infty$ right $B$-module. The left action of $A$
then induces a {\em derived left action}
\[
\mathrm L_A\colon A\to \underline{\mathrm{End}}_{-B}(K),
\]
which is a morphism of $A_\infty$-algebras (the differential graded
algebra $\underline{\mathrm{End}}_{-B}(K)$ is considered as an $A_\infty$-algebra 
with two non-trivial structure maps, the differential and
the product). Similarly we have a morphism of $A_\infty$-algebras
\[
\mathrm R_B\colon B\to \underline{\mathrm{End}}_{A-}(K)^\mathrm{op}.
\]
We say that an $A_\infty$-$A$-$B$-bimodule $K$, for flat $A_\infty$-algebras
$A$, $B$, obeys the Keller condition if $\mathrm L_A$ and $\mathrm R_B$ are
quasi-isomorphisms.

\begin{Lem} The bimodule $K$ of Prop.~\ref{p-Lata} obeys the
Keller condition.
\end{Lem}

An $A_\infty$-version of Keller's theorem \cite{Keller} that we
prove in Section \ref{s-3}, see Theorem \ref{t-keller} states that
if $K$ obeys the Keller condition then $\mathrm p_A$ and $\mathrm p_B$ in
\eqref{e-Saeda} are quasi-isomorphisms. Moreover the Keller
condition is an $A_\infty$-version of the Koszul duality of $A$ and
$B$ and reduces to it in the case of $U=\{0\}$, $V=X$ and quadratic
Poisson brackets, for which both $A$ and $B$ are ordinary
associative algebras. Indeed in this case $\mathrm L_A$ and $\mathrm R_B$ induce
algebra isomorphisms $B\cong\mathrm{Ext}_{A-}^\bullet(K,K)^\mathrm{op}$,
$A\cong\mathrm{Ext}_{-B}^\bullet(K,K)^\mathrm{op}$.

\subsection{The trouble with the curvature} Let us again consider
the simplest case $U=\{0\}$, $V=X$, and suppose that $\pi$ is a
Poisson bivector field on $X$. Then Kontsevich's deformation
quantization gives rise to an associative algebra
$(B_\hbar=\mathrm{S}(X^*)[\![\hbar]\!],\star_B)$ and a possibly curved
$A_\infty$-algebra $(A_\hbar=\wedge(X)[\![\hbar]\!],\mathrm d_{A_\hbar})$, both over
$\mathbb R[\![\hbar]\!]$. The one-dimensional $A$-$B$-bimodule $K$
deforms to an $A_\infty$-$A_\hbar$-$B_\hbar$-bimodule $K_\hbar$. If we restrict the
structure maps of this bimodule to $K_\hbar\otimes\mathrm{T}(B_\hbar)$ we get a deformation $\circ\colon K_\hbar\otimes
B_\hbar\to B_\hbar$ of the right action of $B$ as only non-trivial
map. However, this is not an action: instead we get
\[
(k\circ b_1)\circ b_2-k\circ(b_1\star b_2)=\langle \mathrm F_{A_\hbar},\mathrm d b_1\wedge \mathrm d b_2\rangle k.
\]
The curvature $\mathrm F_{A_\hbar}$ is a formal power series in $\hbar$ whose
coefficients are differential polynomials in the components of the
Poisson bivector field evaluated at zero. Its leading term vanishes
if $\pi(0)=0$ (i.e., if $V$ is coisotropic). The next term is
proportional to $\hbar^3$. It represents an obstruction to the
quantization of the augmentation module over $\mathrm{S}(X^*)$. T.
Willwacher \cite{Willwacher} constructed an example of a zero of a
Poisson bivector field on a five-dimensional space, whose module
over the Kontsevich deformation of the algebra of functions cannot
be deformed. On the other hand, there are several interesting
examples of Poisson structures such that $\mathrm F_{A_\hbar}=0$. Apart from
quadratic Poisson structures there are many examples related to Lie
theory, which we will study elsewhere.

\subsection{Organization of the paper} After fixing our notation
and conventions in Section \ref{s-1}, we recall the basic notions of
$A_\infty$-categories and their Hochschild cochain complex in
Section \ref{s-2}. In Section \ref{s-3} we formulate an
$A_\infty$-version of Keller's condition and extend Keller's theorem
to this case. In Section \ref{s-4} integrals over configuration spaces
of differential forms with brane boundary conditions are described.
The differential graded Lie algebra of Hochschild cochains of an
$A_\infty$-category is discussed in Section \ref{s-5}. Our main
result and its consequences are presented and proved in Section
\ref{s-6}.

\subsection*{Acknowledgements}
We are grateful to Alberto Cattaneo, David Kazhdan, Bernhard Keller, Thomas Tradler and
 Thomas Willwacher for useful comments, discussions and
suggestions. This work been partially supported by SNF Grant
200020-122126 and by the European Union through the FP6 Marie Curie
RTN ENIGMA (contract number MRTN-CT-2004-5652).

\section{Notation and conventions}\label{s-1}
We consider a ground field $k$ of characteristic $0$, e.g.\ $k=\mathbb R$ or $\mathbb C$.

Further, we consider the category $\texttt{GrMod}_k$ of $\mathbb Z$-graded vector spaces over $k$: we only observe that morphisms are meant to be linear maps of degree $0$, and we use the notation $\mathrm{hom}(V,W)$ for the space of morphisms. We denote by $\texttt{Mod}_k$ the full subcategory of $\texttt{GrMod}_k$ with objects being the ones 
concentrated in degree $0$. We denote by $[\bullet]$ the degree-shifting functor on $\texttt{GrMod}_k$.

The category $\texttt{GrMod}_k$ is a symmetric tensor category: the tensor product $V\otimes W$ (where, by abuse of notation, we do not write down the explicit dependence on the ground field $k$), for two general objects of $\texttt{GrMod}_k$, is the tensor product of $V$ and $W$ as $k$-vector spaces, with the grading induced by
\[
(V\otimes W)_p=\bigoplus_{m+n=p}V_m\otimes W_n,\ p\in\mathbb Z\,.
\] 
The symmetry isomorphism $\sigma$ is given by ``signed transposition'' 
$$
\sigma_{V,W}\,:V\,\otimes W\,\longrightarrow\, W\otimes V\,;\,
v\otimes w\,\longmapsto\,(-1)^{|v||w|}w\otimes v\,.
$$

Observe finally that the category $\texttt{GrMod}_k$ has inner Hom's: given two graded vector spaces $V,W$ one can consider the graded vector space $\mathrm{Hom}(V,W)$ 
defined by 
$$
\mathrm{Hom}^i(V,W)=\mathrm{hom}(V,W[-i])=\bigoplus_{k\in\mathbb Z}\mathrm{hom}_{\texttt{Mod}_k}(V_k,W_{k+i}),\ i\in\mathbb Z\,.
$$
Concretely, it will mean that we always assume tacitly Koszul's sign rule when dealing with linear maps between graded vector spaces: e.g. 
\begin{eqnarray*}
(\phi\otimes \psi)(v\otimes w)=(-1)^{|\psi||v|}\phi(v)\otimes\psi(w)\,.
\end{eqnarray*} 

\begin{Rem}\label{rem-completion}
We will sometimes deal with a variant $\widehat{\texttt{GrMod}_k}$ of this category where the Hom-spaces are replaced by their completion: 
$$
\widehat{\mathrm{Hom}}^i(V,W)=\widehat{\mathrm{hom}}(V,W[-i])=\prod_{k\in\mathbb Z}\mathrm{hom}_{\texttt{Mod}_k}(V_k,W_{k+i}),\ i\in\mathbb Z\,.
$$
\end{Rem}

The identity morphism of a general object $V$ of the category $\texttt{GrMod}_k$ induces an isomorphism $s:M\rightarrow M[1]$ of degree $-1$, which is called 
{\it suspension}; its inverse $s^{-1}:M[1]\rightarrow M$, which has obviously degree $1$, called {\it desuspension}. 
It is standard to denote by $|\cdot|$ the degree of homogeneous elements of objects of $\texttt{GrMod}_k$: recalling the definition of suspension and desuspension, 
we get $|s(\bullet)|=|\bullet|-1$. 

For a general object $V$ of $\texttt{GrMod}_k$, we denote by $\mathrm T(V):=\bigoplus_{n\in\mathbb{Z}}V^{\otimes n}$ the graded counital tensor coalgebra cogenerated by $V$: the counit is the canonical projection onto $V^{\otimes 0}=k$ and the coproduct is given by
$$
\Delta(v_1|\cdots|v_n)= 1\otimes (v_1|\cdots|v_n)+\sum_{j=1}^{n-1}(v_1|\cdots|v_{j})\otimes (v_{j+1}|\cdots|v_{n})+(v_1|\cdots|v_n)\otimes 1\,;
$$
where, for the sake of simplicity, we denote by $(v_1|\cdots|v_n)$ the tensor product $v_1\otimes\dots\otimes v_n$ in $V^{\otimes n}$.

Further, the symmetric algebra $\mathrm S(V)$ is defined as $\mathrm S(V)=\mathrm T(V)/\left\langle (v_1|v_2) -(-1)^{|v_1||v_2|} (v_2|v_1):\ v_1,v_2\in V\right\rangle$.
A general, homogeneous element of $\mathrm S(V)$ will be denoted by $v_1\cdots v_n$, $v_i$ in $V$, $i=1,\dots,n$.
The symmetric algebra is endowed with a coalgebra structure, with coproduct given by 
\begin{eqnarray*}
\Delta_{sh}(v_1\cdots v_n)=\sum_{p+q=n}\sum_{\sigma\in\mathfrak S_{p,q}}\epsilon(\sigma,v_1,\dots,v_n)(v_{\sigma(1)}\dots v_{\sigma(p)})\otimes (v_{\sigma(p+1)}\dots v_{\sigma(n)})\,,
\end{eqnarray*}
where $\mathfrak S_{p,q}$ is the set of $(p,q)$-shuffles, i.e.~permutations $\sigma\in\mathfrak S_{p+q}$ such that $\sigma(1)<\dots\sigma(p)$ and $\sigma(p+1)<\dots<\sigma(n)$, with corresponding sign 
\begin{eqnarray}
\epsilon(\sigma,v_1,\dots,v_n)=(-1)^{\sum_{i<j, \sigma(i) >\sigma(j)} |\gamma_i||\gamma_j| }\,, \label{shuffle}
\end{eqnarray}
and counit specified by the canonical projection onto $k$.

We define further the cocommutative coalgebra of invariants on $V$ as $\mathrm C(V)=\bigoplus_{n\geq 0} \mathrm I_n(V)$, with $\mathrm I_n(V)=\{x\in V^{\otimes n}:\ x=\sigma x, \forall \sigma \in\mathfrak S_n \}$: it is a sub-coalgebra of $\mathrm T(V)$, with coproduct given by the restriction of the natural coproduct onto and standard counit. 
We define also the cocommutative coalgebra without counit as $\mathrm C^{+}(V)=\mathrm C(V)/k$: we have an obvious isomorphism of coalgebras $\mathrm{Sym}:\mathrm S(V)\to \mathrm C(V)$, explicitly given by
\[
\mathrm S(V)\ni v_1\cdots v_n\overset{\mathrm{Sym}}\mapsto \frac1{n!}\sum_{\sigma\in S_n}\epsilon(\sigma,v_1,\dots,v_n) (v_{\sigma(1)}|\cdots|v_{\sigma(n)})\in \mathrm C(V)\,.
\]


Finally, we need to consider the category $\texttt{GrMod}_k^{I\times I}$ of $I\times I$-graded objects in $\texttt{GrMod}_k$, where $I$ is a finite set. 
In this category the tensor product is defined by 
$$
(V\otimes_I W)_{i,j}=\bigoplus_{k\in I}V_{i,k}\otimes W_{k,j}
$$
and Hom's are given by 
$$
\mathrm{Hom}_{I\times I}(V,W)_{i,j}=\mathrm{Hom}(V_{i,j},W_{i,j})\,.
$$
This monoidal category is of course NOT symmetric at all ... but we will often allow ourselves to use the symmetry isomorphism $\sigma$ of $\texttt{GrMod}_k$ 
in explicit computations as $V\otimes_IW\subset V\otimes W$ and $\mathrm{Hom}_{I\times I}(V,W)\subset\mathrm{Hom}(V,W)$ for any $I\times I$-graded objects $V,W$ in 
$\texttt{GrMod}_k$. 

E.g.~we have the graded counital tensor coalgebra $T_I(V):=\bigoplus_{n\in\mathbb{N}}V^{\otimes_In}$ cogenerated by $V$ as above. But we do not have 
the symmetric algebra in $\texttt{GrMod}^{I\times I}_k$. 

\section{$A_\infty$-categories}\label{s-2}
In the present Section, we introduce the concept of (small) $A_\infty$-categories and related $A_\infty$-functors.
\begin{Def}\label{d-A_inf}
A (small and finite) {\bf $A_{\infty}$-category} is a triple $\mathcal A=(I,A,\mathrm d_A)$, where 
\begin{itemize}
\item $I$ is a finite set (whose elements are called objects); 
\item $A=(A_{{\bf a},{\bf b}})_{({\bf a},{\bf b})\in I\times I}$ is an element in $\texttt{GrMod}_k^{I\times I}$ ($A_{{\bf a},{\bf b}}$ is called the space of morphisms from ${\bf b}$ to ${\bf a}$); 
\item $\mathrm d_A$ is a codifferential on $\mathrm T_I(A[1])$, i.e.~a degree $1$ endomorphism (in $\texttt{GrMod}_k^{I\times I}$) of $\mathrm T_I(A[1])$, satisfying $\Delta\circ \mathrm d_A = (\mathrm d_A\otimes_I 1 + 1\otimes_I \mathrm d_A)\circ \Delta$, $\varepsilon_A\circ \mathrm d_A=0$ and $(\mathrm d_A)^{2}=0$.
\end{itemize}
This is equivalent to require that $(I,\mathrm T(A[1]),\mathrm d_A)$ is a (small) differential graded cocategory. 
\end{Def}
The fact that $\mathrm d_A$ is a coderivation on $\mathrm T_I(A[1])$ and that it lies in the kernel of the counit implies that $\mathrm d_A$ is uniquely determined by its Taylor components 
$\mathrm d_A^n:A[1]^{\otimes_I n}\to A[1]$, $n\geq 0$, {\em via} 
\[
\mathrm d_A\vert_{\mathrm T_I^n(A[1])}=\sum_{m=0}^n\sum_{l=0}^{n-m}1^{\otimes_I l}\otimes_I \mathrm d_A^m\otimes_I1^{\otimes_I (n-m-l)}\,,
\]
where $1^{\otimes_I l}$ denotes the identity on $A[1]^{\otimes_I l}$. 
Then, the condition $(\mathrm d_A)^2=0$ is equivalent to the following infinite set of quadratic equations w.r.t.\ the Taylor components of $\mathrm d_A$:
\begin{equation}\label{eq-A1}
\sum_{i=0}^{k}\sum_{j=1}^{k-i+1}\mathrm d_A^{k-i+1}\circ\big(1^{\otimes_I(j-1)}\otimes_I \mathrm d_A^i\otimes_I 1^{\otimes_I(k+1-j-i)}\big)=0\,,\quad k\geq 0\,.
\end{equation} 

Equivalently, if we consider the maps $\mu_A^n:A^{\otimes_I n}\to A[2-n]$ obtained by twisting appropriately $\mathrm d_A$ w.r.t.~suspension and desuspension, the quadratic relations~\eqref{eq-A1} become
\begin{equation}\label{eq-A2}
\sum_{i=0}^{k}\sum_{j=1}^{k-i+1}(-1)^{i\sum_{l=1}^{j-1}|a_l|+j(i+1)}\mu_A^{k-i+1}(a_1,\dots,a_{j-1},\mu_A^i(a_j,\dots,a_{i+j-1}),a_{i+j},\dots,a_k)=0,\ k\geq 0,
\end{equation}
$a_i\in A_{{\bf a}_{i-1},{\bf a}_i}$, and ${\bf a_0},\dots,{\bf a}_k\in I$. 

An $A_{\infty}$-category $\mathcal A=(I,A,\mathrm d_A)$ is called {\bf flat}, if $\mathrm d_A^0=0$: in this case, $\mathrm d_A^1$ is a differential on $A$, $\mathrm d_A^2$ is an associative product up to homotopy, etc ... 
Otherwise, $\mathcal A$ is called {\bf curved}.
If a flat $A_{\infty}$-category is such that has $\mathrm d_A^k=0$, for $k\geq 3$, then it is called a differential graded (shortly, from now on, DG) category. 
 
We now assume $\mathcal A=(I,A,\mathrm d_A)$ and $\mathcal B=(J,B,\mathrm d_B)$ are two (possibly curved) $A_\infty$-categories in the sense of Definition~\ref{d-A_inf}, then an $A_\infty$-functor from $\mathcal A$ to $\mathcal B$ is the {\em datum} of a functor $\mathcal F$ between the corresponding DG cocategories. 
More precisely, $\mathcal F$ is given by 
\begin{itemize}
\item a map $f:I\to J$; 
\item an $I\times I$-graded coalgebra morphism $F:\mathrm T_I(A[1])\to \mathrm T_I(B[1])$ of degree $0$ which intertwines the codifferentials 
$\mathrm d_A$ and $\mathrm d_B$, i.e.\ $F\circ\mathrm d_A=\mathrm d_B\circ F$. 
\end{itemize}
From the coalgebra (or better, cocategory) structure on $\mathrm T_I(A[1])$ and $\mathrm T_I(B[1])$ (and since $F$ is compatible with the corresponding counits, whence $F_0(1)=1$), it follows immediately that an $A_\infty$-functor from $A$ to $B$ is uniquely specified by its Taylor components $F_n:A[1]^{\otimes_I n}\to B[1]$ via
\[
F\vert_{A[1]^{\otimes_I n}}=\sum_{k=0}^n \sum_{\mu_1,\dots,\mu_k\geq 0\atop \sum_{i=1}^k \mu_i=n}F_{\mu_1}\otimes_I \cdots\otimes_I F_{\mu_k}.
\]
As a consequence, the condition that $F$ intertwines the codifferentials $\mathrm d_A$ and $\mathrm d_B$ can be re-written as an infinite series of equations w.r.t.\ the Taylor components of $\mathrm d_A$, $\mathrm d_B$ and $F$:
\[
\sum_{m=0}^n\sum_{l=0}^{n-m}F_{n-m+1}\circ\left(1^{\otimes_I l}\otimes_I \mathrm d_A^m\otimes_I 1^{\otimes_I (n-m-l)}\right)=\sum_{k=0}^n \mathrm d_B^k\circ\left(\sum_{\mu_1,\dots,\mu_k\geq 0\atop \sum_{i=1}^k \mu_i=n}F_{\mu_1}\otimes_I \cdots\otimes_I F_{\mu_k}\right).
\]
We finally observe that, twisting the Taylor components $F_n$ of an $A_\infty$-morphism $F$ from $A$ to $B$, we get a semi-infinite series of morphisms 
$\phi_n:A^{\otimes_I n}\to B[1-n]$, of degree $1-n$, $n\geq 0$. The natural signs in the previous relations can be computed immediately using suspension and desuspension.

\begin{Exa}
An $A_\infty$-category with only one object is an $A_\infty$-algebra; a DG algebra is a DG category with only one object. 
\end{Exa}

Given an $A_\infty$-category $\mathcal A=(I,A,\mathrm d_A)$ and a subset $J$ of objects, there is an obvious notion of full $A_\infty$-subcategory w.r.t.~$J$. 
In particular, the space of endomorphisms $A_{{\rm a},{\rm a}}$ of a given object ${\rm a}$ is naturally an $A_\infty$-algebra. 

\begin{Exa}\label{ex-2-3}
We consider an $A_\infty$-category $\mathcal C=(I,C,\mathrm d_C)$ with two objects; $I=\{{\bf a},{\bf b}\}$. We further assume $C_{{\bf b},{\bf a}}=0$. Let us define 
$$
A=C_{{\bf a},{\bf a}}\,,\qquad B=C_{{\bf b},{\bf b}}\,,\qquad K=C_{{\bf a},{\bf b}}\,.
$$
$A$ and $B$ are $A_\infty$-algebras, and we say that $K$ is an $A_\infty$-$A$-$B$-bimodule. 
We observe that we can alternatively define an $A_\infty$-$A$-$B$-bimodule structure on $K$ as a codifferential $\mathrm d_K$ on the cofree $(T(A[1]),T(B[1]))$-bicomodule cogenerated by $K[1]$: we write ${\rm d}_K^{m,n}$ for the restriction of the Taylor component $\mathrm d_C^{m+n+1}$ onto the subspace 
$A[1]^{\otimes m}\otimes K[1]\otimes B[1]^{\otimes m}\subset (C[1]^{\otimes_I m+n+1})_{{\bf a},{\bf b}}$ (which takes values in $K[1]=C_{{\bf a},{\bf b}}[1]$). 
We often denote by $\texttt{Cat}_\infty(A,B,K)$ the corresponding $A_\infty$-category. 
\end{Exa}

\begin{Rem}\label{r-A_inf-mod}
We observe that an $A_\infty$-algebra structure $\mathrm d_A$ on $A$ determines an $A_\infty$-$A$-$A$-bimodule structure on $A$ {\em via} the Taylor components
\begin{eqnarray}
\mathrm d_A^{m,n}:=\mathrm d_A^{m+n+1}. \label{selfbimodule}
\end{eqnarray}   
\end{Rem}

\subsection{The Hochschild cochain complex of an $A_\infty$-category}\label{ss-2-1}
We consider an object $A=(A_{{\bf a},{\bf b}})_{{\bf a},{\bf b}\in I\times I}$ of $\texttt{GrMod}_k^{I\times I}$. 
We associate to it another element $\mathrm{C}^\bullet(A,A)$ of $\texttt{GrMod}_k^{I\times I}$, defined as follows: 
$$
\mathrm{C}^\bullet(A,A)=\bigoplus_{p\geq 0}\mathrm{Hom}_{I\times I}(A^{\otimes_I p+1},A)
=\bigoplus_{p\geq 0}\bigoplus_{{\bf a_0},\dots,{\bf a_{p+1}}\in I}\mathrm{Hom}
\left(A_{{\bf a_0},{\bf a_1}}\otimes\cdots\otimes A_{{\bf a_p},{\bf a_{p+1}}},A_{{\bf a_0},{\bf a_{p+1}}}\right)\,,
$$

The $\mathbb Z$-grading on $\mathrm{C}^\bullet(A,A)$ is given as the total grading of the following $\mathbb Z^2$-grading: 
\[
\mathrm{C}^{(p,q)}(A,A)=\mathrm{Hom}^q_{I\times I}\left(A^{\otimes_I p+1},A\right)\,.
\]
We have the standard brace operations on $\mathrm{C}^\bullet(A,A)$: namely, the brace operations are defined {\em via} the usual higher compositions 
(of course, whenever they make sense), i.e.\
\[
\begin{aligned}
&P\{Q_1,\dots,Q_q\}(a_1,\dots,a_n)=\\
&=\sum_{i_1,\dots,i_q}(-1)^{\sum_{k=1}^q \norm{Q_k}\left(i_k-1+\sum_{j=1}^{i_k-1}|a_j|\right)} P(a_1,\dots,Q_1(a_{i_1},\dots),\dots,Q_q(a_{i_q},\dots),\dots,a_n)\,.
\end{aligned}
\]
In the previous sum, $n=p+\sum_{a=1}^q(q_a-1)$, $1\leq i_1$, $i_k+q_k\leq i_{k+1}$, $k=1,\dots,q-1$, $i_q+q_q-1\leq n$, and $a_i$ is a general element of 
$A$, $i=1,\dots,n$; $|Q_k|$ denotes the degree of $Q_k$, while $q_k$ is the number of entries.
We use the standard notation and sign rules, see e.g.~\cite{GJ,GV,TT,CR2} for more details: in particular, $\norm{\bullet}$ denotes the total degree w.r.t.\ the previous bigrading. We finally recall that the graded commutator of the (non-associative) pairing defined by the brace operations on two elements satisfies the requirements for being a graded Lie bracket (w.r.t.\ the total degree), the so-called {\bf Gerstenhaber bracket}. 
\begin{Rem}
Another (more intrinsic) definition of the Hochschild complex is as the space of $I\times I$-graded coderivations of $\mathrm{T}_I(A[1])$: 
$$
\mathrm{CC}(A):=\mathrm{Coder}_{I\times I}\big(\mathrm{T}_I(A[1])\big)=\mathrm{Hom}_{I\times I}\big(\mathrm{T}_I(A[1]),A[1]\big)\,.
$$
In this description the Gerstenhaber bracket becomes more transparent: it is simply the natural Lie bracket of coderivations. 
The identification between $\mathrm{CC}(A)$ and $\mathrm{C}(A,A)$ is again given by an appropriate twisting w.r.t.~suspension and desuspension. 
\end{Rem}
According to the previous remark, the structure of an $A_\infty$-category with $I$ as set of objects and $A$ as $I\times I$-graded space of morphisms 
then translates into the existence of a Maurer--Cartan (shortly, MC) element $\gamma$ in $\mathrm{C}^\bullet(A,A)$, i.e.\ an element $\gamma$ of 
$\mathrm{C}^\bullet(A,A)$ of (total) degree $1$, satisfying $\frac{1}2[\gamma,\gamma]=\gamma\{\gamma\}=0$.
Finally, the MC element $\gamma$ specifies a degree $1$-differential $\mathrm d_\gamma=[\gamma,\bullet]$, where $[\bullet,\bullet]$ denotes the Gerstenhaber bracket on $\mathrm{C}^\bullet(A,A)$. We obtain this way a DG Lie algebra. 

\begin{Rem}\label{rem-infinite}
The $A_\infty$-structures we have considered so far only have a finite number of Taylor components. 
In order to define $A_\infty$-structures in full generality, where an infinite number of Taylor components is allowed, 
one has to work in the category $\widehat{\texttt{GrMod}_k}$ (see Remark \ref{rem-completion}). More precisely, in this context the degree $n$ part of the completed Hochschild complex is 
$$
\widehat{\mathrm{C}}^n(A,A):=\prod_{p+q=n}\mathrm{Hom}^q_{I\times I}(A^{\otimes_I p+1},A)
=\prod_{p+q=n}\left(\bigoplus_{{\bf a_0},\dots,{\bf a_{p+1}}\in I}\mathrm{Hom}^q
\left(A_{{\bf a_0},{\bf a_1}}\otimes\cdots\otimes A_{{\bf a_p},{\bf a_{p+1}}},A_{{\bf a_0},{\bf a_{p+1}}}\right)\right)\,,
$$
\end{Rem}

\begin{Exa}\label{ex-2-6}
We now make more explicit the case of the $A_\infty$-category $\texttt{Cat}_\infty(A,B,K)$ of Example~\ref{ex-2-3}. 
First of all, the bigrading on $C=\texttt{Cat}_\infty(A,B,K)$ can be read immediately from the above conventions, i.e.
\[
\mathrm{C}^n(C,C)=\bigoplus_{p+q=n}\mathrm{Hom}^q(A^{\otimes(p+1)},A)\oplus\bigoplus_{p+q+r=n}\mathrm{Hom}^r(A^{\otimes p}\otimes K\otimes B^{\otimes q},K)\oplus\bigoplus_{p+q=n}\mathrm{Hom}^q(B^{\otimes(p+1)},B)\,.
\]
The $A_\infty$-structure on $\texttt{Cat}_\infty(A,B,K)$ specifies a MC element $\gamma$, which splits into three pieces according to $\gamma=\mathrm d_A+\mathrm d_K+\mathrm d_B$.
By the very construction of the Hochschild differential $\mathrm d_\gamma$, $\mathrm d_\gamma$ splits into five components, since, for $\varphi=\varphi_A+\varphi_K+\varphi_B$ a general element of $\mathrm{C}^\bullet(C,C)$,
\[
\mathrm d_\gamma \varphi=[\mathrm d_A,\varphi_A]+\mathrm d_K\{\varphi_A\}+[\gamma,\varphi_K]+\mathrm d_K\{\varphi_A\}+[\mathrm d_B,\varphi_B]\,.
\]
We observe that $[\gamma,\varphi_K]=[\mathrm d_K,\varphi_K]-(-1)^{\norm{\varphi_K}}\varphi_K\{\mathrm d_A+\mathrm d_K+\mathrm d_B\}$; we denote by $\mathrm{C}^\bullet(A,K,B)$ the subcomplex which consists of elements $\varphi_K$ in the middle term of the previous splitting. 

We want to explain the meaning of the five components in the alternative description of the Hochschild complex. 
An element $\phi$ in  
$\mathrm{CC}^n(C)$ consists of a triple $(\phi_A,\phi_K,\phi_B)$, where $\phi_A$ (resp.~$\phi_B$) is a coderivation of ${\rm T}(A[1])$ 
(resp.~${\rm T}(B[1])$) and $\phi_K$ is a coderivation of the bicomodule ${\rm T}(A[1])\otimes K[1]\otimes {\rm T}(B[1])$ w.r.t.~$\phi_A$ and $\phi_B$. 
Now the MC element $\gamma$ gives such an element $({\rm d}_A,{\rm d}_K,{\rm d}_B)$, which moreover squares to zero. 
The five components can be then interpreted as 
\[
\mathrm d_\gamma \phi=[\mathrm d_A,\phi_A]+\mathrm L_A\circ \phi_A+[\mathrm d_K,\phi_K]+\mathrm R_B\circ\phi_B+[\mathrm d_B,\phi_B]\,.
\]
The meaning of the morphisms $\mathrm L_A$ and $\mathrm R_B$, the derived left- and right-action, is explained in full details below in Subsection \ref{s-3}.
\end{Exa}

\subsubsection*{Signs considerations}

We now want to discuss the signs appearing in the brace operations, which correspond to the natural Koszul signs appearing when one considers all possible higher compositions between different elements of $\mathrm{Hom}(\mathrm T(A[1]),A[1])$.

Before entering into the details, we need to be more precise on grading conventions: if $\phi$ is a general element of $\mathrm{Hom}^{m}(B[1]^{\otimes n},B[1])$, then we write $|\phi|=m$, and similar notation holds, when $\phi$ is an element of $\mathrm{Hom}^{r}(B[1]^{\otimes p}\otimes K[1]\otimes A[1]^{\otimes q},K[1])$.
On the other hand, we write $\norm{\phi}=m+n-1$, and similarly, if $\phi$ is in $\mathrm{Hom}^{r}(B[1]^{\otimes p}\otimes K[1]\otimes A[1]^{\otimes q},K[1])$, $\norm\phi=p+q+r$.

We consider e.g.\ the Gerstenhaber bracket on $B$: for $\phi_i$, $i=1,2$, in $\mathrm{Hom}^{m_i}(B[1]^{\otimes n_i},B[1])$, we have
\[
[\phi_1,\phi_2 ]:=\sum_{j=1}^{n_1} \phi_1\circ (1^{\otimes(j-1)}\otimes \phi_2\otimes 1^{\otimes(n_1-j)})-(-1)^{|\phi_1||\phi_2|}(\phi_2\leftrightarrow\phi_1),
\] 
Twisting w.r.t.\ suspension and desuspension (we recall that the suspension $s:B\rightarrow B[1]$ has degree -1 and the desuspension $s^{-1}$ degree 1), we introduce the desuspended maps $\widetilde{\phi}_i\in \mathrm{Hom}^{1+m_i-n_i}(B^{\otimes n_i},B)$, and we then set $|\widetilde{\phi}_i|:=1+m_i-n_i$ and $\norm{\widetilde{\phi}_i}=m_1$; in other words, $\phi_i=s\circ \widetilde{\phi}_i\circ (s^{-1})^{\otimes n_i},\ i=1,2$.

We observe that 
\[
\norm{\widetilde{\phi}_i}=|\phi_i|=m_i,\ |\tilde{\phi}_i|=\norm{\phi_i}=m_i+n_i-1\ \text{modulo $2$}.
\]
We then get, by explicit computations,
\[
[\widetilde \phi_1,\widetilde \phi_2]=\widetilde{\phi}_1\bullet \widetilde{\phi}_2-(-1)^{\norm{\widetilde{\phi}_1}\norm{\widetilde{\phi}_2}}\widetilde{\phi}_2\bullet\widetilde{\phi}_2,
\]
where the new desuspended signs for the higher composition $\bullet$ are given by
\begin{equation}\label{eq-des-sign}
\widetilde\phi_1\bullet\widetilde\phi_2=\sum_{j=1}^{n_1}(-1)^{\left(|\widetilde\phi_2|+n_2-1\right)(n_1-1)+(j-1)(n_2-1)}\widetilde\phi_1\circ \left(1^{\otimes (j-1)}\otimes\widetilde\phi_2\otimes 1^{\otimes(n_1-j)}\right).
\end{equation}
We observe that these signs appear also in~\cite{CF}.
Obviously, replacing $B$ by $A$, we repeat all previous arguments to come to the signs for the Gerstenhaber bracket on $\mathrm{C}^\bullet(A,A)$.

Further, assuming e.g.\ $\phi_i$, $i=1,2$, is a general element of $\mathrm{Hom}^{r_i}(B[1]^{\otimes p_i}\otimes K[1]\otimes A[1]^{\otimes q_i},K[1])$, we introduce the desuspended map {\em via} $\phi_i=s\circ \widetilde{\phi}_i \circ (s^{-1})^{\otimes p_i+q_i+1}$, which is an element of $\mathrm{Hom}^{r_i-p_i-q_i}(B^{\otimes p_i}\otimes K\otimes A^{\otimes q_i},K)$.

Setting $|\widetilde\phi_i|=r_i-p_i-q_i$ and $\norm{\widetilde\phi_i}=r_i$, we have
\[
\norm{\widetilde{\phi}_i}=|\phi_i|=r_i,\ |\widetilde{\phi}_i|=\norm{\phi_i}=r_i+p_i+r_i\ \text{modulo $2$}.
\]
We further get the higher composition $\bullet$ between $\widetilde\phi_1$ and $\widetilde\phi_2$, coming from the natural brace operations, with corresponding signs
\[
\widetilde\phi_1\bullet\widetilde\phi_2=(-1)^{\left(|\widetilde\phi_2|+p_2+q_2\right)(p_1+q_1)+p_1(p_2+q_2)}\widetilde{\phi}_1\circ \left(1^{\otimes p_1}\otimes \widetilde\phi_2\otimes 1^{\otimes q_1}\right)
\]
If now $\phi_1$ is in $\mathrm{Hom}^{r_1}(B[1]^{\otimes p_1}\otimes K[1]\otimes A[1]^{\otimes q_1},K[1])$ and $\phi_2$ is in $\mathrm{Hom}^{m_2}(B[1]^{\otimes n_2},B[1])$, and by introducing the desuspended maps $\widetilde\phi_i$, $i=1,2$, whose (total) degrees satisfy the same relations as above, we get the higher composition with corresponding signs between $\widetilde\phi_1$ and $\widetilde\phi_2$, coming from the previously described brace operations:
\[
\widetilde\phi_1\bullet\widetilde\phi_2=\sum_{j=1}^{p_1}(-1)^{\left(|\widetilde\phi_2|+n_2-1\right)(p_1+q_1)+(j-1)(n_2-1)}\widetilde\phi_2\circ\left(1^{\otimes(j-1)}\otimes\widetilde\phi_2\otimes 1^{\otimes (p_1+q_1+1-j)}\right).
\]
Finally, if $\phi_1$, resp.\ $\phi_2$, lies in $\mathrm{Hom}^{r_1}(B[1]^{\otimes p_1}\otimes K[1]\otimes A[1]^{\otimes q_1},K[1])$, resp.\ $\mathrm{Hom}^{m_2}(A[1]^{\otimes n_1},A[1])$, then the higher composition between the desuspended maps $\widetilde\phi_1$ and $\widetilde\phi_2$ with corresponding signs, coming from the brace operations, has the explicit form
\[
\widetilde\phi_1\bullet\widetilde\phi_2=\sum_{j=1}^{q_1}(-1)^{\left(|\widetilde\phi_2|+n_2-1\right)(p_1+q_1)+(p_1+j)(n_2-1)}\widetilde\phi_2\circ\left(1^{\otimes(p_1+j)}\otimes\widetilde\phi_2\otimes 1^{\otimes (q_1-j)}\right).
\]

\section{Keller's condition in the $A_\infty$-framework}\label{s-3}
We now discuss some cohomological features of the Hochschild cochain complex of the $A_\infty$-category $\texttt{Cat}_\infty(A,B,K)$ from Example~\ref{ex-2-3}, Section~\ref{s-2}: in particular, we will extend to this framework the classical result of Keller for DG categories~\cite{Keller}, which is a central piece in the proof of the main result of~\cite{Sh}.
\begin{Rem}\label{r-hom}
Unless otherwise specified, $\mathrm{Hom}$ and $\mathrm{End}$ have to be understood in the category $\texttt{GrMod}_k$. 
\end{Rem} 

\subsection{The derived left- and right-actions}\label{ss-3-1}
We consider $A$, $B$ and $K$ as in Example~\ref{ex-2-3}, Section~\ref{s-2}, borrowing the same notation.

We consider the restriction $\mathrm d_{K,B}$ of $\mathrm d_K$ to $K[1]\otimes \mathrm T(B[1])$, i.e.\ the map 
\[
\mathrm d_{K,B}=\mathrm P_{K,B}\circ \mathrm d_K,
\]
where $\mathrm P_{K,B}$ denotes the natural projection from $\mathrm T(A[1])\otimes K[1]\otimes\mathrm T(B[1])$ onto $K[1]\otimes \mathrm T(B[1])$.

A direct check implies that $\mathrm P_{K,B}$ is a morphism of right $\mathrm T(B[1])$-comodules, whence it follows directly that $\mathrm d_{K,B}$ is a coderivation on $K[1]\otimes\mathrm T(B[1])$.
\begin{Rem}\label{r-left}
Similarly, the restriction of $\mathrm d_K$ on $\mathrm T(A[1])\otimes K[1]$ defines a left coderivation $\mathrm d_{A,K}$ on $\mathrm T(A[1])\otimes K[1]$.
\end{Rem}
For $A$, $B$ and $K$ as above, we set 
\[
\underline{\mathrm{End}}_{-B}(K)=\mathrm{End}^c_{\mathrm{comod}-\mathrm T(B[1])}(K[1]\otimes \mathrm T(B[1])),
\]
where the superscript $c$ means that we consider coalgebra endomorphisms of the right cofree $\mathrm T(B[1])$-comodule $K[1]\otimes \mathrm T(B[1])$, for which only finitely many Taylor components are non-trivial.
Obviously, $\underline{\mathrm{End}}_{-B}(K)$ becomes, w.r.t.\ the composition, a graded algebra (shortly, GA).

Further, there is an obvious identification
\[
\underline{\mathrm{End}}_{-B}(K)=\bigoplus_{p\in\mathbb Z\atop q\geq 0}\mathrm{Hom}^p(K[1]\otimes \mathrm B[1]^{\otimes q},K[1])=\bigoplus_{p\in\mathbb Z\atop q\geq 0} \mathrm{Hom}^{p-q}(K\otimes B^{\otimes q},K),
\]
in the category $\texttt{GrMod}_k$.
We will sometimes refer to $p$, resp.\ $q$, as to the total, resp.\ cohomological, degree: their difference $p-q$ is the internal grading. 

The derived left action of $A$ on $K$, denoted by $\mathrm L_A$, is defined as a coalgebra morphism from $\mathrm T(A[1])$ to $\mathrm T(\underline{\mathrm{End}}_{-B}(K)[1])$, both endowed with the obvious coalgebra structures, whose $m$-th Taylor component, viewed as an element of $\underline{\mathrm{End}}_{-B}(K)[1]$, decomposes as 
\begin{equation}\label{eq-left-comp}
\mathrm L_{A}^m(a_1|\cdots|a_m)^n(k|b_1|\cdots|b_n)=\mathrm d_K^{m,n}(a_1|\cdots|a_m|k|b_1|\cdots|b_n),\ m\geq 1,\ n\geq 0. 
\end{equation}
In a more formal way, the Taylor component $\mathrm L_A^m$ may be defined as 
\[
\mathrm L_A^m(a_1|\cdots|a_m)=(\mathrm P_{K,B}\circ \mathrm d_K)(a_1|\cdots|a_m|\cdots). 
\]
It is not difficult to check that $\mathrm L_A^m(a_1|\cdots|a_m)$ is an element of $\underline{\mathrm{End}}_{-B}(K)$.

The grading conditions on $\mathrm d_K$ imply, by direct computations, that $\mathrm L_A^m$ is a morphism from $A[1]^{\otimes n}$ to $\underline{\mathrm{End}}_{-B}(K)[1]$ of degree $0$.

For later computations, we write down explicitly the Taylor series of the derived left action up to order $2$, namely,
\[
\mathrm L_A(a_1|\cdots|a_n)=\mathrm L_A^n(a_1|\cdots|a_n)+\sum_{n_1+n_2=n\atop n_i\geq 1,\ i=1,2}\left(\mathrm L_A^{n_1}(a_1|\cdots|a_{n_1})|\mathrm L_A^{n_2}(a_{n_1+1}|\cdots|a_n)\right)+\cdots
\] 

We now want to discuss an $A_\infty$-algebra structure on $\underline{\mathrm{End}}_{-B}(K)$.
For this purpose, we first consider $\mathrm d_{K,B}^2$: since $\mathrm d_{K,B}$ is a right coderivation on $K[1]\otimes \mathrm T(B[1])$, its square is easily verified to be an element of $\underline{\mathrm{End}}_{-B}(K)$.
\begin{Lem}\label{l-curv-end}
The operator $\mathrm d_{K,B}^2$ satisfies 
\[
\mathrm d_{K,B}^2=-\mathrm L_A^1(\mathrm d_A^0(1)).
\]
\end{Lem}
\begin{proof}
By its very definition, $\mathrm d_{K,B}$ obeys
\[
\mathrm d_{K,B}^2=\mathrm P_{K,B}\circ \mathrm d_K\circ \mathrm P_{K,B}\circ \mathrm d_K\big\vert_{K[1]\otimes \mathrm T(B[1])}.
\]
Since $\mathrm d_K$ is a bicomodule morphism, then, taking into account the definition of the left and right coactions $\Delta_L$ and $\Delta_R$ on $\mathrm T(A[1])\otimes K[1]\otimes \mathrm T(B[1])$, we get
\[
\mathrm P_{K,B}\circ \mathrm d_K\big\vert_{K[1]\otimes\mathrm T(B[1])}=\mathrm d_K\big\vert_{K[1]\otimes \mathrm T(B[1])}-\left(\mathrm d_A^0(1)|\bullet\right).
\]
Since $\mathrm d_K^2=0$, the claim follows directly.
\end{proof}

Therefore, $\underline{\mathrm{End}}_{-B}(K)$ inherits a structure of $A_\infty$-algebra, i.e.\ there is a degree $1$ codifferential $Q$, whose only non-trivial Taylor components are 
\[
Q^0(1)=\mathrm L_A^1(\mathrm d_A^0(1)),\ Q^1(\varphi)=-\left[\mathrm d_{K,B},\varphi\right],\ Q^2(\varphi_1|\varphi_2)=(-1)^{|\varphi_1|}\varphi_1\circ\varphi_2.
\]
\begin{Rem}\label{r-compl-end}
If we consider an $A_\infty$-$A$-$B$-bimodule structure on $K$ with finitely many non-trivial Taylor components, then the derivation $Q^1$ is well-defined on $\underline{\mathrm{End}}_{-B}(K)$.
When considering a more general $A_\infty$-$A$-$B$-bimodule structure, see Remark~\ref{rem-infinite}, Subsection~\ref{ss-2-1}, then we have to consider a completed version $\widehat{\underline{\mathrm{End}}}_{-B}(K)$ of $\underline{\mathrm{End}}_{-B}(K)$, allowing comodule morphisms with infinitely many Taylor components (this translates into switching from the category $\texttt{GrMod}_k$ to $\widehat{\texttt{GrMod}}_k$), in order to make $Q^1$ well-defined.
\end{Rem}
\begin{Rem}\label{r-left-mod}
In a similar way, we may introduce the $A_\infty$-algebra $\underline{\mathrm{End}}_{A-}(K)=\mathrm{End}_{\mathrm T(A[1])-\mathrm{comod}}(\mathrm T(A[1])\otimes K[1])$ and the derived right action $\mathrm R_B$: accordingly, $\underline{\mathrm{End}}_{A-}(K)$ is an $A_\infty$-algebra, with $A_\infty$-structure given by the curvature $Q^0(1)=\mathrm R_B(\mathrm d_B^0(1))$, degree $1$ derivation $\left[\mathrm d_{A,K},\bullet\right]$, and composition as product.
Needless to mention, Remark~\ref{r-compl-end} has to be taken into accout, with due modifications, also for $\underline{\mathrm{End}}_{A-}(K)$.
\end{Rem}
It is clear that, if $A$ and $B$ are flat $A_\infty$-algebras, $\underline{\mathrm{End}}_{-B}(K)$ and $\underline{\mathrm{End}}_{A-}(K)$ are DG algebras.
\begin{Rem}\label{r-end}
The DG algebras $\underline{\mathrm{End}}_{-B}(K)$ and $\underline{\mathrm{End}}_{A-}(K)$ have been introduced by B.~Keller in~\cite{Kel-A}.
\end{Rem}
\begin{Lem}\label{l-left}
The derived left action $\mathrm L_A$ is an $A_\infty$-morphism from $A$ to $\underline{\mathrm{End}}_{-B}(K)$
\end{Lem}
\begin{proof}
The condition for $\mathrm L_A$ to be an $A_\infty$-morphism can be checked by means of its Taylor components of $\mathrm L_A$: namely, recalling that the $A_\infty$-structure on $\underline{\mathrm{End}}_{-B}(K)$ has only three non-trivial components, we have to check the two identities
\begin{equation}\label{eq-left-cond}
\begin{aligned}
&(\mathrm L_A\circ \mathrm d_A)(1)=(Q\circ \mathrm L_A)(1),\\
&\sum_{k=0}^m \sum_{i=1}^{m-k+1}(-1)^{\sum_{j=1}^{i-1}(|a_j|-1)}\mathrm L_A^{m-k+1}\left(a_1|\cdots|\mathrm d_A^k(a_i|\cdots|a_{i+k-1})|a_{i+k}|\cdots|a_m\right)=\\
&=-\left[\mathrm d_{K,B},\mathrm L_A^m(a_1|\cdots|a_m)\right]+\sum_{m_1+m_2=m\atop m_i\geq 1,\ i=1,2}(-1)^{\sum_{k=1}^{m_1}(|a_k|-1)}\mathrm L_A^{m_1}(a_1|\cdots|a_{m_1})\circ \mathrm L_A^{m_2}(a_{m_1+1}|\cdots|a_m). 
\end{aligned}
\end{equation}
The first identity in~\eqref{eq-left-cond} follows immediately from the construction of the $A_\infty$-structure on $\underline{\mathrm{End}}_{-B}(K)$.
In order to prove the second one, we evaluate both sides of the second expression explicitly on a general element of $K[1]\otimes \mathrm T(B[1])$, projecting down to $K[1]$: writing down the natural signs arising from Koszul's sign rule and the differential $[\mathrm d_{K,B},\bullet]$, we see immediately that it is equivalent to the condition that $K$ is an $A_\infty$-$A$-$B$-bimodule. 
\end{proof}
Of course, similar arguments imply that there is an $A_\infty$-morphism $\mathrm R_B$ from $B$ to $\underline{\mathrm{End}}_{A-}(K)^{\mathrm{op}}$, where the suffix ``op'' refers to the fact that we consider the opposite product on $\underline{\mathrm{End}}_{A-}(K)$: again, the condition that $\mathrm R_B$ is an $A_\infty$-morphism is equivalent to the fact that $K$ is an $A_\infty$-$A$-$B$-bimodule.

Furthermore, $\mathrm L_A$, resp.\ $\mathrm R_B$, endow $\underline{\mathrm{End}}_{-B}(K)$, resp.\ $\underline{\mathrm{End}}_{A-}(K)^\mathrm{op}$, with a structure of $A_\infty$-$A$-$A$-bimodule, resp.\ -$B$-$B$-bimodule.

In a more conceptual way, given two $A_\infty$-algebras $A$ and $B$ and an $A_\infty$-morphism $F$ from $A$ to $B$, we first view both $A$ and $B$ as $A_\infty$-bimodules in the sense of Remark~\ref{r-A_inf-mod}, Section~\ref{s-2}.
Then, we define an $A_\infty$-$A$-$A$-bimodule structure on $B$ simply {\em via} the codifferential $\mathrm d_B\circ (F\otimes 1\otimes F)$, where $\mathrm d_B$ denotes improperly the codifferential inducing the $A_\infty$-$B$-$B$-bimodule structure on $B$.

Explicitly, we write down the Taylor components of the $A_\infty$-$A$-$A$-bimodule structure on $\underline{\mathrm{End}}_{-B}(K)$: since the $A_\infty$-structure on $\underline{\mathrm{End}}_{-B}(K)$ has only three non-trivial Taylor components, a direct computation shows 
\begin{equation}\label{eq-end_A-inf}
\begin{aligned}
\mathrm Q^{0,0}(\varphi)&=-[\mathrm d_{K,B},\varphi],\ & \ \mathrm Q^{m,n}&=0,\ n,m\geq 1,\\
\mathrm Q^{m,0}(a_1|\cdots|a_m|\varphi)&=(-1)^{\sum_{k=1}^m(|a_k|-1)}\mathrm L_A^m(a_1|\cdots|a_m)\circ\varphi,\ & \ \mathrm Q^{0,n}(\varphi|a_1|\cdots |a_n)&=(-1)^{\varphi} \varphi\circ\mathrm L_A^m(a_1|\cdots|a_n),
\end{aligned}
\end{equation}
where $m$, resp.\ $n$, is bigger than $1$ in the third, resp.\ fourth, formula. 
Similar formul\ae\ hold true for the derived right action.

\subsection{The Hochschild cochain complex of an $A_\infty$-algebra}\label{ss-3-2}
For the $A_\infty$-algebra $A$, we consider its Hochschild cochain complex with values in itself: as we have already seen in Subsection~\ref{ss-2-1}, it is defined as 
\[
\mathrm{C}^\bullet(A,A)=\mathrm{Coder}(\mathrm T(A[1]))=\mathrm{Hom}(\mathrm T(A[1]),A[1]),
\]
the vector space of coderivations of the coalgebra $\mathrm T(A[1])$ (with the obvious coalgebra structure), and differential $[\mathrm d_A,\bullet]$.

If we now consider a general $A_\infty$-$A$-$A$-bimodule $M$, we define the Hochschild cochain complex of $A$ with values in $M$, which we denote by $\mathrm{C}^\bullet(A,M)$, as the vector space of morphisms $\varphi$ from $\mathrm T(A[1])$ to the bicomodule $\mathrm T(A[1])\otimes M[1]\otimes \mathrm T(A[1])$, for which
\[
\Delta_L\circ\varphi=(1\otimes \varphi)\circ\Delta_A,\ \Delta_R\circ\varphi=(\varphi\otimes 1)\circ\Delta_A.
\]
The differential is then simply given by $\mathrm d_M\varphi=\mathrm d_M\circ\varphi-(-1)^{|\varphi|}\varphi\circ\mathrm d_A$.
It is clear that $\mathrm{C}^\bullet(A,M)=\mathrm{Hom}(\mathrm T(A[1]),M[1])$.
\begin{Rem}\label{r-hoch}
The previous definition of the Hochschild cochain complex $\mathrm C^\bullet(A,M)$, in the case where $M=A$, agrees with the definition of $\mathrm C^\bullet(A,A)$: this is because, in both cases, $C^\bullet(A,A)=\mathrm{Hom}(\mathrm T(A[1]),A[1])$, and because $A$ becomes an $A_\infty$-$A$-$A$-bimodule in the sense of Remark~\ref{r-A_inf-mod}, Section~\ref{s-2}, which implies that the differentials on the two complexes coincide.
\end{Rem}
We further consider the complex $\mathrm{C}^\bullet(A,B,K)$, with differential $[\mathrm d_K,\bullet]$ as in Subsection~\ref{ss-2-1}.

Finally, for $A$, $B$ and $K$ as above, we consider the $A_\infty$-$A$-$A$-bimodule $\underline{\mathrm{End}}_{-B}(K)$; similar arguments work for the $A_\infty$-$B$-$B$-bimodule $\underline{\mathrm{End}}_{A-}(K)^\mathrm{op}$.

\begin{Lem}\label{l-hoch-mix}
The complexes $\left(\mathrm{C}^\bullet(A,B,K),[\mathrm d_K,\bullet]\right)$ and $\left(\mathrm{C}^\bullet(A,\underline{\mathrm{End}}_{-B}(K)),\mathrm d_{\underline{\mathrm{End}}_{-B}(K)}\right)$ are isomorphic.
\end{Lem}
\begin{proof}
It suffices to give an explicit formula for the isomorphism: a general element $\varphi$ of $\mathrm{C}^\bullet(A,B,K)$ is uniquely determined by its Taylor components $\varphi^{m,n}$ from $A[1]^{\otimes m}\otimes K[1]\otimes B[1]^{\otimes n}$ to $K[1]$.

On the other hand, a general element $\psi$ of $\mathrm{C}^\bullet(A,\underline{\mathrm{End}}_{-B}(K))$ is also uniquely determined by its Taylor components $\psi^m$ from $A[1]^{\otimes m}$ to $\underline{\mathrm{End}}_{-B}(K)$; in turn, any Taylor component $\psi^m(a_1|\cdots|a_m)$ is, by definition, completely determined by its Taylor components $(\psi^m(a_1|\cdots|a_m))^{n}$ from $K[1]\otimes B[1]^{\otimes n}$ to $K[1]$.

Thus, the isomorphism from $\mathrm{C}^\bullet(A,B,K)$ to $\mathrm{C}^\bullet(A,\underline{\mathrm{End}}_{-B}(K))$ is explicitly described {\em via}
\[
\left(\widetilde{\varphi}^m(a_1|\cdots|a_m)\right)^n(k|b_1|\cdots|b_n)=\varphi^{m,n}(a_1|\cdots|a_m|k|b_1|\cdots|b_n),\ m,n\geq 0.
\]
It remains to prove that the previous isomorphism is a chain map: for the sake of simplicity, we omit the signs here, since they can be all deduced quite easily from our previous conventions and from Koszul's sign rule, and we only write down the formul\ae, from which we deduce immediately the claim.
It also suffices, by construction, to prove the claim on the corresponding Taylor components.

Thus, we consider 
\[
\begin{aligned}
&\left(\widetilde{[\mathrm d_K,\varphi]}^m(a_1|\cdots|a_m)\right)^{n}(k|b_1|\cdots|b_n)=\left([\mathrm d_K,\varphi]\right)^{m,n}(a_1|\cdots|a_m|k|b_1|\cdots|b_n)=\\
=&(\mathrm d_K\circ\varphi)^{m,n}(a_1|\cdots|a_m|k|b_1|\cdots|b_n)-(-1)^{|\varphi|}(\varphi\circ\mathrm d_K)^{m,n}(a_1|\cdots|a_m|k|b_1|\cdots|b_n).
\end{aligned}
\]
The first term in the last expression can be re-written as a sum of terms of the form
\begin{equation}\label{eq-hoch-1}
\mathrm d_K^{i-1,n-j}\left(a_1|\cdots|a_{i-1}|\varphi^{(m-i+1,j)}(a_i|\cdots|a_m|k|b_1|\cdots|b_j)|b_{j+1}|\cdots|b_n\right), 1\leq i\leq m+1,\ 0\leq j\leq n.
\end{equation}
On the other hand, the second term in the last expression is the sum of three types of terms, which are listed here:
\begin{align}
\label{eq-hoch-2}&\mathrm \varphi^{i-1,n-j}\left(a_1|\cdots|a_{i-1}|\mathrm d_K^{(m-i+1,j)}(a_i|\cdots|a_m|k|b_1|\cdots|b_j)|b_{j+1}|\cdots|b_n\right),\ 1\leq i\leq m+1,\ 0\leq j\leq n,\\
\label{eq-hoch-3}&\mathrm \varphi^{m-j+1,n}\left(a_1|\cdots|a_{i-1}|\mathrm d_B^j(a_i|\cdots|a_{i+j-1})|a_{i+j}|\cdots|a_m|k|b_1|\cdots|b_n\right),\ 0\leq i\leq m+1,\ 0\leq j\leq m,\\
\label{eq-hoch-4}&\mathrm \varphi^{m,n-j+1}\left(a_1|\cdots|a_m|k|b_1|\cdots|b_{i-1}|\mathrm d_A^j(b_i|\cdots|b_{i+j-1})|b_{i+j}|\cdots|b_n\right),\ 0\leq i\leq n+1,\ 0\leq j\leq n.
\end{align}

We now consider the expression 
\[
\left(\mathrm d_{\underline{\mathrm{End}}_{-B}(K)}\widetilde\varphi\right)^m(a_1|\cdots|a_m)=\left(\mathrm d_{\underline{\mathrm{End}}_{-B}(K)}\circ\widetilde\varphi\right)^m(a_1|\cdots|a_m)-(-1)^{|\widetilde\varphi|}\left(\widetilde\varphi\circ \mathrm d_A\right)^m(a_1|\cdots|a_m).
\]
If we further consider the previous identity applied to an element $(k|b_1|\cdots|b_n)$ as above, then the second term on the right-hand side is, by definition, a sum of terms of the type~\eqref{eq-hoch-3}.

On the other hand, we consider the first term on the right-hand side of the previous expression: we recall the Taylor components~\eqref{eq-end_A-inf} of the $A_\infty$-$A$-$A$-bimodule structure on $\underline{\mathrm{End}}_{-B}(K)$, whence 
\begin{equation}\label{eq-d-end}
\begin{aligned}
\left(\mathrm d_{\underline{\mathrm{End}}_{-B}(K)}\circ\widetilde\varphi\right)^m(a_1|\cdots|a_m)&=-\left[\mathrm d_{K,B},\widetilde\varphi^m(a_1|\cdots|a_m)\right]+\\
&\phantom{=}+\sum_{m_1+m_2=m\atop m_i\geq 1,\ i=1,2}(-1)^{|\widetilde\varphi|+\sum_{k=1}^{m_1}(|a_k|-1)}\widetilde\varphi^{m_1}(a_1|\cdots|a_{m_1})\circ\mathrm L_A^{m_2}(a_{m_1+1}|\cdots|a_m)+\\
&\phantom{=}+\sum_{m_1+m_2=m\atop m_i\geq 1,\ i=1,2}(-1)^{(|\widetilde\varphi|+1)\left(\sum_{k=1}^{m_1}(|a_k|-1)\right)}\mathrm L_A^{m_1}(a_{1}|\cdots|a_{m_1})\circ\widetilde\varphi^{m_2}(a_{m_1+1}|\dots|a_m).
\end{aligned}
\end{equation}

The sum of expressions~\eqref{eq-hoch-2}, for which $i=m+1$, and~\eqref{eq-hoch-4}, equals, by definition, the first term on the right-hand side of~\eqref{eq-d-end}; expressions~\eqref{eq-hoch-2}, resp.~\eqref{eq-hoch-1}, for which $i\leq m$, sum up to the second, resp.\ third, term on the right-hand side of~\eqref{eq-d-end}.
\end{proof}
The same arguments, with obvious due changes, imply that the complex $\left(\mathrm{C}^\bullet(A,B,K),[\mathrm d_K,\bullet]\right)$ is isomorphic to the Hochschild chain complex $\left(\mathrm{C}^\bullet(B,\underline{\mathrm{End}}_{A-}(K)^\mathrm{op}),\mathrm d_{\underline{\mathrm{End}}_{A-}(K)^\mathrm{op}}\right)$, replacing $\mathrm L_A$ by $\mathrm R_B$.

Finally, composition with $\mathrm L_A$ and $\mathrm R_B$ defines morphisms of complexes
\[
\begin{aligned}
&\mathrm L_A:\mathrm{C}^\bullet(A,A)\to \mathrm{C}^\bullet(A,B,K)\cong\mathrm{C}^\bullet(A,\underline{\mathrm{End}}_{-B}(K)),\\
&\mathrm R_B:\mathrm{C}^\bullet(B,B)\to \mathrm{C}^\bullet(A,B,K)\cong \mathrm{C}^\bullet(B,\underline{\mathrm{End}}_{A-}(K)^\mathrm{op}).
\end{aligned}
\] 
More precisely, composition with $\mathrm L_A$ on $\mathrm C^\bullet(A,A)$ is defined {\em via} the assignment
\[
\begin{aligned}
&\left(\mathrm L_A\circ\varphi\right)^{m,n}(a_1|\cdots|a_m|k|b_1|\cdots|b_n)=\\
&=\sum_{i=0}^m\sum_{j=0}^{m+1}(-1)^{|\varphi|\left(\sum_{k=1}^{j-1}(|a_k|-1)\right)}\mathrm d_K^{m-i+1,n}(a_1|\cdots|\varphi^i(a_j|\cdots|a_{j+i-1})|\cdots|a_m|k|b_1|\cdots|b_n),
\end{aligned}
\]
and a similar formula defines composition with $\mathrm R_B$.
The fact that composition with $\mathrm L_A$ and $\mathrm R_A$ is a map of complexes is a direct consequence of the computations in the proof of Lemma~\ref{l-left}, Subsection~\ref{ss-3-1}, and of Lemma~\ref{l-hoch-mix}.
\begin{Rem}\label{r-l-comp}
We observe that the previous formula coincides with $\mathrm d_K\{\varphi\}$, using the notation of Subsection~\ref{ss-2-1}.
\end{Rem}
\begin{Rem}\label{r-compl-hoch}
There are obvious changes to be made when switching from the category $\texttt{GrMod}_k$ to the category $\widehat{\texttt{GrMod}}_k$; all results in this Subsection can be translated almost {\em verbatim} to $\widehat{\texttt{GrMod}}_k$, with obvious due modifications.
\end{Rem}

\subsection{Keller's condition}\label{ss-3-3}
From the arguments of Subsection~\ref{ss-2-1}, it is easy to verify that the natural projections $\mathrm p_A$ and $\mathrm p_B$ from $\mathrm{C}^\bullet(\texttt{Cat}_\infty(A,B,K))$ onto $\mathrm{C}^\bullet(A,A)$ and $\mathrm{C}^\bullet(B,B)$, respectively, are well-defined morphisms of complexes.

A natural question for our purposes is the following one: under which conditions are the projections $\mathrm p_A$ and $\mathrm p_B$ quasi-isomorphisms? 
The previous question generalizes, in the framework of $A_\infty$-algebras and modules, a similar problem for DG algebras and DG modules, solved by Keller in~\cite{Keller}, and recently brought to attention in the framework of deformation quantization by Shoikhet~\cite{Sh}.

In fact, when $A$ and $B$ are DG algebras and $K$ is a DG $A$-$B$-bimodules, we may consider the DG category $\texttt{Cat}(A,B,K)$ as in Section~\ref{s-2}.
Analogously, we may consider the Hochschild cochain complex of $\texttt{Cat}(A,B,K)$ with values in itself: again, it splits into three pieces, and the Hochschild differential $\mathrm d_\gamma$, uniquely determined by the DG structures on $A$, $B$, and $K$, splits into five pieces.

Again, the two natural projections $\mathrm p_A$ and $\mathrm p_B$ from $\mathrm{C}^\bullet(\texttt{Cat}(A,B,K))$ onto $\mathrm{C}^\bullet(A,A)$ and $\mathrm{C}^\bullet(B,B)$ are morphisms of complexes: Keller~\cite{Keller} has proved that both projections are quasi-isomorphisms, if the derived left- and right-actions $\mathrm L_A$ and $\mathrm R_B$ from $A$ and $B$ to $\mathrm{RHom}_{-B}^\bullet(K,K)$ and $\mathrm{RHom}_{A-}^\bullet(K,K)^\mathrm{op}$ respectively are quasi-isomorphisms.
Here, e.g.\ $\mathrm{RHom}_{-B}(K,K)$ denotes the right-derived functor of $\mathrm{Hom}_{-B}(\bullet,K)$ in the derived category $\mathcal D(\texttt{Mod}_B)$ of the category $\texttt{Mod}_B$ of graded right $B$-modules, whose spaces of morphisms are specified by
\[
\mathrm{Hom}_{-B}(V,W)=\bigoplus_{p\in\mathbb Z}\mathrm{Hom}^p_{-B}(V,W)=\bigoplus_{p\in\mathbb Z}\mathrm{hom}_{-B}(V,W[p]).
\]
The cohomology of the complex $\mathrm{RHom}_{-B}^\bullet(K,K)$ computes the derived functor $\mathrm{Ext}_{-B}^\bullet(K,K)$; accordingly, $\mathrm L_A$ denotes the derived right action of $A$ on $K$ in the framework of derived categories.

We observe that the DG algebras $\underline{\mathrm{End}}_{-B}(K)$ and $\underline{\mathrm{End}}_{A-}(K)$ represent respectively $\mathrm{RHom}_{-B}(K,K)$ and $\mathrm{RHom}_{A-}(K,K)^\mathrm{op}$, taking the Bar resolution of $K$ in $\texttt{Mod}_B$ and ${}_A\texttt{Mod}$ respectively (of course, the product structure on $\mathrm{RHom}_{A-}(K,K)^\mathrm{op}$ is induced by the opposite of Yoneda product).
Thus, the derived left- and right-action in the $A_\infty$-framework truly generalize the corresponding derived left- and right-action in the case of a DG category, with the obvious advantage of providing explicit formul\ae\ involving homotopies.   
Furthermore, in the framework of derived categories, the derived left- and right-actions $\mathrm L_A$ and $\mathrm R_B$ induce structures of DG bimodule on $\mathrm{RHom}_{-B}(K,K)$ and $\mathrm{RHom}_{A-}^\bullet(K,K)^\mathrm{op}$ in a natural way; further, two components of the Hochschild differential $\mathrm d_\gamma$ are determined by composition with $\mathrm L_A$ and $\mathrm R_B$.
\begin{Thm}\label{t-keller}
Given $A$, $B$ and $K$ as above, where $A$ and $B$ are assumed to be flat, if $\mathrm L_A$, resp.\ $\mathrm R_B$, is a quasi-isomorphism, the canonical projection 
$$
\mathrm p_B\,:\,\mathrm{C}^\bullet(\texttt{Cat}_\infty(A,B,K))\,\twoheadrightarrow\,\mathrm{C}^\bullet(B,B),\ \text{resp.}\ \mathrm p_A\,:\,\mathrm{C}^\bullet(\texttt{Cat}_\infty(A,B,K))\,\twoheadrightarrow\,\mathrm{C}^\bullet(A,A),
$$
is a quasi-isomorphism. 
\end{Thm}
\begin{proof}
We prove the claim for the derived left-action; the proof of the claim for the derived right-action is almost the same, with obvious due changes.

Since $\mathrm p_B$ is a chain map, that it is a quasi-isomorphism tantamounts to the acyclicity of $\mathrm{Cone}^\bullet(\mathrm p_B)$, the cone of $\mathrm p_B$. 
First of all, $\mathrm{Cone}^\bullet(\mathrm p_B)$ is quasi-isomorphic to the subcomplex $\mathrm{Ker}(\mathrm p_B)$.\footnote{Namely, as in \cite{Sh}, we regard $\mathrm{Cone}(\mathrm p_B)$ as a bicomplex 
with vertical differential being the sum of the corresponding Hochschild differentials of the two complexes involved and horizontal differential being $\mathrm p_B[1]$. It has only $2$ columns, hence the associated spectral sequence stabilizes at 
$E_2$, and moreover, $E_1$ coincides with $\mathrm{Ker}(\mathrm p_B)$.}

We observe that $\mathrm{Ker}(\mathrm p_B)=\mathrm{C}^\bullet(A,B,K)\oplus \mathrm{C}^\bullet(A,A)$: by the arguments of Subsection~\ref{ss-2-1}, $\mathrm{C}^\bullet(A,B,K)$ is a subcomplex thereof. 
Lemma~\ref{l-hoch-mix}, Subsection~\ref{ss-3-2} yields the isomorphism of complexes 
$$
\mathrm{C}^\bullet(A,B,K)\cong \mathrm{C}^\bullet(A,\underline{\mathrm{End}}_{-B}(K)).
$$
As already observed, composition with the derived left action $\mathrm L_A$ defines a morphism of complexes from $\mathrm{C}^\bullet(A,A)$ to $\mathrm{C}^\bullet(A,\underline{\mathrm{End}}_{-A}(K))$: from this, and by the arguments of Subsection~\ref{ss-2-1}, it is easy to see that $\mathrm{C}^\bullet(A,B,K)\oplus \mathrm{C}^\bullet(A,A)$ is precisely the cone of the morphism induced by composition with $\mathrm L_A$, which we denote improperly by $\mathrm{Cone}(\mathrm L_A)$.

It is now a standard fact that, for any $A_\infty$-quasi-isomorphism of $A_\infty$-algebras $A\to B$, the induced cochain map $\mathrm C^\bullet (A,A)\to \mathrm C^\bullet(A,B)$ is a quasi-isomorphism, where $B$ is viewed as an $A_\infty$-$A$-$A$-bimodule as explained at the end of Subsection~\ref{ss-3-1}.

Therefore, $\mathrm{Cone}(\mathrm L_A)$ is quasi-isomorphic to the cone of the identity map on $\mathrm{C}^\bullet(A,A)$, which is obviously acyclic. 
\end{proof}

\begin{Rem}
If the $A_\infty$-structures on $A$, $B$ or $K$ are allowed to have infinitely many Taylor components, then Theorem \ref{t-keller} remains true for completed variants of Hochschild complexes. 
E.g., if the $A_\infty$-morphism ${\rm L}_A:A\to\widehat{\underline{\rm End}}_{-B}(K)$ is a quasi-isomorphism, then the canonical projection 
${\rm p}_B:\widehat{\mathrm{C}}({\tt Cat}_\infty(A,B,K))\to\widehat{\mathrm{C}}(B,B)$ is a quasi-isomorphism. 
The proof is completely parallel, with due modifications (e.g., $\mathrm{C}({\underline{\rm End}}_{-B}(K))$ being replaced by $\widehat{\mathrm{C}}(\widehat{\underline{\rm End}}_{-B}(K))$). 
\end{Rem}

\section{Configuration spaces, their compactifications and colored propagators}\label{s-4}

In this Section we discuss in some details compactifications of configuration spaces of points in the complex 
upper-half plane $\mathbb H$ and on the real axis $\mathbb R$.

We will focus our attention on Kontsevich's Eye $\mathcal C_{2,0}$ and on the I-cube $\mathcal C_{2,1}$, in order to better formulate the properties of the $2$-colored and $4$-colored propagators, which will play a central role in the proof of the main result.

\subsection{Configuration spaces and their compactifications}\label{ss-4-1}
In this Subsection we recall compactifications of configuration spaces of points in the complex 
upper-half plane $\mathcal H$ and on the real axis $\mathbb R$.

We consider a finite set $A$ and a finite (totally) ordered set $B$.
We define the open configuration space $C_{A,B}^+$ as 
$$
C_{A,B}^+:=\mathrm{Conf}_{A,B}^+/G_2=\left\{(p,q)\in \mathbb H^A\times \mathbb R^B\,|\,p(a)\neq p(a')\textrm{ if }a\neq a'\,,\,
q(b)< q(b')\textrm{ if }b<b'\right\}/G_2,
$$
where $G_2$ is the semidirect product $\mathbb R^+\ltimes \mathbb R$, which acts diagonally on $\mathbb H^A\times \mathbb R^B$ {\em via}
\[
(\lambda,\mu)(p,q)=(\lambda p+\mu,\lambda q+\mu)\qquad(\lambda\in\mathbb R^+,\ \mu\in\mathbb R)\,.
\]
The action of the $2$-dimensional Lie group $G_2$ on such $n+m$-tuples is free, precisely when $2|A|+|B|-2\geq 0$: 
in this case, $C_{A,B}^+$ is a smooth real manifold of dimension $2|A|+|B|-2$. 
(Of course, when $|B|$ is either $0$ or $1$, then we may simply drop the suffix $+$, as no ordering of $B$ is involved.)

The configuration space $C_A$ is defined as
$$
C_A:=\left\{p\in\mathbb C^A\,|\,p(a)\neq p(a')\textrm{ if }a\neq a'\right\}/G_3,
$$
where $G_3$ is the semidirect product $\mathbb R^+\ltimes \mathbb C$, which acts diagonally on $\mathbb C^A$ {\em via} 
\[
(\lambda,\mu)p=\lambda p+\mu\qquad(\lambda\in\mathbb R^+,\ \mu\in\mathbb C)\,.
\] 
The action of $G_3$, which is a real Lie group of dimension $3$, is free precisely when $2|A|-3\geq 0$, 
in which case $C_A$ is a smooth real manifold of dimension $2|A|-3$. 

The configuration spaces $C_{A,B}^+$, resp.\ $C_A$, admit compactifications \`a la Fulton--MacPherson, 
obtained by successive real blow-ups: we will not discuss here the construction of their compactifications 
$\mathcal C_{A,B}^+$, $\mathcal C_A$, which are smooth manifolds with corners, referring to~\cite{K,CR} for more details, but we focus mainly on their stratification, in particular on the boundary strata 
of codimension $1$ of $\mathcal C_{A,B}^+$. 

Namely, the compactified configuration space $\mathcal C_{A,B}^+$ is a stratified space, and its boundary 
strata of codimension $1$ look like as follows:
\begin{enumerate}
\item[$i)$] there are a subset $A_1$ of $A$ and an ordered subset $B_1$ of successive elements of $B$, such that 
\begin{equation}\label{eq-upbound1}
\partial_{A_1,B_1}\mathcal C_{A,B}^+\cong
\mathcal C_{A_1,B_1}^+\times \mathcal C_{A\smallsetminus A_1,B\smallsetminus B_1\sqcup \{*\}}^+:
\end{equation}
intuitively, this corresponds to the situation, where points in $\mathbb H$, labelled by $A_1$, and successive points in $\mathbb R$ 
labelled by $B_1$, collapse to a single point labelled by $*$ in $\mathbb R$.
Obviously, we must have $2|A_1|+|B_1|-2\geq 0$ and $2(|A|-|A_1|)+(|B|-|B_1|+1)-2\geq 0$.
\item[$ii)$] there is a subset $A_1$ of $A$, such that
\begin{equation}\label{eq-upbound2}
\partial_{A_1}\mathcal C_{A,B}^+\cong \mathcal C_{A_1}\times \mathcal C_{A\smallsetminus A_1\sqcup \{*\},B}^+:
\end{equation}
this corresponds to the situation, where points in $\mathbb H$, labelled by $A_1$, collapse together to a single point $*$ in $\mathbb H$, labelled by $*$.
Again, we must have $2|A_1|-3\geq 0$ and $2(|A|-|A_1|+1)+|B|-2\geq 0$.
\end{enumerate}

\subsection{Orientation of configuration spaces}\label{ss-4-2}
We now spend some words for the description of the orientation of (compactified) configuration spaces $\mathcal C^+_{A,B}$ and of their boundary strata of codimension $1$.

For this purpose, we follow the patterns of~\cite{AMM}: we consider the (left) principal $G_2$-bundle $\mathrm{Conf}^+_{A,B}\to C^+_{A,B}$, and we define an orientation on the (open) configuration space $C^+_{A,B}$ in such a way that any trivialization of the $G_2$-bundle $\mathrm{Conf}^+_{A,B}$ is orientation-preserving.

We observe that $i)$ the real, $2$-dimensional Lie group $G_2$ is oriented by the volume form $\Omega_{G_2}=\mathrm db\mathrm d a$, where a general element of $G_2$ is denoted by $(a,b)$, $a\in \mathbb R^+$, $b\in \mathbb R$, and $ii)$ the real, $2n+m$-dimensional manifold $\mathrm{Conf}^{+}_{n,m}$ is oriented by the volume form $\Omega_{\mathrm{Conf}^{+}_{n,m}}=\mathrm d^2z_1\cdots\mathrm d^2 z_n\mathrm d x_1\cdots\mathrm d x_m$, $\mathrm d^2 z_i=\mathrm d\mathrm{Re}z_i\mathrm d\mathrm{Im}z_i$, $z_i$ in $\mathbb H$, $x_j$ in $\mathbb R$.

We only recall, without going into the details, that there are three possible choices of global sections of $\mathrm{Conf}^+_{n,m}$, to which correspond three orientation forms on $C_{n,m}^+$ and on $\mathcal C_{n,m}^+$.

\subsubsection{Orientation of boundary strata of codimension $1$}\label{sss-4-2-1}
We recall the discussion on the boundary strata of codimension $1$ of $\mathcal C_{A,B}^+$, for a finite subset $A$ of $\mathbb N$ and a finite, ordered subset $B$ of $\mathbb N$ at the end of Subsection~\ref{ss-4-1}.

Therefore, we are interested in determining the induced orientation on the two types of boundary strata~\eqref{eq-upbound1} and~\eqref{eq-upbound2}.
In fact, we want to compare the natural orientation of the boundary strata of codimension $1$, induced from the orientation of $\mathcal C_{A,B}^+$, with the product orientation coming from the identifications~\eqref{eq-upbound1} and~\eqref{eq-upbound2}.

We may quote the following results of~\cite{AMM}, Section I.2.
\begin{Lem}\label{l-or}
Borrowing notation and convention from Subsection~\ref{ss-4-1},
\begin{itemize}
\item[$i)$] for boundary strata of type~\eqref{eq-upbound1}, 
\begin{equation}\label{eq-b1-or}
\Omega_{\partial_{A_1,B_1}\mathcal C^{+}_{A,B}}=(-1)^{j(|B_1|+1)-1}\Omega_{\mathcal C^{+}_{A_1,B_1}}\wedge \Omega_{C^{+}_{A\smallsetminus A_1,B\smallsetminus B_1\sqcup\{*\}}},
\end{equation}
where $j$ is the minimum of $B_1$;
\item[$ii)$] for boundary strata of type~\eqref{eq-upbound2}, 
\begin{equation}\label{eq-b2-or}
 \Omega_{\partial_{A_1}\mathcal C^{+}_{A,B}}=-\Omega_{\mathcal C_{A_1}}\wedge \Omega_{\mathcal C^{+}_{A\smallsetminus A_1\sqcup\{*\},B}}.
\end{equation}
\end{itemize}
\end{Lem}

\subsection{Explicit formul\ae\ for the colored propagators}\label{ss-4-3}
In the present Subsection we define and discuss the main properties of $i)$ the $2$-colored propagators and $ii)$ the $4$-colored propagators, which play a fundamental role in the constructions of Sections~\ref{s-5} and~\ref{s-6}. 

\subsubsection{The $2$-colored propagators}\label{sss-4-3-1}
We need first an explicit description of the compactified configuration space $\mathcal C_{2,0}$, known as {\bf Kontsevich's Eye}. 
Here is a picture of it, with all boundary strata of codimension $1$, labelled by Greek letters:
\begin{center}
\includegraphics[scale=0.36]{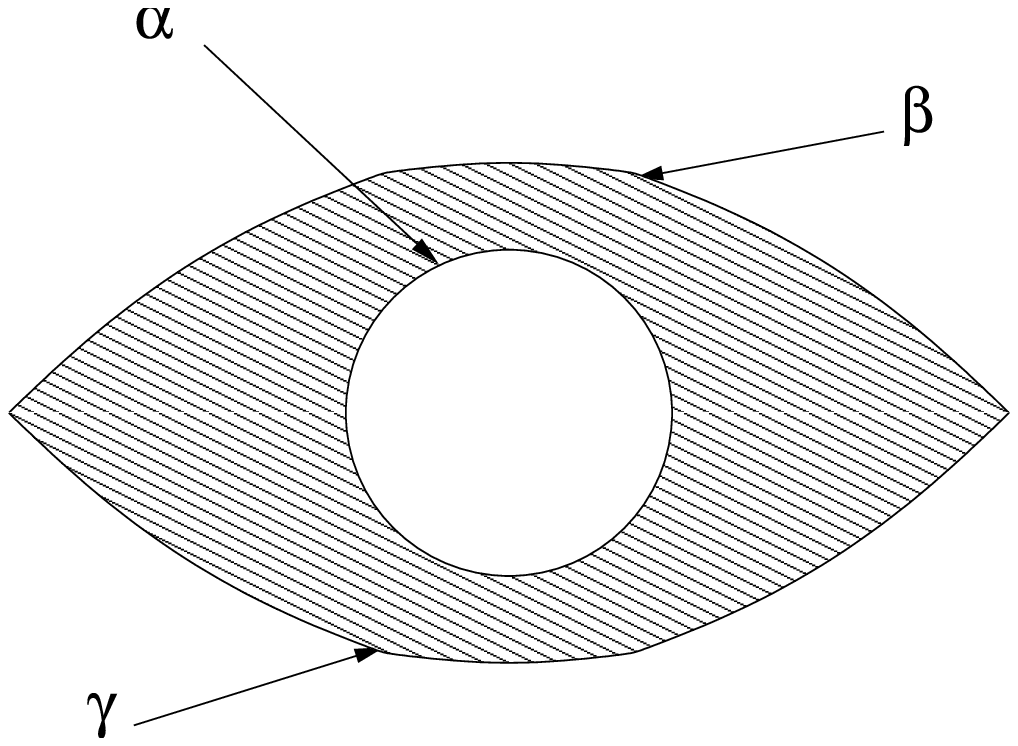} \\
\text{Figure 1 - Kontsevich's eye} \\
\end{center}
We now describe the boundary strata of $\mathcal C_{2,0}$ of codimension $1$, namely
\begin{itemize}
\item[$i)$] The stratum labelled by $\alpha$ corresponds to $\mathcal C_2=S^1$: intuitively, it describes the situation, where the two points collapse to a single point in $\mathbb H$; 
\item[$ii)$] the stratum labelled by $\beta$ corresponds to $\mathcal C_{1,1}\cong[0,1]$: it describes the situation, where the first point goes to $\mathbb R$; 
\item[$iii)$] the stratum labelled by $\gamma$ corresponds to $\mathcal C_{1,1}\cong[0,1]$: it describes the situation, where the second point goes to $\mathbb R$. 
\end{itemize}
For any two distinct points $z,w$ in $\mathbb H\sqcup\mathbb R$, we set
\begin{equation*}
\varphi^+(z,w)=\frac1{2\pi}\mathrm{arg}\!\left(\frac{z-w}{\overline z-w}\right),\ \varphi^-(z,w)=\varphi^+(w,z).
\end{equation*}
We observe that the real number $\varphi^+(z,w)$ represents the (normalized) angle from the geodesic from $z$ to the point $\infty$ on the positive imaginary axis to the geodesic between $z$ and $w$ w.r.t.\ the hyperbolic metric of 
$\mathcal H\sqcup\mathbb R$, measured in counterclockwise direction.
Both functions are well-defined up to the addition of constant terms, therefore $\omega^\pm:={\rm d}\varphi^\pm$ are well-defined $1$-forms, which are obviously basic w.r.t.\ the action of $G_2$: in summary, $\omega^\pm$ are well-defined $1$-forms on the open configuration space $C_{2,0}$. 
\begin{Lem}\label{l-angle}
The $1$-forms $\omega^\pm$ extend to smooth $1$-forms on Kontsevich's eye $\mathcal C_{2,0}$, with the 
following properties: 
\begin{itemize}
\item[$i)$] 
\begin{equation}\label{eq-pupil}
\omega^\pm\vert_\alpha=\pi_1^*(\mathrm d\varphi),
\end{equation}
where $\mathrm d \varphi$ denotes the (normalized) angle measured in counterclockwise direction from the positive imaginary axis, and $\pi_1$ is the projection from $\mathcal C_2\times\mathcal C_{1,0}$ onto the first factor.
\item[$ii)$] 
\begin{equation}\label{eq-eyelids}
\omega^+\vert_\beta=0,\ \omega^-\vert_\gamma=0.
\end{equation}
\end{itemize}
\end{Lem}
\begin{proof}
We first observe that $\omega^+$ is the standard angle form of Kontsevich, see e.g.~\cite{K}, whence it is a smooth form on $\mathcal C_{2,0}$, which enjoys the properties~\eqref{eq-pupil} and~\eqref{eq-eyelids}.

On the other hand, by definition, $\omega^-=\tau^*\omega^+$, where $\tau$ is the involution of $\mathcal C_{2,0}$, which extends smoothly the involution $(z,w)\mapsto (w,z)$ on $\mathrm{Conf}_{2,0}$: then, the smoothness of $\omega^-$ as well as properties~\eqref{eq-pupil} and~\eqref{eq-eyelids} follow immediately.
\end{proof}
We refer to~\cite{CF-br} for the physical origin of the $2$-colored Kontsevich propagators: we only mention that they arise from the Poisson Sigma model in the presence of a brane (i.e.\ a coisotropic submanifold of the target Poisson manifold) dictating boundary conditions for the fields.

\subsubsection{The $4$-colored propagators}\label{sss-4-3-2}
We now want to describe the so-called $4$-colored propagators: for an explanation of their physical origin, which is traced back to boundary conditions for the Poisson Sigma model dictated by two branes (i.e.\ two coisotropic submanifolds of the target Poisson manifold), we refer once again to~\cite{CF-br}.

Here, we are mainly interested in their precise construction and their properties: for this purpose, we find an appropriate compactified configuration space, to which the {\em na\"ive} definition of the $4$-colored propagators extend smoothly.

\

\paragraph{{\bf Description of the I-cube}}\label{p-4-3-2-1}

We shortly describe the compactified configuration space $\mathcal C_{2,1}$ of $2$ distinct points in the complex upper half-plane $\mathbb H$ and one point on the real axis $\mathbb R$: by construction, it is a smooth manifold with corners of real dimension $3$, which will be called the {\bf I-cube}.

Pictorially, the I-cube looks like as follows:
\bigskip
\begin{center}
\resizebox{0.35 \textwidth}{!}{\input{I-cube_Kosz.pstex_t}}\\
\text{Figure 2 - The I-cube $\mathcal C_{2,1}$} \\
\end{center}
\bigskip
Its boundary stratification consists of $9$ strata of codimension $1$, $20$ strata of codimension $2$ and $12$ 
strata of codimension $3$: we will explicitly describe only the boundary strata of codimension $1$, the boundary strata of higher codimensions can be easily characterized by inspecting the former strata. 

Before describing the boundary strata of $\mathcal C_{2,1}$ of codimension $1$ mathematically, it is better to depict them:
\bigskip
\begin{center}
\resizebox{0.85 \textwidth}{!}{\input{I-cube_Kosz_1.pstex_t}}\\
\text{Figure 3 - Boundary strata of the I-cube of codimension $1$} \\
\end{center}
\bigskip
The boundary stratum labelled by $\alpha$ factors as $\mathcal C_2\times \mathcal C_{1,1}$: since $\mathcal C_2=S^1$ and $\mathcal C_{1,1}$ is a closed interval, $\alpha$ is a cylinder.
\begin{Rem}\label{r-angle}
We consider the open configuration space $C_{1,1}\cong \{e^{\mathrm i t}:\ t\in (0,\pi)\}$: on it, we take the closed $1$-form $\frac{1}{2\pi}\mathrm d t$.
It extends smoothly to a closed $1$-form $\rho$ on the compactified configuration space $\mathcal C_{1,1}$, which vanishes on its two boundary strata of codimension $1$: these properties will play a central role in subsequent computations.
\end{Rem}
The boundary strata labelled by $\beta$ and $\gamma$ are both described by $\mathcal C_{2,0}\times \mathcal C_{0,2}^+$, the only difference being the position of the cluster corresponding to $\mathcal C_{2,0}$ w.r.t.\ the point $x$ on $\mathbb R$: since $\mathcal C_{0,2}^+$ is $0$-dimensional, the strata $\alpha$ and $\beta$ are two copies of Kontsevich's eye $\mathcal C_{2,0}$.

The boundary strata labelled by $\delta$ and $\varepsilon$ are both described by $\mathcal C_{1,1}\times\mathcal C_{1,1}$, depending on whether the point labelled by $1$ or $2$ collapses to the point $x$ on the real axis: since $\mathcal C_{1,1}$ is a closed interval, both $\delta$ and $\varepsilon$ are two squares.

Finally, the boundary strata labelled by $\eta$, $\theta$, $\zeta$ and $\xi$ factor as $\mathcal C_{1,2}^+\times\mathcal C_{1,0}$, depending on whether the point labelled by $1$ or $2$ goes to the real axis either on the left or on the right of $x$: since $\mathcal C_{1,0}$ is $0$-dimensional, these boundary strata correspond to $\mathcal C_{1,2}^+$.
The latter compactified configuration space is a hexagon: this is easily verified by direct inspection of its boundary stratification. 

\

\paragraph{{\bf Explicit formul\ae\ for the $4$-colored propagators}}\label{p-4-3-2-2}
First of all, we observe that there is a projection $\pi_{2,0}$ from $\mathcal C_{2,1}$ onto $\mathcal C_{2,0}$, which extends smoothly the obvious projection from $C_{2,1}$ onto $C_{2,0}$ forgetting the point $x$ on the real axis.
It makes therefore sense to set
\[
\omega^{+,+}=\pi_{2,0}^*(\omega^+),\ \omega^{-,-}=\pi_{2,0}^*(\omega^-).
\]
We further consider a triple $(z,w,x)$, where $z$, $w$ are two distinct points in $\mathbb H$ and $x$ is a point on $\mathbb R$.
We recall that the complex function $z\mapsto \sqrt z$ is a well-defined holomorphic function on $\mathbb H$, mapping $\mathbb H$ to the first quadrant $\mathbb Q^{+,+}$ of the complex plane, whence it makes sense to consider the $1$-forms
\begin{align*}
\omega^{+,-}(z,w,x)&=\frac{1}{2\pi}\mathrm d\,\mathrm{arg}\!\left(\frac{\sqrt{z-x}-\sqrt{w-x}}{\sqrt{z-x}-\overline{\sqrt{w-x}}}\frac{\sqrt{z-x}+\overline{\sqrt{w-x}}}{\sqrt{z-x}+\sqrt{w-x}}\right),\\
\omega^{-,+}(z,w,x)&=\frac{1}{2\pi}\mathrm d\,\mathrm{arg}\!\left(\frac{\sqrt{z-x}-\sqrt{w-x}}{\sqrt{z-x}+\overline{\sqrt{w-x}}}\frac{\sqrt{z-x}-\overline{\sqrt{w-x}}}{\sqrt{z-x}+\sqrt{w-x}}\right).
\end{align*}
Thus, $\omega^{+,-}$ and $\omega^{-,+}$ are smooth forms on the open configuration space $\mathrm{Conf}_{2,1}$.
We recall that there is an action of the $2$-dimensional Lie group $G_2$ on $\mathrm{Conf}_{2,1}$: it is not difficult to verify that both $1$-forms $\omega^{+,-}$ and $\omega^{-,+}$ are basic w.r.t.\ the action of $G_2$, hence they both descend to smooth forms on the open configuration spaces $C_{2,1}$.

In the following Lemma, we use the convention that the point in $\mathbb H$ labelled by $1$, resp.\ $2$, corresponds to the initial, resp.\ final, argument in $\mathbb H$ of the forms under consideration.
\begin{Lem}\label{l-CF}
The $1$-forms $\omega^{+,+}$, $\omega^{+,-}$, $\omega^{-,+}$ and $\omega^{-,-}$ extend smoothly to the I-cube $\mathcal C_{2,1}$.
They further enjoy the following properties:
\begin{itemize}
\item[$i)$] \begin{equation}\label{eq-alpha}
\omega^{+,+}\vert_\alpha=\pi_1^*(\mathrm d\varphi),\ \omega^{+,-}\vert_\alpha=\pi_1^*(\mathrm d\varphi)-\pi_2^*(\rho),\  \omega^{-,+}\vert_\alpha=\pi_1^*(\mathrm d\varphi)-\pi_2^*(\rho),\ \omega^{-,-}\vert_\alpha=\pi_1^*(\mathrm d\varphi),
\end{equation} 
where $\pi_i$, $i=1,2$, denotes the projection onto the $i$-th factor of the decomposition $\mathcal C_2\times \mathcal C_{1,1}$ of the boundary stratum $\alpha$, and $\rho$ is the smooth $1$-form on $\mathcal C_{1,1}$ discussed in Remark~\ref{r-angle}. 
\item[$ii)$] \begin{equation}\label{eq-alpha-beta}
\begin{aligned}
\omega^{+,+}\vert_\beta&=\omega^+,\ & \omega^{+,-}\vert_\beta&=\omega^+,\ & \omega^{-,+}\vert_\beta&=\omega^-,\ & \omega^{-,-}\vert_\beta&=\omega^-\quad \text{and}\\
\omega^{+,+}\vert_\gamma&=\omega^+,\ & \omega^{+,-}\vert_\gamma&=\omega^-,\ & \omega^{-,+}\vert_\gamma&=\omega^+,\ & \omega^{-,-}\vert_\gamma&=\omega^-.,
\end{aligned}
\end{equation}  
where we implicitly identify both boundary strata with Kontsevich's Eye, see also Subsubsection~\ref{sss-4-3-1}.
\item[$iii)$] \begin{equation}\label{eq-gamma-delta}
\begin{aligned}
&\omega^{+,+}\vert_\delta=\omega^{+,-}\vert_\delta=\omega^{-,+}\vert_\delta=0,\\
&\omega^{+,-}\vert_\varepsilon=\omega^{-,+}\vert_\varepsilon=\omega^{-,-}\vert_\varepsilon=0.
\end{aligned}
\end{equation} 
\item[$iv)$] \begin{equation}\label{eq-eps}
\begin{aligned}
&\omega^{+,-}\vert_\eta=\omega^{-,-}\vert_\eta=0,\ & &\omega^{+,+}\vert_\theta=\omega^{-,+}\vert_\theta=0,\\
&\omega^{-,+}\vert_\zeta=\omega^{-,-}\vert_\zeta=0,\ & &\omega^{+,+}\vert_\xi=\omega^{+,-}\vert_\xi=0.
\end{aligned}
\end{equation}
\end{itemize}
\end{Lem}
\begin{proof}
First of all, since the projection $\pi_{2,0}:\mathcal C_{2,1}\to \mathcal C_{2,0}$ is smooth, Lemma~\ref{l-angle}, Subsubsection~\ref{sss-4-3-1}, implies immediately that $\omega^{+,+}$ and $\omega^{-,-}$ are smooth $1$-forms on $\mathcal C_{2,1}$.
Lemma~\ref{l-angle}, Subsection~\ref{sss-4-3-1}, yields also immediately Properties $i)$, $ii)$, $iii)$ and $iv)$ of $\omega^{+,+}$ and $\omega^{-,-}$.

It remains to prove smoothness of $\omega^{+,-}$ and $\omega^{-,+}$ on $\mathcal C_{2,1}$ and Properties $i)$, $ii)$, $iii)$ and $iv)$.
We prove the statements e.g.\ for $\omega^{-,+}$: similar computations lead to the proof of the statements for $\omega^{-,+}$.

In order to prove all statements, we make use of local coordinates of $\mathcal C_{2,1}$ near the boundary strata of codimension $1$ in all cases.

We begin by considering the boundary stratum labelled by $\alpha$: local coordinates of $\mathcal C_{2,1}$ near $\alpha$ are specified {\em via}
\[
\mathcal C_2\times \mathcal C_{1,1}\cong S^1\times [0,\pi]\ni (\varphi,t)\mapsto \left[\left(e^{\mathrm i t},e^{\mathrm i t}+\varepsilon e^{\mathrm i \varphi},0\right)\right]\in \mathcal C_{2,1},\ \varepsilon>0,
\]
where $\alpha$ is recovered, when $\varepsilon$ tends to $0$.
We have implicitly used local sections of $C_{1,1}$ and $C_{2,1}$: the point on $\mathbb R$ has been put at $0$, and the first point in $\mathbb H$ has been put on a circle of radius $1$ around $0$.

Then, using the standard notation $[(z,w,x)]$ for a point in $\mathcal C_{2,1}$, we have
\[
\sqrt{w-x}=\sqrt{z-x+\varepsilon e^{\mathrm i\varphi}}=\sqrt{z-x}+\varepsilon \frac{1}{2|z-x|}e^{\mathrm i\left(\varphi-\frac{1}2 t\right)}+\mathcal O(\varepsilon^2),\ z=e^{\mathrm i t},\ x=0.
\]
Substituting the right-most expression in the definition of $\omega^{-,+}$ and taking the limit as $\varepsilon$ tends to $0$, we get
\[
\omega^{-,+}\vert_\alpha=\frac{1}{2\pi}\left(\mathrm d\varphi-\mathrm dt\right)=\pi_1^*\mathrm d\varphi-\pi_2^*(\rho),
\]  
where $\rho$ is the smooth $1$-form discussed in Remark~\ref{r-angle}.
We observe that, in the last equality, we have abused the notation $\mathrm d\varphi$, in order to be consistent with the notation of Lemma~\ref{l-angle}, Subsubsection~\ref{sss-4-3-1}.

We now consider e.g.\ the boundary strata labelled by $\beta$ and $\gamma$.
Local coordinates of $\mathcal C_{2,1}$ near $\beta$, resp.\ $\gamma$, are specified {\em via}
\[
\begin{aligned}
\mathcal C_{2,0}\times \mathcal C_{0,2}^+\cong \mathcal C_{2,0}\times \{-1,0\}\ni \left((\mathrm i,\mathrm i+\rho e^{\mathrm i\varphi}),(-1,0)\right)&\mapsto \left[\left(-1+\varepsilon \mathrm i,-1+\varepsilon (\mathrm i+\rho e^{\mathrm i\varphi}),0)\right)\right]\in \mathcal C_{2,1},\quad \text{resp.}\\
\mathcal C_{2,0}\times \mathcal C_{0,2}^+\cong \mathcal C_{2,0}\times \{0,1\}\ni \left((\mathrm i,\mathrm i+\rho e^{\mathrm i\varphi}),(0,1)\right)&\mapsto \left[\left(1+\varepsilon \mathrm i,1+\varepsilon (\mathrm i+\rho e^{\mathrm i\varphi}),0)\right)\right]\in \mathcal C_{2,1},\ \rho,\varepsilon>0,
\end{aligned}
\]
where again $\beta$ and $\gamma$ are recovered, when $\varepsilon$ tends to $0$.
(Once again, we have made use of local sections of the interior of $\mathcal C_{2,0}$ and $\mathcal C_{2,1}$.)

Using the standard notation for a general point in (the interior of) $\mathcal C_{2,1}$, we have, near the boundary stratum $\beta$, resp.\ $\gamma$,
\[
\begin{aligned}
&\sqrt{z-x}=\sqrt{y-x+\varepsilon \widetilde z}=\mathrm i-\varepsilon \frac{\mathrm i\widetilde z}2 +\mathcal O(\varepsilon^2),\ & &\sqrt{w-x}=\sqrt{y-x+\varepsilon \widetilde w}=\mathrm i-\varepsilon \frac{\mathrm i\widetilde w}2 +\mathcal O(\varepsilon^2),\quad \text{resp.}\\
&\sqrt{z-x}=\sqrt{y-x+\varepsilon \widetilde z}=1+\varepsilon \frac{\widetilde z}2 +\mathcal O(\varepsilon^2),\ & &\sqrt{w-x}=\sqrt{y-x+\varepsilon \widetilde w}=1+\varepsilon \frac{\widetilde w}2+\mathcal O(\varepsilon^2),
\end{aligned}
\]
where $\tilde z=\mathrm i$ and $\tilde w=\mathrm i+\rho e^{\mathrm i \varphi}$, $y=-1$ for $\beta$ and $y=1$ for $\gamma$, and $x=0$.

Substituting the right-most expressions on all previous chains of identities in $\omega^{+,-}$ and $\omega^{-,+}$, and letting $\varepsilon$ tend to $0$, we obtain $ii)$: in particular, the restrictions of $\omega^{+,-}$ and $\omega^{-,+}$ to $\beta$ and $\gamma$ are smooth $1$-forms.

We now consider the boundary strata labelled by $\delta$ and $\varepsilon$.
Local coordinates of $\mathcal C_{2,1}$ near $\delta$, resp.\ $\varepsilon$, are specified {\em via}
\[
\begin{aligned}
\mathcal C_{1,1}\times \mathcal C_{1,1}\cong [0,\pi]\times [0,\pi]\ni (s,t)&\mapsto \left[\left(\rho e^{\mathrm i s},e^{\mathrm i t},0\right)\right]\in \mathcal C_{2,1},\quad \text{resp.}\\
\mathcal C_{1,1}\times \mathcal C_{1,1}\cong [0,\pi]\times [0,\pi]\ni (s,t)&\mapsto \left[\left(e^{\mathrm i s},\rho e^{\mathrm i t},0\right)\right]\in \mathcal C_{2,1},\ 
\end{aligned}
\] 
where $\delta$, resp.\ $\varepsilon$, is recovered, when $\rho$ tends to $0$.

Using again standard notation for a point in (the interior of) $\mathcal C_{2,1}$, we then get
\[
\sqrt{z-x}=\sqrt{\rho}\sqrt{\widetilde z},\quad \text{resp.}\quad \sqrt{w-x}=\sqrt{\rho}\sqrt{\widetilde w},
\]
where $\widetilde z=e^{\mathrm i s}$ and $\widetilde w=e^{\mathrm i t}$.
The remaining square roots do not contain $\rho$.

In particular, if we substitute the previous expressions in $\omega^{+,-}$ and $\omega^{-,+}$ and let $\rho$ tend to $0$, we easily obtain
\[
\begin{aligned}
\omega^{-,+}\vert_\delta=\omega^{+,-}\vert_\delta=0,\ \omega^{-,+}\vert_\varepsilon=\omega^{+,-}\vert_\varepsilon=0,
\end{aligned}
\]
which in particular implies that the restrictions of $\omega^{+,-}$ and $\omega^{-,+}$ to $\delta$ and $\varepsilon$ are smooth $1$-forms.

Finally, we consider e.g.\ the boundary stratum labelled by $\eta$.
Local coordinates nearby are specified {\em via}
\[
\mathcal C_{1,0}\times \mathcal C_{1,2}^+\cong \{\mathrm i\}\times \mathcal C_{1,2}^+\ni (\mathrm i,(z,0,1))\mapsto \left[\left(z,1+\varepsilon \mathrm i,0\right)\right]\in \mathcal C_{2,1},
\]
where $\eta$ is recovered, when $\varepsilon$ tends to $0$.
Here, we have used global sections of $C_{1,1}$, $C_{1,2}^+$ and $C_{2,1}$, using the action of $G_2$ to put the point in $\mathbb H$ to $\mathrm i$, to put the first and the second point on $\mathbb R$ to $0$ and $1$, and to put the point on $\mathbb R$ to $0$ and the real part of the second point in $\mathbb H$ to $1$.
 
Computations similar in spirit to the previous ones permit to compute explicit expressions for the restrictions of $\omega^{+,-}$ ad $\omega^{-,+}$ to $\eta$: in particular, we see that $\omega^{+,-}$ and $\omega^{-,+}$ restrict to smooth $1$-forms on $\mathcal C_{1,2}^+$, and we also get $iv)$.
\end{proof}

\

\paragraph{{\bf The $4$-colored propagators on the first quadrant}}\label{p-4-3-2-3}
We observe that the complex function $z\mapsto \sqrt z$ restricts on $\mathbb H\sqcup \mathbb R\smallsetminus\{0\}$ to a holomorphic function, whose image is $\mathcal Q^{+,+}\sqcup \mathbb R^+\sqcup \mathrm i \mathbb R^+$: the negative real axis is mapped to $\mathrm i\mathbb R^+$, the positive real axis is mapped to itself, and $\mathbb H$ is mapped to $\mathcal Q^{+,+}$.
Further, $z\mapsto \sqrt z$ is multi-valued, when considered as a function on $\mathbb C$, with $0$ as a branching point.

There is an explicit global section of the projection $\mathrm{Conf}_{2,1}\to C_{2,1}$, namely
\[
C_{2,1}\ni [(z,w,x)]\mapsto \left(\frac{z-x}{|z-x|},\frac{w-x}{|z-x|},0\right)\in \mathrm{Conf}_{2,1}.
\]
Setting $\widetilde z=\frac{z-x}{|z-x|}$ and $\widetilde w=\frac{w-x}{|z-x|}$, we get two point in $\mathbb H$: hence, setting $u=\sqrt{\widetilde z}$ and $v=\sqrt{\widetilde w}$, $u$ and $v$ lie in $\mathcal Q^{+,+}$.
We then find the alternative descriptions of the $4$-colored propagators:
\[
\begin{aligned}
\omega^{+,+}(u,v)&=\frac{1}{2\pi}\mathrm d\,\mathrm{arg}\left(\frac{u-v}{\overline u-v}\frac{u+v}{\overline u+v}\right),\ & \omega^{+,-}(u,v)&=\frac{1}{2\pi}\mathrm d\,\mathrm{arg}\left(\frac{u-v}{u-\overline v}\frac{u+\overline v}{u+v}\right),\\ 
\omega^{-,+}(u,v)&=\frac{1}{2\pi}\mathrm d\,\mathrm{arg}\left(\frac{u-v}{u+\overline v}\frac{u-\overline v}{u+v}\right),\ & \omega^{+,+}(u,v)&=\frac{1}{2\pi}\mathrm d\,\mathrm{arg}\left(\frac{u-v}{u-\overline v}\frac{u+v}{u+\overline v}\right).
\end{aligned}
\]
We observe that the previous formul\ae\ descend to the quotient of the configuration space of two points in $\mathcal Q^{+,+}$ w.r.t.\ the action of $G_1\cong \mathbb R^+$ by rescaling.

In fact, the present description of the $4$-colored propagators is the original one, see~\cite{CF-br}: we have preferred to work with the previous (apparently more complicated) description, because it is more well-suited to work with compactified configuration spaces.

We finally observe that all previous formul\ae\ are special cases of the main result in~\cite{Ferrario}, where general (super)propagators for the Poisson $\sigma$-model in the presence of $n$ branes, $n\geq 1$, are explicitly produced.

\section{$L_\infty$-algebras and morphisms}\label{s-5}

In the present Section, we briefly discuss the concept of $L_\infty$-algebra and $L_\infty$-morphism; further, we describe explicitly the two $L_\infty$-algebras (which are actual genuine DG Lie algebras) 
which will be central in the constructions of Section~\ref{s-6}. 

A DG Lie algebra $\mathfrak g$ is an object of $\widehat{\texttt{GrMod}_k}$, endowed with an endomorphism $\mathrm d_\mathfrak g:\mathfrak g\to \mathfrak g$ of degree $1$ and with a graded anti-symmetric, 
bilinear map $[\bullet,\bullet]:\mathfrak g\otimes\mathfrak g\to\mathfrak g$ of degree $0$, such that $\mathrm d_\mathfrak g$ squares to $0$, and such that 
\[
\begin{aligned}
&\mathrm d_\mathfrak g([x,y])=\left[\mathrm d_\mathfrak g(x),y\right]+(-1)^{|x|}\left[x,\mathrm d_\mathfrak g(y)\right],\\
&(-1)^{|x||z|}[[x,y],z]+(-1)^{|x||y|}[[y,z],x]+(-1)^{|z||y|}[[z,x],y]=0,
\end{aligned}
\]
for homogeneous elements $x$, $y$, $z$ of $\mathfrak g$. 
The first Identity is the graded Leibniz rule, while the second is the graded Jacobi identity.

A formal pointed $Q$-manifold is an object $V$ of $\widehat{\texttt{GrMod}_k}$, such that $\mathrm C^{+}(V)\cong \mathrm S^+(V)$ is endowed with a codifferential $Q$.
A morphism $U$ between $Q$-manifolds $(U,Q_U)$ and $(V,Q_V)$ is a coalgebra morphism $\mathrm C^{+}(V)\rightarrow \mathrm C^{+}(V^{'})$ of degree 0, intertwining $Q_U$ and $Q_V$.
\begin{Def}\label{d-L_inf}
An $L_{\infty}$-structure on an object $\mathfrak g$ of the category $\widehat{\texttt{GrMod}_k}$ is a $Q$-manifold structure on $\mathfrak g[1]$; the pair $(\mathfrak g,Q)$ is called an $L_{\infty}$-algebra.
Accordingly, a morphism $F$ between $L_\infty$-algebras $(\mathfrak g_1,Q_1)$ and $(\mathfrak g_2,Q_2)$ is a morphism between the corresponding pointed $Q$-manifolds. 
\end{Def}
The fact that $Q$ is a coderivation on $\mathrm S^+(\mathfrak g[1])$ implies that $Q$ is uniquely determined by its Taylor components $Q_n:\mathrm S^n(\mathfrak g[1])\to \mathfrak g[1]$: an explicit formula for recovering $Q$ from its Taylor components may be found e.g.\ in~\cite{K,Dol}, we only mention that it is similar in spirit to the formul\ae\ appearing in the case of $A_\infty$-structures, although the fact that we consider the symmetric algebra causes the arising of shuffles.

Furthermore, the fact that an $L_\infty$-morphism $F:\mathfrak g_1\to\mathfrak g_2$ is a coalgebra morphism, implies that $F$ is also uniquely determined by its Taylor components $F_n:\mathrm S^n(\mathfrak g_1[1])\to \mathfrak g_2[1]$.
\begin{Rem}\label{r-DGLA}
If $(\mathfrak g,\mathrm d_\mathfrak g,[\bullet,\bullet])$ is a DG Lie algebra, $\mathfrak g$ has a structure of $L_\infty$-algebra, which we now describe explicitly: the Taylor components of the coderivation all vanish, except $Q_1$ and $Q_2$, specified {\em via}
\begin{eqnarray}
Q_1=\mathrm d_\mathfrak g,\ Q_2(x_1,x_2)=(-1)^{|x_1|}[x_1,x_2],\ x_i\in \mathfrak g^{|x_i|}=(\mathfrak g[1])^{|x_i|-1}. \label{noi}
\end{eqnarray}
In fact, it is easy to verify that $Q^2=0$ is equivalent to the compatibility between $\mathrm d_\mathfrak g$ and $[\bullet,\bullet]$ (graded Leibniz rule) and the graded Jacobi identity.
\end{Rem}
We consider an $L_\infty$-morphism $F:\mathfrak g_1\to \mathfrak g_2$ between $L_\infty$-algebras: the condition that $F$ intertwines the codifferentials $Q_1$ and $Q_2$ can be re-written as an infinite set of quadratic relations involving the Taylor coefficients of $Q_1$, $Q_2$ and $F$.

Exemplarily, assuming $\mathfrak g_i$, $i=1,2$, are DG Lie algebras, the quadratic identities of order $1$ and $2$ take the form
\begin{eqnarray}
Q^1_2(F_1(x))&=&F_1(Q_1^1(x)),\label{mm1} \\
Q^2_2(F_1(x),F_1(y))-F_1(Q_1^2(x,y)))&=&F_2(Q_1^1(x),y)+(-1)^{|x|-1}F_2(x,Q_1^1(y))-Q^1_2(F_2(x,y)),\ x,y\in\mathfrak g_1[1]. \label{mm2}
\end{eqnarray}
Equation~\eqref{mm1} is equivalent to the fact that $F_1$ is a morphism of complexes, while Equation~\eqref{mm2}) expresses the fact that $F_1$ is a morphism of GLAs up to a homotopy expressed by the Taylor component $F_2$. 

More generally, we have the following Proposition, for whose proof we refer to~\cite{AMM}.
\begin{Prop}
We consider two DG Lie algebras $(\mathfrak g_1,\mathrm d_1,[\bullet,\bullet]_1)$ and $(\mathfrak g_2,\mathrm d_2,[\bullet,\bullet]_2)$, which we also view as $L_\infty$-algebras as in Remark~\ref{r-DGLA}.

Then, a coalgebra morphism $F:\mathrm S^{+}(\mathfrak g_1[1])\to \mathrm S^{+}(\mathfrak g_2[1])$ is an $L_{\infty}$-morphism, if and only if it satisfies 
\begin{eqnarray}
&&Q_1'\left(F_n(\alpha_1,\dots,\alpha_n)\right)+\frac{1}{2}\sum_{I\sqcup J=\{1,\dots,n\}, I,J\neq \emptyset}\epsilon_{\alpha}(I,J)Q_2'(F_{|I|}(\alpha_I),F_{|J|}(\alpha_J))=\nonumber \\
&&\sum_{k=1}^n\sigma_{\alpha}(k,1,\dots,\hat{k},\dots,n)F_n(Q_1(\alpha_k),\alpha_1,\dots,\widehat{\alpha_k},\dots,\alpha_n)+\nonumber \\
&&+\frac{1}{2}\sum_{k\neq l}\sigma_{\alpha}(k,l,1,\dots,\hat{k},\dots,\hat{l},\dots,n)F_{n-1}(Q_2(\alpha_k,\alpha_l),\alpha_1,\dots,\widehat{\alpha_k},\dots,\widehat{\alpha_l},\dots,\alpha_n),\label{morph}
\end{eqnarray}
where $\epsilon_{\alpha}(I,J)$ denotes the sign associated to the shuffle relative to the decomposition $I\sqcup J=\{1,\dots,n\}$, and $\sigma_{\alpha}(\dots)$ denotes the sign associated to the permutation in $(\dots)$, see Section~\ref{s-1}.
\end{Prop}

\subsection{The DG Lie algebras $T_\mathrm{poly}(X)$, for $X=k^d$}\label{ss-5-1}
We consider now a ground field $k$ of characteristic $0$, which contains $\mathbb R$ or $\mathbb C$; we further set $X=k^d$.

To $X$, we associate the DG Lie algebra $T_\mathrm{poly}(X)$ of poly-vector fields on $X$ with shifted degree.
More precisely, the degree-$p$-component $T_\mathrm{poly}^p(X)$, $p\geq -1$, is $\Gamma(X,\wedge^{p+1}\mathrm TX)={\rm S}(X^*)\otimes\wedge^{p+1}(X)$, with trivial differential and Schouten-Nijenhuis bracket, 
determined by extending the Lie bracket between vector fields on $X$ as a (graded) biderivation.

Hence, $T_\mathrm{poly}(X)$ is an $L_\infty$-algebra, whose $Q$-manifold structure is
\begin{eqnarray*}
Q_1=0,\ Q_2(\alpha_1,\alpha_2):=-(-1)^{(k_1-1)(k_2)}[\alpha_2,\alpha_1]_{SN}=\alpha_1\bullet\alpha_2+(-1)^{k_1k_2}\alpha_2\bullet\alpha_1,
\end{eqnarray*}
for general elements $\alpha_1\in T_\mathrm{poly}^{k_1-1}(X)$, $\alpha_2\in T_\mathrm{poly}^{k_2-1}(X)$, where the composition $\bullet$ is 
\begin{eqnarray}
\alpha_1\bullet\alpha_2=\sum_{l=1}^{k_1}(-1)^{l-1}\alpha_1^{i_1\dots i_{k_1}}\partial_l\alpha_2^{j_1,\dots j_{k_2}}
\partial_{i_1}\wedge\dots\wedge\widehat{\partial_{i_l}}\wedge\dots\wedge\partial_{i_{k_1}}\wedge\partial_{j_1}\wedge\dots\wedge\partial_{j_{k_2}}. \label{composition}
\end{eqnarray}

We will also consider a completed version $\widehat{T}_{\rm poly}(X)$, consisting of formal poly-vector fields near the origin on $X=k^d$: 
$\widehat{T}^\bullet_\mathrm{poly}(X)=\widehat{\mathrm{S}}(X^*)\otimes\wedge^{\bullet+1}(X)$, and the previously defined Schouten--Nijenhuis bracket extends naturally to the completion $\widehat{T}^\bullet_\mathrm{poly}(X)$.

\subsection{The DG Lie algebra $\mathrm{C}^\bullet(\texttt{Cat}_{\infty}(A,B,K))$}\label{ss-5-2}
We consider again the $d$-dimensional $k$-vector space $X=k^d$; we further assume $X$ to be endowed with an inner product (hence, we may safely assume here $k=\mathbb R$ or $k=\mathbb C$).
We consider two vector subspaces $U$ and $V$ thereof, such that, w.r.t.\ the previously introduced inner product, the following decomposition holds true:
\begin{equation}\label{eq-orth-split}
X=(U\cap V)\overset{\perp}\oplus (U^\perp\cap V)\overset{\perp}\oplus (U\cap V^\perp)\overset{\perp}\oplus (U+V)^\perp.
\end{equation} 
It follows immediately from~\ref{eq-orth-split} that
\begin{equation*}
U=(U\cap V)\overset{\perp}\oplus (U\cap V^\perp),\ V=(U\cap V)\overset{\perp}\oplus (U^\perp\cap V).
\end{equation*}
To $X$, $U$ and $V$, we may associate three graded vector spaces, namely
\begin{align*}
A&=\Gamma(U,\wedge (\mathrm NU))=\mathrm S(U^*)\otimes \wedge (X/U)=\mathrm S(U^*)\otimes \wedge (U^\perp\cap V)\otimes \wedge (U+V)^\perp,\\
B&=\Gamma(V,\wedge (\mathrm NV))=\mathrm S(V^*)\otimes \wedge (X/V)=\mathrm S(V^*)\otimes \wedge (U\cap V^\perp)\otimes \wedge (U+V)^\perp,\\
K&=\Gamma(U\cap V,\wedge \left(\mathrm TX/(\mathrm TU+\mathrm TV)\right))=\mathrm S((U\cap V)^*)\otimes \wedge (U+V)^\perp,
\end{align*}
where $\mathrm T X$, resp.\ $\mathrm NU$, denotes the tangent bundle of $X$, resp.\ the normal bundle of $U$ in $\mathrm T X$.

We define a (cohomological) grading on $A$, $B$ and $K$: on $A$ and $B$, we define a grading analogously to the grading on $T_\mathrm{poly}(X)$ as in Subsection~\ref{ss-4-1}.
On the other hand, the (cohomological) grading on $K$ is defined without shifting. 

Therefore, $A$ and $B$, endowed with the trivial differential, both admit a (trivial) structure of $A_\infty$-algebra.
We now construct on $K$ a non-trivial $A_\infty$-$A$-$B$-bimodule structure.

We consider a set of linear coordinates $\{x_i\}$ on $X$, which are adapted to the orthogonal decomposition~\eqref{eq-orth-split} in the following sense: there are two non-disjoint subsets $I_i$, $i=1,2$, of $[d]$, such that 
\[
[d]=\left(I_1\cap I_2\right)\sqcup\left(I_1\cap I_2^c\right)\sqcup\left(I_1^c\cap I_2\right)\sqcup\left(I_1^c\cap I_2^c\right),
\]
and such that $\{x_i\}$ is a set of linear coordinates on $U\cap V$, $U\cap V^\perp$, $U^\perp\cap V$, $(U+V)^\perp$, if the index $i$ belongs to $I_1\cap I_2$, $I_1\cap I_2^c$, $I_1^c\cap I_2$ and $I_1^c\cap I_2^c$ respectively. 


To a general pair $(n,m)$ of non-negative integers, we associate the set $\mathcal G_{n,m}$ of admissible graphs of type $(n,m)$: a general element $\Gamma$ thereof is a directed graph (i.e.\ every edge of $\Gamma$ has an orientation), with $n$, resp. $m$, vertices of the first, resp.\ second, type.
We denote by $\mathrm E(\Gamma)$ and $\mathrm V(\Gamma)$ the set of edges and vertices of an admissible graph $\Gamma$ respectively.
\begin{Rem}\label{r-admgr}
We observe that, {\em a priori}, the admissible graphs considered here admit multiple edges (i.e.\ between any two distinct vertices there may be more than one edge) and loops (i.e.\ edges connecting a vertex of the first type to itself): as we will see, multiple edges and loops do not arise in the construction of the $A_\infty$-$A$-$B$-bimodule structure on $K$ below, but arise in Section~\ref{s-6} in the construction of a formality morphism, see later on. 
\end{Rem}
We now consider any pair of non-negative integers $(m,n)$, and to it we associate the compactified configuration space $\mathcal C_{0,m+1+n}^+$: we have $m+1+n$ ordered points on $\mathbb R$, one of which, the $m+1$-st point, plays a central role, whence the notation.
E.g.\ using the action of $G_2$ on $C_{0,m+1+n}^+$, we may put it at $x=0$.

Accordingly, we consider the set $\mathcal G_{0,m+1+n}$ of admissible graphs of type $(0,m+1+n)$: to any edge $e=(i,j)$ of a general admissible graph $\Gamma$, where the label $i$, resp.\ $j$, refers to the initial, resp.\ final, point of $e$, we associate a projection $\pi_e:\mathcal C_{0,m+1+n}^+\to \mathcal C_{0,3}^+\subset \mathcal C_{2,1}$ or $\pi_e:\mathcal C_{0,m+1,n}^+\to \mathcal C_{2,0}^+\times \mathcal C_{1,1}\subset \mathcal C_{2,1}$.

In order to define the projection $\pi_e$ precisely, we need to identify $\mathcal C_{0,3}^+$ and $\mathcal C_{2,0}^+\times \mathcal C_{1,1}$ with certain boundary strata of codimension $2$ of the I-cube $\mathcal C_{2,1}$: it is better to do this pictorially, 
\bigskip
\begin{center}
\resizebox{0.85 \textwidth}{!}{\input{I-cube_Kosz_2.pstex_t}}\\
\text{Figure 4 - Boundary strata of the I-cube of codimension $2$ needed to construct $\pi_e$} \\
\end{center}
\bigskip
Thus, for any edge $e=(i,j)$ of $\Gamma$, $i,j=1,\dots,m+1+n$, we have the following possibilities: $a)$ $1\leq i<m+1<j\leq m$, $b)$ $1\leq i<j\leq m$, $c)$ $m+1<i<j\leq m+1+n$, $d)$ $1\leq j<i\leq m$, $e)$ $1\leq j<m+1<i\leq m+1+n$, $f)$ $m+1<j<i<m+1+n$, $g)$ $m+1=j<i$, $h)$ $1\leq i<m+1=j$, $i)$ $m+1=i<j\leq m+1+n$, $j)$ $1\leq j<m+1=i$. 
We observe that the labelling of the ten cases under inspection corresponds to the labelling of the boundary strata of codimension $2$ listed above.
It is then obvious how to define the projection $\pi_e$ in all ten cases: we only observe that the vertex of the second type labelled $i$, resp.\ $j$, resp.\ $m+1$, corresponds via the projection $\pi_e$ to the vertex labelled by $1$, resp.\ $2$, resp.\ $x$ in the above picture.

This way, to every edge $e$ of an admissible graph $\Gamma$ in $\mathcal G_{0,m+1+n}$ we may associate an element $\omega^K_e$ of $\Omega^1(\mathcal C_{0,m+1+n}^+)\otimes \mathrm{End}(\mathrm T_\mathrm{poly}(X)^{\otimes m+1+n})$ {\em via} 
\[
\omega^K_e=\pi_e^*(\omega^{+,+})\otimes \tau^{I_1\cap I_2}_e+\pi_e^*(\omega^{+,-})\otimes \tau^{I_1\cap I_2^c}_e+\pi_e^*(\omega^{-,+})\otimes \tau^{I_1^c\cap I_2}_e+\pi_e^*(\omega^{-,-})\otimes \tau^{I_1^c\cap I_2^c}_e,
\]
where now 
\[
\tau^I_e=\sum_{k\in I}1^{\otimes (i-1)}\otimes \iota_{\mathrm dx_k}\otimes 1^{\otimes (m-i)}\otimes 1^{\otimes (j-1)}\otimes \partial_{x_k}\otimes 1^{\otimes (m+1+n-j)}.
\]
The degree of the operator $\tau^I_e$ is readily computed to be $-1$, because of the contraction operators.

To a general admissible graph $\Gamma$ in $\mathcal G_{0,m+1+n}$, to $m$, resp.\ $n$, general elements $a_i$ of $A$, resp.\ $b_j$ of $B$, and $k$ of $K$, we associate an element of $K$ by
\[
\mathcal O_\Gamma^K(a_1|\cdots|a_m|k|b_1|\cdots|b_n)=\mu_{m+1+n}^K\left(\int_{\mathcal C_{0,m+1+n}^+}\prod_{e\in \mathrm E(\Gamma)}\omega^K_e(a_1|\cdots|a_m|k|b_1|\cdots|b_n)\right),
\] 
where $\mu^K_{m+1+n}:T_\mathrm{poly}(X)^{\otimes m+1+n}\to K$ is the $k$-multi-linear map given by multiple products in $T_\mathrm{poly}(X)$, followed by restriction on $K$.
Of course, we implicitly regard $A$, $B$ and $K$ as subalgebras of $T_\mathrm{poly}(X)$ w.r.t.\ the wedge product.

First of all, we observe that the product over all edges of $\Gamma$ does not depend on the ordering of the factors: namely, $\omega^K_e$ is a smooth $1$-form, but is also an endomorphism of $T_\mathrm{poly}(X)^{\otimes m+1+n}$ of degree $-1$, because of the contraction.
Furthermore, since $\omega^K_e$ is a smooth $1$-form on the compactified configuration space $\mathcal C_{0,m+1+n}^+$, the integral exists.

Finally, we define the Taylor component $\mathrm d_K^{m,n}:A[1]^{\otimes m}\otimes K[1]\otimes B[1]^{\otimes n}\to K[1]$ {\em via}
\begin{equation}\label{eq-A_inf-bimod}
\mathrm d_K^{m,n}(a_1|\cdots|a_m|k|b_1|\cdots|b_n)=\sum_{\Gamma\in\mathcal G_{0,m+1+n}}\mathcal O_\Gamma^K(a_1|\cdots|a_m|k|b_1|\cdots|b_n),\ a_i\in A,\ b_j\in B,\ k\in K.
\end{equation}
We first observe that the map~\eqref{eq-A_inf-bimod} has degree $1$: namely, for a general admissible graph $\Gamma$ of type $(0,m+1+n)$, the operator $\mathcal O_\Gamma(a_1|\cdots|a_m|k|b_1|\cdots|b_n)$ does not vanish, only if $|\mathrm E(\Gamma)|=m+n-1$, which is the dimension of $\mathcal C_{0,m+1+n}^+$.
Since to each edge is associated a contraction operator, which lowers degrees by $1$, it follows immediately that $\mathrm d_K^{m,n}$ has degree $1$: of course, if we omit the degree-shifting, the degree of $\mathrm d_K^{m,n}$ is equivalently $1-m-n$.

For later purposes, we also observe that $A$, $B$ and $K$ factor into a product of a symmetric algebra and an exterior algebra, and we focus our attention to the symmetric part: assuming the arguments are all homogeneous w.r.t.\ the grading on the symmetric algebra, we now want to determine the corresponding grading of the map~\eqref{eq-A_inf-bimod}.
For this purpose, we introduce the following notation: a general element $a$ of $A$ has degree $\mathrm{deg}(a)$ w.r.t.\ the symmetric part, and similarly for $b$ in $B$ and $k$ in $K$.
Again, for a general admissible graph $\Gamma$ of type $(0,m+1+n)$, $\mathcal O_\Gamma(a_1|\cdots|a_m|k|b_1|\cdots|b_n)$ does not vanish, only if $\Gamma$ has exactly $m+n-1$ edges, and, since to each edge is associated a derivative, it follows easily that the polynomial degree of $\mathcal O_\Gamma(a_1|\cdots|a_m|k|b_1|\cdots|b_n)$ equals
\[
\sum_{i=1}^m\mathrm{deg}(a_i)+\mathrm{deg}(k)+\sum_{j=1}^n \mathrm{deg}(b_j)-(m+n-1)=\sum_{i=1}^m\mathrm{deg}(a_i)+\mathrm{deg}(k)+\sum_{j=1}^n\mathrm{deg}(b_j)+1-m-n.
\]

Lemma~\ref{l-CF}, Subsubsection~\ref{sss-4-3-2}, implies that the operator $\omega^K_e$ is non-trivial, only if the edge $e$ is as in $a)$ and $e)$, in which cases we have
\[
\omega^K_e=\begin{cases}
\pi_e^*(\omega^{+,-})\otimes \tau^{I_1\cap I_2^c},& \text{$e$ as in $a)$}\\
\pi_e^*(\omega^{-,+})\otimes \tau^{I_1^c\cap I_2},& \text{$e$ as in $e)$,}\\
\end{cases}
\]
hence a general admissible graph of type $(0,m+1+n)$ appearing in Formula~\eqref{eq-A_inf-bimod} has the form
\bigskip
\begin{center}
\resizebox{0.5 \textwidth}{!}{\input{adm-bimod.pstex_t}}\\
\text{Figure 5 - An admissible graph of type $(0,6)$ appearing in $\mathrm d_K^{2,3}$ } \\
\end{center}
\bigskip
In view of Remark~\ref{r-admgr}, we observe that admissible graphs with multiple edges yield trivial contributions: namely, if any two distinct vertices (both necessarily of the second type) are connected by more than $1$ edge, the corresponding integral weight vanishes, since it contains the square of a $1$-form $\omega^{+,-}$ or $\omega^{-,+}$.
\begin{Prop}\label{p-A_inf-bimod}
For a field $k$ of characteristic $0$, containing $\mathbb R$ or $\mathbb C$, we consider $A$, $B$ and $K$ as above.

Then, the Taylor components~\eqref{eq-A_inf-bimod} endow $K$ with an $A_\infty$-$A$-$B$-bimodule structure, where $A$ and $B$ are viewed as GAs with their natural product, hence, in particular, $A$ and $B$ have a (trivial) $A_\infty$-algebra structure.
\end{Prop}
If we denote by $\mathrm d_A$, resp.\ $\mathrm d_B$ and $\mathrm d_K$ the $A_\infty$-structures on $A$, $B$ and $K$ respectively described in Proposition~\ref{p-A_inf-bimod}, then we may regard the formal sum $\gamma=\mathrm d_A+\mathrm d_B+\mathrm d_K$ as a MCE for the graded Lie algebra $\widehat{\mathrm{C}}^\bullet(\texttt{Cat}_\infty(A,B,K))$: thus, the triple 
$\left(\widehat{\mathrm{C}}^\bullet(\texttt{Cat}_\infty(A,B,K)),[\gamma,\bullet],[\bullet,\bullet]\right)$ defines a DG Lie algebra, where $[\bullet,\bullet]$ denotes the Gerstenhaber bracket. 
In the case when ${\rm d}_K$ has only finitely many non-trivial Taylor components\footnote{For example, this happens when $X=U\oplus V$.} one may instead consider the DG Lie algebra 
$\left(\mathrm{C}^\bullet(\texttt{Cat}_\infty(A,B,K)),[\gamma,\bullet],[\bullet,\bullet]\right)$. 
\begin{proof}[Proof of Proposition~\ref{p-A_inf-bimod}]
The Taylor components~\eqref{eq-A_inf-bimod} define an $A_\infty$-$A$-$B$-bimodule structure, if the following identities hold true: 
\begin{equation}\label{eq-str-A_inf-bi}
\begin{aligned}
&\sum_{j=1}^{m-1}(-1)^{j}\mathrm d_K^{m-1,n}(a_1|\cdots|a_{j-1}|a_ja_{j+1}|a_{j+2}|\cdots|a_m|k|b_1|\cdots|b_n)+\\
&\sum_{j=1}^{n-1}(-1)^{m+j+1}\mathrm d_K^{m,n-1}(a_1|\cdots|a_m|k|b_1|\cdots|b_{j-1}|b_jb_{j+1}|b_{j+2}|\cdots|b_n)+\\
&\sum_{i=0}^m\sum_{j=0}^{n}(-1)^{(m-i+1)(i+j)+(1-i-j)\sum_{k=1}^{m-i}|a_k|}\mathrm d_K^{m-i,n-j}(a_1|\cdots|a_{m-i}|\mathrm d_K^{i,j}(a_{m-i+1}|\cdots|a_{m}|k|b_1|\cdots|b_{j})|b_{j+1}|\cdots|b_n)=0.
\end{aligned}
\end{equation}
The proof of Identity~\eqref{eq-str-A_inf-bi} is based on Stokes' Theorem in the same spirit of the proof of the main result of~\cite{K}: namely, the quadratic relations in \eqref{eq-str-A_inf-bi} are equivalent to quadratic relations between the corresponding integral weights, recalling~\eqref{eq-A_inf-bimod}. 

For this purpose, we consider 
\begin{equation}\label{eq-stokess}
\sum_{\widetilde\Gamma\in\mathcal G_{0,m+1+n}}\int_{\mathcal C_{0,m+1+n}^+}\mathrm d\widetilde{\mathcal O}_{\widetilde\Gamma}(b_1|\cdots|b_m|k|a_1|\cdots|a_n)=\sum_i\sum_{\widetilde\Gamma\in\mathcal G_{0,m+1+n}}\int_{\partial_i \mathcal C_{0,m+1+n}^+}\widetilde{\mathcal O}_{\widetilde\Gamma}(b_1|\cdots|b_m|k|a_1|\cdots|a_n)=0,
\end{equation}
where the first summation in the second expression in~\eqref{eq-stokess} is over boundary strata of $\mathcal C^+_{0,m+1+n}$ of codimension $1$, and
\[
\begin{aligned}
&\widetilde{\mathcal O}_{\widetilde\Gamma}(a_1|\cdots|a_m|k|b_1|\cdots|b_n)=\mu_{m+1+n}^K\left(\prod_{e\in \mathrm V(\widetilde\Gamma)}\omega^K_e(a_1|\cdots|a_m|k|b_1|\cdots|b_n)\right)=\\
&=\mu_{m+1+n}^K\left(\omega^K_{\widetilde\Gamma}(a_1|\cdots|a_m|k|b_1|\cdots|b_n)\right)
\end{aligned}
\]
which is viewed as a smooth $K$-valued form on $\mathcal C_{0,m+1+n}^+$ of form degree equal to $\mathrm E(\widetilde\Gamma)$.  
We observe that, by construction, a contribution indexed by a graph $\widetilde{\Gamma}$ in $\mathcal G_{0,m+1+n}$ is non-trivial, only if $\mathrm E(\widetilde\Gamma)=m+n-2$. 

Boundary strata of $\mathcal C_{0,m+1+n}^+$ of codimension $1$ are all of type~\eqref{eq-upbound1}, Subsection~\ref{ss-3-1}, with no points in $\mathbb H$: furthermore, we distinguish three cases
\begin{itemize}
\item[$i)$] $\partial_{\emptyset,B}\mathcal C_{0,m+1+n}^+\cong \mathcal C_{0,B}^+\times \mathcal C_{0,[m+1+n]\smallsetminus\{B\}\sqcup \{*\}}$, where $B$ is an ordered subset of $[m]$ of consecutive elements;
\item[$ii)$] $\partial_{\emptyset,B}\mathcal C_{0,m+1+n}^+\cong \mathcal C_{0,B}^+\times \mathcal C_{0,[m+1+n]\smallsetminus\{B\}\sqcup \{*\}}$, where $B$ is an ordered subset of $\{m+1,\dots,n\}$ of consecutive elements; 
\item[$iii)$] $\partial_{\emptyset,B}\mathcal C_{0,m+1+n}^+\cong \mathcal C_{0,B}^+\times \mathcal C_{0,[m+1+n]\smallsetminus\{B\}\sqcup \{*\}}$, where $B$ is an ordered subset of $[m+1+n]$ of consecutive elements, containing $m+1$. 
\end{itemize}
We begin by considering a general boundary stratum of type $i)$: it corresponds to the situation, where $|B|$ consecutive points on $\mathbb R$, labelled by $B$, collapse to a single point on $\mathbb R$, which lies on the left of the special point labelled by $m+1$. 

Recalling Lemma~\ref{l-or}, Subsection~\ref{ss-4-2},~\eqref{eq-b1-or}, and Lemma~\ref{l-CF}, $ii)$, Subsubsection~\ref{sss-4-3-2}, we get
\begin{equation}\label{d21}
\underset{\partial_{\emptyset,B}\mathcal C_{0,m+1+n}^+}\int\omega_{\widetilde{\Gamma}}^K=(-1)^{j(|B|+1)+1}\left(\int_{\mathcal C^+_{0,B}}\omega_{\Gamma_B}^A\right)\left(\underset{\mathcal C^{+}_{0,[m+1+n]\smallsetminus\{B\}\sqcup \{*\}}}\int\omega_{\Gamma^B}^K\right), 
\end{equation}
where $\Gamma_B$, resp. $\Gamma^B$, is the subgraph of $\widetilde\Gamma$, whose edges have both endpoints belonging to $B$, resp.\ the graph obtained from $\widetilde\Gamma$ by collapsing $\Gamma_B$ to a single vertex; $j$ is the minimum of $B$.

The operator-valued form $\omega_{\Gamma_B}^A$ will be defined precisely later on, since, as we will soon see, we will not actually need its form for the present computations.
We recall namely the general form of an element of $\mathcal G_{0,m+1+n}$: in particular, since all vertices labelled by $B$ lie on the left of the vertex labelled by $m+1$, the degree of the form $\omega_{\Gamma_B}^A$ equals $0$, since the graph $\Gamma_B$ does not contain any edge, whence, by dimensional reasons, its weight does not vanish only if $|B|=2$, i.e.\ $B=\{j,j+1\}$, for $1\leq j\leq m-1$, and equals to $1$.
As a consequence, $\Gamma^B$ is an admissible graph in $\mathcal G_{0,m+n}$.

We do not get any further sign other than the sign in Identity~\eqref{d21} coming from the orientation, when moving a copy of the standard multiplication on $T_\mathrm{poly}(X)$ to act on the factors $a_j$, $a_{j+1}$, since the standard multiplication has degree $0$. 
Therefore, the sum in Identity~\eqref{eq-stokess} over boundary strata of codimension $1$ of type $i)$ gives exactly the first term on the left-hand side of Identity~\eqref{eq-str-A_inf-bi}.

Second, we consider a general boundary stratum of codimension $1$ of type $ii)$: it describes the situation, where $|B|$ consecutive points on $\mathbb R$, labelled by $B$, collapse to a single point of $\mathbb R$, which lies on the right of the special point labelled by $m+1$.

Once again, we recall the orientation formul\ae~\eqref{eq-b1-or} from Lemma~\ref{l-or}, Subsection~\ref{ss-4-2}, to find a factorization as~\eqref{d21}.
We may now repeat almost {\em verbatim} the arguments in the analysis of the previous case: namely, $|B|=2$, and the minimum $j$ of $B$ satisfies, by assumption, $m+1<j$, which we also re-write, by abuse of notation, as $m+1+j$, for $1\leq j\leq n-1$.
Thus, the sum in Identity~\eqref{eq-stokess} over boundary strata of codimension $1$ of type $ii)$ produces the second term on the left-hand side of Identity~\eqref{eq-str-A_inf-bi}. 

It remains to discuss boundary strata of type $iii)$: in this case, the situation describes the collapse of $|B|$ consecutive points on $\mathbb R$, labelled by $B$, among which is the special point labelled by $m+1$, to a single point on $\mathbb R$, which will become the new special point. 

Recalling the orientation formul\ae~\eqref{eq-b1-or} from Lemma~\ref{l-or}, Subsection~\ref{ss-4-2}, we find a factorization of the type~\eqref{d21}.

First of all, we observe that, in this case, the subgraph $\Gamma_B$ is disjoint from $\Gamma\smallsetminus \Gamma_B$: this follows immediately from Lemma~\ref{l-CF}, Subsubsection~\ref{sss-4-3-2},~\eqref{eq-gamma-delta} and~\eqref{eq-eps}, and from the discussion on the shape of admissible graphs appearing in Formula~\eqref{eq-A_inf-bimod} (in other words, there are no edges connecting $\Gamma_B$ with its complement $\Gamma^B\smallsetminus\Gamma_B$).
In particular, $\widetilde\Gamma$ factors out as $\widetilde\Gamma=\Gamma_B\sqcup \Gamma^B$, and $\Gamma_B$ and $\Gamma^B$ are both admissible.
We also observe that, in general, $|B|\geq 2$ in this case: namely, $\Gamma_B$ can be non-empty.

The orientation sign is $j(|B|+1)+1$, where $j$ is the minimum of $B$: since $1\leq j\leq m$, we may rewrite it as $m-i+1$, for $i=1,\dots,m$.
The maximum of $B$ is bigger or equal than $m+1$, hence we may write it as $j$, for $0\leq j\leq n$, shifting w.r.t.\ $m+1$.

Plus, we get an additional sign $(1-i-j)\left(\sum_{k=1}^{m-i}|a_k|\right)$, when moving $\int_{\mathcal C_{0,B}^+}\omega^K_{\Gamma_B}$ through $a_k$, $k=1,\dots,m-i$.

Finally, the fact that $\Gamma_B$ and $\Gamma^B$ are disjoint implies that we may safely restrict the product of the $B$-factors in $\int_{\mathcal C_{0,B}^+}\omega_{\Gamma_B}^K(a_{m-i+1}|\cdots|a_m|k|b_1|\cdots|b_j)$ to $K$, since no derivative acts on it and departs from it.
As a consequence, the sum in Identity~\eqref{eq-stokess} over boundary strata of codimension $1$ of type $iii)$ yields the third term on the left-hand side of Identity~\eqref{eq-str-A_inf-bi}. 

(We now observe that the signs coming from orientations in the previous calculations agree with the signs in Identity~\eqref{eq-str-A_inf-bi} up to an overall $-1$-sign, which is of no influence.)
\end{proof}

\begin{Rem}
Observe that ${\rm d}_K$ may have infinitely many non-trivial Taylor components. In this case we have to deal with the completed Hochschild cochain complex 
$\widehat{\rm C}^\bullet({\tt Cat}_\infty(A,B,K))$ if we want $[{\rm d}_K,\bullet]$ to be well-defined. Nevertheless, this will not be sufficient to 
make our main Theorem \ref{t-form-cat} work in this context. To do so we will have to consider the completed Hochschild complex 
$\widehat{\rm C}^\bullet({\tt Cat}_\infty(\widehat{A},\widehat{B},\widehat{K}))$ associated to completed (or formal) versions of $A$, $B$ and $K$: 
\begin{align*}
A&=\widehat{\mathrm S}(U^*)\otimes \wedge (X/U)=\widehat{\mathrm S}(U^*)\otimes \wedge (U^\perp\cap V)\otimes \wedge (U+V)^\perp,\\
B&=\widehat{\mathrm S}(V^*)\otimes \wedge (X/V)=\widehat{\mathrm S}(V^*)\otimes \wedge (U\cap V^\perp)\otimes \wedge (U+V)^\perp,\\
K&=\widehat{\mathrm S}((U\cap V)^*)\otimes \wedge (U+V)^\perp,
\end{align*}
\end{Rem}

\section{Formality for the Hochschild cochain complex of an $A_\infty$-category}\label{s-6}

In this Section we may assume that ${\rm d}_K$ has finitely many non-trivial taylor components. 

We consider the $A_\infty$-algebras $A$, $B$ and the $A_\infty$-$A$-$B$-bimodule $K$ from Subsection~\ref{ss-5-2}, to which we associate the $A_\infty$-category $\texttt{Cat}_\infty(A,B,K)$, and the corresponding 
Hochschild cochain complex $\mathrm{C}^\bullet(\texttt{Cat}_\infty(A,B,K))$: in particular, we are interested in the DG Lie algebra-structure on 
$\left(\mathrm{C}^\bullet(\texttt{Cat}_\infty(A,B,K)),[\mu,\bullet],[\bullet,\bullet]\right)$, where $\mu$ denotes the $A_\infty$-$A$-$B$-bimodule structure on $\texttt{Cat}_\infty(A,B,K)$. 

We construct an $L_\infty$-quasi-isomorphism $\mathcal U$ from the DG Lie algebra $(T_\mathrm{poly}(X),0,[\bullet,\bullet])$ to the DG Lie algebra 
$\left(\mathrm{C}^\bullet(\texttt{Cat}_\infty(A,B,K)),[\mu,\bullet],[\bullet,\bullet]\right)$.
The proof of the main result is divided into two parts: first, we construct explicitly $\mathcal U$, and we prove, by means of Stokes' Theorem, that $\mathcal U$ is an $L_\infty$-morphism, and second, 
we will prove that $\mathcal U$ is a quasi-isomorphism.
The proof of the second statement is a consequence of Keller's condition.

In the general situation (when ${\rm d}_K$ may not have finitely many non-trivial Taylor components), the construction below produces an $L_\infty$-quasi-isomorphism $\mathcal U$ from the DG Lie algebra 
$(\widehat{T}_\mathrm{poly}(X),0,[\bullet,\bullet])$ to the DG Lie algebra 
$\left(\widehat{\mathrm{C}}^\bullet(\texttt{Cat}_\infty(\widehat{A},\widehat{B},\widehat{K})),[\mu,\bullet],[\bullet,\bullet]\right)$. 
The only substantial modification in the proof is in the Koszul duality argument, that we will make explicit in both situations. 

\subsection{The explicit construction}\label{ss-6-1}

We now produce an explicit formula for the $L_\infty$-quasi-isomorphism $\mathcal U$: first of all, by the results of Section~\ref{s-5}, to construct an $L_\infty$-morphism from $T_\mathrm{poly}(X)$ to 
$\mathrm{C}^\bullet(\texttt{Cat}_\infty(A,B,K))$ is equivalent to constructing three distinct maps $\mathcal U_A$, $\mathcal U_B$ and $\mathcal U_K$, where
\[
\begin{aligned}
\mathcal U_A:T_\mathrm{poly}(X)&\to \mathrm{C}^\bullet(A,A),\ & \mathcal U_B:T_\mathrm{poly}(X)&\to \mathrm{C}^\bullet(B,B),\\
\mathcal U_K:T_\mathrm{poly}(X)&\to \mathrm{C}^\bullet(A,B,K).
\end{aligned}
\]
We fix an orthogonal decomposition~\eqref{eq-orth-split} of $X$ as in Subsection~\ref{ss-4-2}, and an adapted coordinate system $\{x_i\}$, in the sense of Subsection~\ref{ss-4-2}; we also recall from Subsubsections~\ref{sss-4-3-1} and~\ref{sss-4-3-2} the $2$-colored and the $4$-colored propagators.

To a pair of non-negative integers $(n,m)$, we associate the set $\mathcal G_{n,m}$ of admissible graphs of type $(n,m)$; further, we may write $(n,m)=(n,p+1+q)$, if $m\geq 1$, for some non-negative integers $p$, $q$.

To an admissible graph $\Gamma$ in $\mathcal G_{n,m}$ and general elements $\gamma_i$ of $T_\mathrm{poly}(X)$, $i=1,\dots,n$, general elements $a_j$ of $A$, $j=1,\dots,m$, we associate an element of $A$ by the assignment
\begin{equation}\label{eq-A-comp}
\mathcal O^A_\Gamma(\gamma_1|\cdots|\gamma_n|a_1|\cdots|a_m)=\mu_{n+m}^B\left(\int_{\mathcal C_{n,m}^+}\omega^A_\Gamma(\gamma_1|\cdots|\gamma_n|a_1|\cdots|a_m)\right),
\end{equation}
where $\mu_{n+m}^A$ is the multiplication operator from $T_\mathrm{poly}(X)^{\otimes n+m}$ to $T_\mathrm{poly}(X)$, followed by restriction to $A$, viewed (in a non-canonical way) as a sub-algebra of $T_\mathrm{poly}(X)$.
Further, the $\Omega^{|\mathrm E(\Gamma)|}(\mathcal C_{n,m}^+)$-valued endomorphism of $T_\mathrm{poly}(X)^{\otimes n+m}$ is defined as
\begin{equation}\label{eq-A-form}
\omega^A_\Gamma=\prod_{e\in\mathrm E(\Gamma)}\omega^A_e,\ \omega^A_e=\pi_e^*(\omega^+)\otimes\left(\tau^{I_1\cap I_2}_e+\tau^{I_1\cap I_2^c}_e\right)+\pi_e^*(\omega^-)\otimes\left(\tau^{I_1^c\cap I_2}_e+\tau^{I_1^c\cap I_2^c}_e\right),
\end{equation} 
$\pi_e$ being the natural projection from $\mathcal C_{n,m}^+$ onto $\mathcal C_{2,0}$ or its boundary strata of codimension $1$ (in fact, $\omega^+$ and $\omega^-$ vanish on all strata of codimension $2$ of $\mathcal C_{2,0}$, thanks to Lemma~\ref{l-angle}, Subsubsection~\ref{sss-4-3-1}), and the operator $\tau^I_e$, for $I\subset [d]$, has been defined in Subsection~\ref{ss-4-2}.

Once again, we observe that the product~\eqref{eq-A-form} is well-defined, since the $2$-colored propagators are $1$-forms, while $\tau^I_e$ is an endomorphism of $T_\mathrm{poly}(X)^{\otimes n+m}$ of degree $-1$.
Further, since the dimension of $\mathcal C_{n,m}^+$ is $2n+m-2$, the element~\eqref{eq-A-comp} is non-trivial, precisely when $|\mathrm E(\Gamma)|=2n+m-2$.

We then set 
\begin{equation}\label{eq-B-mor}
\mathcal U_A^n(\gamma_1|\cdots|\gamma_n)(a_1|\dots|a_m)=(-1)^{\left(\sum_{i=1}^n|\gamma_i|-1\right)m}\sum_{\Gamma\in\mathcal{G}_{n,m}}\mathcal{O}_{\Gamma}^A(\gamma_1|\cdots|\gamma_n|a_1|\cdots|a_m).
\end{equation}

Similar formul\ae, with due changes, specify the Taylor components $\mathcal U_B^n$, $n\geq 1$: we only observe that 
\[
\omega^B_e=\pi_e^*(\omega^+)\otimes\left(\tau_e^{I_1\cap I_2}+\tau_e^{I_1^c\cap I_2}\right)+\pi_e^*(\omega^-)\otimes\left(\tau_e^{I_1\cap I_2^c}+\tau_e^{I_1^c\cap I_2^c}\right),
\]
for an edge $e$ of a general admissible graph $\Gamma$ as above.

Finally, we define the Taylor components $\mathcal U_K^n$ {\em via}
\begin{equation}\label{eq-K-mor}
\mathcal U_K^n(\gamma_1|\cdots|\gamma_n)(a_1|\cdots|a_p|k|b_1|\cdots|b_q)=(-1)^{\left(\sum_{i=1}^n|\gamma_i|-1\right)(p+q+1)}\sum_{\Gamma\in\mathcal{G}_{n,p+1+q}}\mathcal{O}_{\Gamma}^K(\gamma_1|\cdots|\gamma_n|a_1|\cdots|a_p|k|b_1|\cdots|b_q).
\end{equation}

We want to point out now, before entering into the details, that $i)$ Formula~\eqref{eq-B-mor} contains admissible graphs with multiple edges and no loops (i.e.\ whenever an admissible graph contains at least $1$ loop, the corresponding contribution to Formula~\eqref{eq-B-mor} is set ot be $0$), and that $ii)$ Formula~\eqref{eq-K-mor} contains admissible graphs with multiple edges and (possibly) loops.

Since in the usual constructions in Deformation Quantization multiple edges and loops are not present, we need to discuss how to deal with both of them separately.

If $\Gamma$ is admissible and contains multiple edges, we consider a pair $(i,j)$ of distinct vertices of the first type of $\Gamma$, such that the cardinality of the set $\mathrm E_{(i,j)}=\{e\in \mathrm E(\Gamma):\ e=(i,j)\}$ is bigger than $1$.
Then, to $(i,j)$ we associate the smooth, operator-valued $|\mathrm E_{(i,j)}|$-form given by
\[
\omega_{(i,j)}^A=\frac{1}{(|\mathrm E_{(i,j)}|)!}\prod_{e\in \mathrm E_{(i,j)}}\omega_e^A=\frac{(\omega_{(i,j)}^A)^{|\mathrm E_{(i,j)}|}}{(|\mathrm E_{(i,j)}|)!},\ \omega_{(i,j)}^K=\frac{1}{(|\mathrm E_{(i,j)}|)!}\prod_{e\in \mathrm E_{(i,j)}}\omega_e^K=\frac{(\omega_{(i,j)}^K)^{|\mathrm E_{(i,j)}|}}{(|\mathrm E_{(i,j)}|)!}
\] 
(when replacing $A$ by $B$, obvious due changes have to be performed).

In particular, by abuse of notation, we denote by $\omega_e^A$, resp.\ $\omega_e^K$, the normalized operator-valued form associated to a (multiple) edge $e$ of $\Gamma$ in Formula~\eqref{eq-B-mor}, resp.~\eqref{eq-K-mor}: of course, if the edge $e$ appears only once in $\Gamma$, then $\omega_e^A$, resp.\ $\omega_e^K$, coincides with the standard expression, otherwise, it is given by the previous formula. 

We now recall from Subsubsection~\ref{sss-4-3-2} the closed $1$-form $\rho$ on $\mathcal C_{1,1}$.
The vertex $v_\ell$ of the first type, corresponding to a loop $\ell$ of $\Gamma$, specifies a natural projection $\pi_{v_\ell}:\mathcal C_{n,p+1+q}^+\to \mathcal C_{1,1}$, which extends to the corresponding compactified configuration spaces the projection onto the vertex $v_\ell$ and the special vertex $p+1$.
Further, we consider also the restricted divergence operator 
\[
\mathrm{div}^{(I_1\cap I_2)\sqcup (I_1^c\cap I_2^c)}=\sum_{k\in (I_1\cap I_2)\sqcup (I_1^c\cap I_2^c)}\iota_{\mathrm d x_k}\partial_{x_k}
\]
on $T_\mathrm{poly}(X)$; by $\mathrm{div}^{(I_1\cap I_2)\sqcup (I_1\cap I_2)}_{(r)}$, for $1\leq r\leq n$, we denote the endomorphism of $T_\mathrm{poly}(X)^{\otimes (n+p+1+q)}$ of degree $-1$ given by
\[
\mathrm{div}^{(I_1\cap I_2)\sqcup (I_1^c\cap I_2^c)}_{(r)}=1^{\otimes (r-1)}\otimes \mathrm{div}^{(I_1\cap I_2)\sqcup (I_1^c\cap I_2^c)}\otimes 1^{(n-r+p+1+q)}.
\]
Finally, for a loop $\ell$ of $\Gamma$, we set  
\begin{equation}
\omega_\ell=\pi_{v_\ell}^*(\rho)\otimes\mathrm{div}^{(I_1\cap I_2)\sqcup (I_1^c\cap I_2^c)}_{(v_\ell)}:
\end{equation}
it is clear that $\rho_\ell$ is a closed $1$-form on $\mathcal C_{n,p+1+q}^+$ with values in $\mathrm{End}(T_\mathrm{poly}(X)^{\otimes (n+p+1+q)})$ of degree $-1$, whence $\omega_\ell$ has total degree $-1$.

\begin{Rem}\label{r-loop}
We observe that loops are trivial, when $U\oplus V=X$, because the restricted divergence operator vanishes by construction.
\end{Rem}

We want to examine in some detail the admissible graphs and their colorings yielding (possibly) non-trivial contributions to Formul\ae~\eqref{eq-B-mor} and~\eqref{eq-K-mor}. 

We begin with Formula~\eqref{eq-B-mor}: in this case, we recall Lemma~\ref{l-angle}, Subsubsection~\ref{sss-4-3-1}, $ii)$, which implies that the $2$-colored propagator $\omega^+$, resp.\ $\omega^-$, vanishes on the boundary stratum $\beta$, resp.\ $\gamma$.
This, in turn, implies that edges of an admissible graph $\Gamma$ of type $(n,m)$, whose initial, resp.\ final, point lies in $\mathbb R$, are colored by propagators of type $\omega^-$, resp.\ $\omega^+$: according to the definition of $\omega^A_e$, for $e$ an edge of $\Gamma$, since to vertices of the second type are associated to elements of $A$, this is coherent with the fact that such elements may be differentiated only w.r.t.\ to coordinates $\{x_i\}$, for $i$ in $(I_1\cap I_2)\sqcup (I_1\cap I_2^c)$, and can be contracted only w.r.t.\ differentials of coordinates $\{x_i\}$, for $i$ in $(I_1^c\cap I_2)\sqcup (I_1^c\cap I_2^c)$.
Pictorially,
\bigskip
\begin{center}
\resizebox{0.5 \textwidth}{!}{\input{koszul-B.pstex_t}}\\
\text{Figure 6 - A general admissible graph of type $(4,4)$ appearing in $\mathcal U_A$} \\
\end{center}
\bigskip
Similar arguments hold, when replacing $A$ by $B$.

We now consider Formula~\eqref{eq-K-mor}, in particular, an admissible graph $\Gamma$ of type $(n,p+1+q)$.

The point $k+1$ on $\mathbb R$ plays a very special role in subsequent computations: in fact, it corresponds, w.r.t.\ the natural projections from $\mathcal C_{n,k+1+l}^+$ onto $\mathcal C_{2,1}$, to the single point on $\mathbb R$ in $\mathcal C_{2,1}$.
 
First of all, we recall Lemma~\ref{l-CF}, Subsubsection~\ref{sss-4-3-2}, $iii)$: as a consequence, if $e$ is an edge, whose initial, resp.\ final, point is $p+1$, then $e$ is colored by the propagator $\omega^{-,-}$, resp.\ $\omega^{+,+}$, and according to the definition of $\omega_e^K$, this is coherent with the fact that an element $k$ of $K$ can be only differentiated w.r.t.\ coordinates $\{x_i\}$, $i$ in $I_1\cap I_2$, and can be contracted only w.r.t.\ differentials of coordinates $\{x_i\}$, for $i$ in $I_1^c\cap I_2^c$.

As a consequence of the very same arguments of Subsection~\ref{ss-4-2}, $\Gamma$ cannot contain any edge $e$, which joins two vertices of the second type, both lying either on the left-hand side of $p+1$ or on the right-hand side of $p+1$; similarly, there is no edge joining $p+1$ to any other vertex of the second type.

It is also clear that, if $\Gamma$ possesses a vertex of the first type with more than $1$ loop attached to it, then the corresponding contribution to Formula~\eqref{eq-K-mor} vanishes, since it contains the square of the $1$-form $\rho$ on $\mathcal C_{1,1}$.

Finally, we observe that, if $\Gamma$ has more than $4$ multiple edges between the same two distinct vertices (obviously of the first type), the corresponding contribution to Formula~\eqref{eq-K-mor} is trivial: namely, since to any edge is associated a sum of $4$ distinct $1$-forms, any power of at least $5$ identical operator-valued forms contains at least a square of $1$ of the $4$-colored propagators.

Pictorially, 
\bigskip
\begin{center}
\resizebox{0.5 \textwidth}{!}{\input{koszul-K.pstex_t}}\\
\text{Figure 7 - A general admissible graph of type $(4,4)$ appearing in $\mathcal U_K$} \\
\end{center}
\bigskip
Of course, once again, we recall that loops do not appear in the special case $U\oplus V=X$.

\subsection{The main result}\label{ss-6-2}
We now state and prove the main result of the paper, namely
\begin{Thm}\label{t-form-cat}
We consider $X=k^d$, and we denote collectively by $\mu$ the $A_\infty$-structure on the category $\texttt{Cat}_\infty(A,B,K)$, defined as in Subsection~\ref{ss-5-2}. 

The morphisms $\mathcal U_A^n$, $\mathcal U_B^n$ and $\mathcal U_K^n$, $n\geq 1$, are the Taylor components of an $L_\infty$-quasi-isomorphism
\[
\mathcal U:\left(T_\mathrm{poly}(X),0,[\bullet,\bullet]\right)\to \left(\mathrm{C}^\bullet(\texttt{Cat}_\infty(A,B,K)),[\mu,\bullet],[\bullet,\bullet]\right).
\]
\end{Thm}
\begin{proof}
First of all, $T_\mathrm{poly}(X)$ and $\mathrm{C}^\bullet(\texttt{Cat}_\infty(A,B,K))$ are $L_\infty$-algebras {\em via} 
\begin{eqnarray}
&&Q_1=0,~~Q_{2}(\gamma_1,\gamma_2)=(-1)^{|\gamma_2|}[\gamma_1,\gamma_2], ~~\gamma_i\in (T_\mathrm{poly}(W)[1])_{|\gamma_i|},\ i=1,2 \nonumber \\
&&Q'_1=[\mu,\bullet],~~Q'_{2}(\phi_1,\phi_2)=(-1)^{|\phi_1|}[\phi_1,\phi_2], ~~\phi_i\in (\mathrm{C}^\bullet(\texttt{Cat}_\infty(A,B,K)))[1])_{|\phi_i|},\ i=1,2.\label{signs11}
\end{eqnarray}
For the sake of simplicity, we set $\mathcal U^n=\mathcal U_B^n+\mathcal U_K^n+\mathcal U_A^n$.

The conditions for $\mathcal U$ to be an $L_{\infty}$-morphism translate into the semi-infinite family of relations
\begin{eqnarray}
&&[\mu,\mathcal U^n(\gamma_1|\cdots|\gamma_n)]+\frac{1}{2}\sum_{I\sqcup J=\{1,\dots,n\}, I,J\neq \emptyset}\epsilon_{\gamma}(I,J)Q_2'(\mathcal U^{|I|}(\gamma_I),\mathcal U^{|J|}(\gamma_J))=\nonumber \\
&&=\frac{1}{2}\sum_{k\neq l}\sigma_{\gamma}(k,l,1,\dots,\hat{k},\dots,\hat{l},\dots,n)\mathcal U^{n-1}(Q_2(\gamma_k,\gamma_l),\gamma_1,\dots,\widehat{\gamma_k},\dots,\widehat{\gamma_l},\dots,\gamma_n).\label{morph2}
\end{eqnarray}
We denote by $\gamma_I$ the element $(\gamma_{i_1},\dots,\gamma_{i_I})\in C^{+|I|}(T_{poly}(W)[1])$, for every index set $I=\{i_1,\dots,i_I\}\subseteq\{1,\dots,n\}$ of cardinality $|I|$. The same notation holds for $\gamma_J$.

The infinite set of identities~\eqref{morph2} consists of three different infinite sets of identities, corresponding to the three projections of~\eqref{morph2} onto $A$, $B$ and $K$.
It is easy to verify that the projections onto $A$ or $B$ of~\eqref{morph2} define infinite sets of identities, which correspond to the identities satisfied by $L_\infty$-morphisms from $T_\mathrm{poly}(X)$ to $\mathrm{C}^\bullet(A,A)$ or $\mathrm{C}^\bullet(B,B)$, which have been proved in~\cite{CF} (in a slightly different form). 

Thus, it remains to prove Identity~\eqref{morph2} for the $K$-component.

First of all, we observe that
\[
[\mu,\mathcal U^n(\gamma_1|\cdots|\gamma_n)]=\mu\bullet \mathcal U^n(\gamma_1|\cdots|\gamma_n)-(-1)^{\sum_{i=1}^n|\gamma_i|+2-n}\mathcal U^n(\gamma_1|\cdots|\gamma_n)\bullet \mu.
\]
By setting $\mathcal U^0=\mu$, and recalling the higher compositions $\bullet$ from Subsection~\ref{ss-2-1}, the product $\bullet$ on $T_\mathrm{poly}(X)$, and by finally projecting down onto $K$ Identity~\eqref{morph2}, we find 
\begin{equation}\label{eq-L_inf-K}
\begin{aligned}
&\sum_{I\sqcup J=[n]}\epsilon_{\gamma}(I,J)\left(\mathcal U_K^{|I|}(\gamma_I)\bullet \mathcal U_B^{|J|}(\gamma_J)+\mathcal U_K^{|I|}(\gamma_I)\bullet \mathcal U_K^{|J|}(\gamma_J)+\mathcal U_K^{|I|}(\gamma_I)\bullet \mathcal U_A^{|J|}(\gamma_J)\right)=\\
&=\sum_{k\neq l}\sigma_{\gamma}(k,l,1,\dots,\hat{k},\dots,\hat{l},\dots,n)\mathcal U_K^{n-1}(\gamma_k\bullet\gamma_l,\gamma_1,\dots,\widehat{\gamma_k},\dots,\widehat{\gamma_l},\dots,\gamma_n).
\end{aligned}
\end{equation}

The proof of Identity~\eqref{eq-L_inf-K} relies on Stokes' Theorem: namely, for any two non-negative integers $p$, $q$, we consider the Identity for elements of $\mathrm{Hom}(T_\mathrm{poly}(X)^{\otimes (n+p+1+q)},K)$, 
\begin{equation}\label{eq-stok-L}
\sum_{\widetilde\Gamma\in\mathcal G_{n,p+1+q}}\int_{\mathcal C_{n,p+1+q}^+}\mathrm d\widetilde{\mathcal O}^K_{\widetilde\Gamma}=\sum_i\sum_{\widetilde\Gamma\in\mathcal G_{n,p+1+q}}\int_{\partial_i \mathcal C_{n,p+1+q}^+}\widetilde{\mathcal O}^K_{\widetilde\Gamma}=0,
\end{equation}
where the first summation in the second expression in~\eqref{eq-stok-L} is over boundary strata of $\mathcal C^+_{n,p+1+q}$ of codimension $1$, and
\[
\widetilde{\mathcal O}^K_{\widetilde\Gamma}=\mu_{n+p+1+q}^K\circ \prod_{e\in \mathrm V(\widetilde\Gamma)}\omega^K_e=\mu_{n+p+1+q}^K\circ \omega^K_{\widetilde\Gamma},
\]
regarded as a smooth $K$-valued form on $\mathcal C_{n,p+1+q}^+$ of form degree equal to $|\mathrm E(\widetilde\Gamma)|$.  
Then, by construction, a contribution indexed by a graph $\widetilde{\Gamma}$ in $\mathcal G_{n,p+1+q}$ is non-trivial, only if $|\mathrm E(\widetilde\Gamma)|=2n+p+q-2$. 

Boundary strata of $\mathcal C_{n,p+1+q}^+$ of codimension $1$ are either of type~\eqref{eq-upbound1} or~\eqref{eq-upbound2}, Subsection~\ref{ss-4-1}: 
\begin{itemize}
\item[$i)$] $\partial_{A}\mathcal C_{n,p+1+q}^+\cong \mathcal C_{A}\times \mathcal C^+_{[n]\smallsetminus A\sqcup\{*\},p+1+q}$, where $A$ is a subset of $[n]$ with $|A|\geq 2$;
\item[$ii_1)$] $\partial_{A_1,B_1}\mathcal C_{n,p+1+q}^+\cong \mathcal C_{A_1,B_1}^+\times \mathcal C_{[n]\smallsetminus A_1,[p+1+q]\smallsetminus B_1\sqcup \{*\}}^+$, where $A_1$ is a subset of $[n]$ with $|A_1|\geq 1$ and $B_1$ is an ordered subset of $[p]$ of consecutive elements with $|B_1|\geq 1$;
\item[$ii_2)$] $\partial_{A_1,B_1}\mathcal C_{n,p+1+q}^+\cong \mathcal C_{A_1,B_1}^+\times \mathcal C_{[n]\smallsetminus A_1,[p+1+q]\smallsetminus B_1\sqcup \{*\}}^+$, where $A_1$ is a subset of $[n]$ with $|A_1|\geq 1$ and $B_1$ is an ordered subset of $\{p+2,\dots,p+q+1\}$ of consecutive elements with $|B_1|\geq 1$;
\item[$ii_3)$] $\partial_{A_1,B_1}\mathcal C_{n,p+1+q}^+\cong \mathcal C_{A_1,B_1}^+\times \mathcal C_{[n]\smallsetminus A_1,[p+1+q]\smallsetminus B_1\sqcup \{*\}}^+$, where $A_1$ is a subset of $[n]$ with $|A_1|\geq 1$ and $B_1$ is an ordered subset of $[p+1+q]$ of consecutive elements with $|B_1|\geq 1$ and containing $p+1$.
\end{itemize}
We begin by considering a general boundary stratum of type $i)$: it corresponds to the situation, where points in $\mathbb H$, labelled by $A$, collapse to a single point again in $\mathbb H$. 

For a boundary stratum as in $i)$, we need Lemma~\ref{l-CF}, Subsubsection~\ref{sss-4-3-2}, $i)$, to find the following factorization, for a general admissible graph $\widetilde \Gamma$ of type $(n,p+1+q)$ as in Identity~\eqref{eq-stok-L}, recalling the orientations~\eqref{eq-b2-or} from Lemma~\ref{l-or}, Subsection~\ref{ss-4-2}:
\begin{equation}\label{eq-b-i}
\int_{\partial_A \mathcal C_{n,p+1+q}^+}\omega^K_{\widetilde \Gamma}=-\left(\int_{\mathcal C_A}\omega^K_{\Gamma_A}\right)\left(\underset{\mathcal C_{[n]\smallsetminus A\sqcup \{*\},p+1+q}^+}\int \omega^K_{\Gamma^A}\right),
\end{equation}
where $\Gamma_A$, resp.\ $\Gamma^A$, is the subgraph of $\widetilde\Gamma$, whose edges have both endpoints in $A$, resp.\ $\Gamma^A$ is the graph obtained by collapsing the subgraph $\Gamma_A$ to a point.

We now focus on the first factor on the right-hand side of Identity~\eqref{eq-b-i}.

Recalling Lemma~\ref{l-CF}, $i)$, from Subsubsection~\ref{sss-4-3-2}, the restriction to $\mathcal C_A$ of $\omega_e^K$, for $e$ a edge of the subgraph $\Gamma_A$ (not counted with multiplicities, in the case of a multiple edge), may be re-written as
\[
\omega_e^K\big\vert_{\mathcal C_A}=\left(\pi_e^*(\mathrm d\varphi)-\pi_{v_A}^*(\rho)\right)\otimes \tau_e^{[d]}+\pi_{v_A}^*(\rho)\otimes \tau_e^{(I_1\cap I_2)\sqcup (I_1^c\cap I_2^c)}=\widetilde\omega_e+\rho_{v_A,e},
\]
where $\pi_e$ is the (smooth extension to compactified configuration spaces of the) natural projection from $\mathcal C_A$ onto $\mathcal C_2$, and $\pi_{v_A}$ is the (smooth extension to compactified configuration spaces of the) natural projection from $\mathcal C_{[n]\smallsetminus A\sqcup \{v_A\},p+1+q}^+$ onto $\mathcal C_{1,1}$, and $v_A$ denotes the vertex corresponding to the collapse of the subgraph $\Gamma_A$.
Of course, when $U\oplus V=X$, the second term in the rightmost expression vanishes, see also Remark~\ref{r-loop}.

Therefore, we may re-write 
\begin{equation}\label{eq-b-i-left}
\int_{\mathcal C_A}\omega_{\Gamma_A}^K=\int_{\mathcal C_A}\prod_{e\in\mathrm E(\Gamma_A)}(\widetilde\omega_e+\rho_{v_A,e})\prod_{\ell\ \text{loop of $\Gamma_A$}}\omega_\ell,
\end{equation}
where, of course, the contributions to multiple edges are normalized as above.

We now first observe that the form-part of any loop contribution and of any operator-valued form $\rho_{v_A,e}$ is simply $\rho$ evaluated at the vertex corresponding to the collapse: hence, there can be at most $1$ such contribution, and in particular, if $\Gamma_A$ contains more than $1$ loop, the corresponding boundary contribution vanishes.

We first consider $\Gamma_A$ to be loop-free: because of the previous argument, we may re-write the right-hand side of~\eqref{eq-b-i-left} as 
\[
\int_{\mathcal C_A}\omega_{\Gamma_A}^K=\int_{\mathcal C_A}\prod_{e\in\mathrm E(\Gamma_A)}\widetilde\omega_e+\sum_{e\in\mathrm E(\Gamma_A)}\left(\int_{\mathcal C_A}\prod_{e'\neq e}\widetilde\omega_{e'}\right)\rho_{v_A,e}.
\]
The two integral contributions on the right-hand side vanish, if $|A|\geq 3$, either because of dimensional reasons or in virtue of Kontsevich's Lemma: therefore, we need only consider the case $|A|=2$.
The integral contributions are non-trivial in this case, only if the degree of the integrand equals $1$, which happens only $\Gamma_A$ has at most $2$ edges: graphically, we find the contributions 
\bigskip
\begin{center}
\resizebox{0.6 \textwidth}{!}{\input{b-i-contr.pstex_t}}\\
\text{Figure 8 - The four possible loop-free subgraphs $\Gamma_A$ yielding non-trivial boundary contributions of type $i)$} \\
\end{center}
\bigskip

The contribution from the first graph, in view of the previous expression, is given by
\[
\int_{\mathcal C_2}\omega^K_{\Gamma_A}=\left(\int_{S_1}\mathrm d\varphi\right)\otimes\tau_e^{[d]}=\tau_e^{[d]},
\] 
as the integral over $\mathcal C_2=S^1$ of the second term in $\widetilde{\omega}_e$ is basic w.r.t.\ the fiber integration.

Taking into account the fact that the second graph has $2$ multiple edges and recalling thus the normalization factor $2$, its contribution equals
\[
\int_{\mathcal C_2}\omega_{\Gamma_A}^K=\pi_{v_A}^*(\rho)\otimes \tau^{[d]}_e\tau^{(I_1\cap I_2)\sqcup (I_1^c\cap I_2^c)}_e,
\]
where $e=(i,j)$.
The very same computations yield for the fourth graph
\[
\int_{\mathcal C_2}\omega_{\Gamma_A}^K=\pi_{v_A}^*(\rho)\otimes \tau^{[d]}_e\tau^{(I_1\cap I_2)\sqcup (I_1^c\cap I_2^c)}_e,
\] 
$e=(j,i)$ in this case.

Finally, the third graph yields the contribution
\[
\int_{\mathcal C_2}\omega_{\Gamma_A}^K=\pi_{v_A}^*(\rho)\otimes\tau^{[d]}_{e_1}\tau^{(I_1\cap I_2)\sqcup (I_1^c\cap I_2^c)}_{e_2}+\pi_{v_A}^*(\rho)\otimes\tau^{[d]}_{e_2}\tau^{(I_1\cap I_2)\sqcup (I_1^c\cap I_2^c)}_{e_1}, 
\] 
where $e_1=(i,j)$, $e_2=(j,i)$.

Now, we assume the subgraph $\Gamma_A$ to have exactly one loop: in this case, the right-hand side of~\eqref{eq-b-i-left} can be re-written as
\[
\int_{\mathcal C_A}\omega_{\Gamma_A}^K=\left(\int_{\mathcal C_A} \prod_{e\in \mathrm E(\Gamma_A)}\widetilde\omega_e\right)\omega_\ell,
\] 
because the $1$-form associated to the loop is basic w.r.t.\ the projection onto $\mathcal C_A$.
Again, dimensional reasons or Kontsevich's Lemma imply that the above contribution is non-trivial, only if $|A|=2$: in this case, the subgraph $\Gamma_A$ yields non-trivial contributions, only if it is as in the picture
\bigskip
\begin{center}
\resizebox{0.6 \textwidth}{!}{\input{b-i-loop.pstex_t}}\\
\text{Figure 9 - The four possible subgraphs $\Gamma_A$ with one loop yielding non-trivial boundary contributions of type $i)$} \\
\end{center}
\bigskip

We write down explicitly only the contribution coming from the first graph
\[
\int_{\mathcal C_2}\omega_{\Gamma_A}^K=\pi_{v_A}^*(\rho)\otimes\mathrm{div}_{(v_A)}^{(I_1\cap I_2)\sqcup (I_1^c\cap I_2^c)}\tau_e^{[d]}, 
\]
where $e=(i,j)$, and, by the very construction of $\omega_\ell$, $v_\ell=v_A$.

We now recall the sign conventions previously discussed, which imply that sign issues can be dealt in this framework exactly as in the proof of Theorem A.7,~\cite{CF}.
We only observe that $i)$ the endomorphism $\tau_e^{[d]}$, which appears in all contributions, leads to the Schouten--Nijenhuis bracket between the poly-vector fields associated to the two distinct vertices of $\Gamma_A$, and that $ii)$ the contributions involving the restricted divergence and the endomorphism $\tau_e^{(I_1^c\cap I_2)\sqcup (I_1\cap I_2^c)}$ sum up, by Leibniz's rule, to the restricted divergence applied to the Schouten--Nijenhuis bracket between the aforementioned poly-vector fields.  

Thus, the sum in~\eqref{eq-stok-L} involving boundary strata of type $i)$ contribute to the right-hand side of Identity~\eqref{eq-L_inf-K} in all cases.

We then consider boundary strata of type $ii_1)$: such strata describe the collapse of points in $\mathbb H$ labelled by $A_1$ and of consecutive points on $\mathbb R$ labelled by $B_1$, where the maximum of $B_1$ lies on the left-hand side of the special point labelled by $p+1$, to a single point in $\mathbb R$ (the point resulting from the collapse lies obviously on the left-hand side of $p+1$), graphically
\bigskip
\begin{center}
\resizebox{0.4 \textwidth}{!}{\input{koszul-b-ii_1.pstex_t}}\\
\text{Figure 10 - A general configuration of points in a boundary stratum of $\mathcal C^+_{n,p+1+q}$ of type $ii_1)$} \\
\end{center}
\bigskip
We recall, in particular, Lemma~\ref{l-CF}, Subsubsection~\ref{sss-4-3-2}, $ii)$, for the restriction of the $4$-colored propagators on the boundary stratum $\beta$ of $\mathcal C_{2,1}$, and the orientations~\ref{eq-b1-or} from Lemma~\ref{l-or}, Subsection~\ref{ss-4-2}: hence, we get the factorization
\begin{equation}\label{eq-b-ii_1}
\int_{\partial_{A_1,B_1} \mathcal C_{n,p+1+q}^+}\omega^K_{\widetilde \Gamma}=(-1)^{j(|B_1|+1)+1}\left(\int_{\mathcal C_{A_1,B_1}^+}\omega^A_{\Gamma_{A_1,B_1}}\right)\left(\underset{\mathcal C_{[n]\smallsetminus A_1,[p+1+q]\smallsetminus B_1\sqcup\{*\}}^+}\int \omega^K_{\Gamma^{A_1,B_1}}\right),
\end{equation} 
where $\Gamma_{A_1,B_1}$, resp.\ $\Gamma^{A_1,B_1}$, denotes the subgraph of $\widetilde\Gamma$, whose edges have both endpoints labelled by $A_1\sqcup B_1$, resp.\ the graph obtained by collapsing $\Gamma_{A,B}$ to a single point.

We observe first that $\Gamma_{A_1,B_1}$ cannot have edges connecting vertices labelled by $A_1\sqcup B_1$ to vertices on $\mathbb R$, not labelled by $A_1\sqcup B_1$, which lie on the left of the vertex labelled by $p$, because of Lemma~\ref{l-CF}, Subsubsection~\ref{sss-4-3-2}, $iv)$. 
It thus follows that $\Gamma_{A_1,B_1}$, as well as $\Gamma^{A_1,B_1}$, is an admissible graph.

Second, we notice that, if $\Gamma_{A_1,B_1}$ has at least one loop, the corresponding contribution vanishes, because the $1$-form $\rho$ vanishes on the boundary strata of codimension $1$ of $\mathcal C_{1,1}$.

Once again, the sign conventions for the higher compositions $\bullet$ we have previously elucidated, see Subsection~\ref{ss-2-1}, imply that all signs arising in this situation are the same signs, with due modifications, owing to the different algebraic setting, appearing in the proof of Theorem A.7,~\cite{CF}: due to the appearance of operators of the form $\omega_{\Gamma_{A_1,B_1}^B}$ in Identity~\eqref{eq-b-ii_1}, it follows that the sum in~\eqref{eq-stok-L} over all boundary strata of type $ii_1)$ yields the first term on the left-hand side of Identity~\eqref{eq-L_inf-K}.

Now, we consider boundary strata of type $ii_2)$: in this case, such a boundary stratum describes the collapse of points in $\mathbb H$, labelled by $A_1$, and of consecutive points on $\mathbb R$, labelled by $B_1$, where the minimum of $B_1$ lies on the right-hand side of $p+1$, to a single point on $\mathbb R$ (clearly, the point resulting from the collapse lies on the right-hand side of $p+1$).
\bigskip
\begin{center}
\resizebox{0.55 \textwidth}{!}{\input{koszul-b-ii_2.pstex_t}}\\
\text{Figure 11 - A general configuration of points in a boundary stratum of $\mathcal C^+_{n,p+1+q}$ of type $ii_2)$} \\
\end{center}
\bigskip
In this situation, we recall Lemma~\ref{l-CF}, Subsubsection~\ref{sss-4-3-2}, $ii)$, when dealing with the restriction of the $4$-colored propagators on the boundary stratum $\gamma$ of $\mathcal C_{2,1}$, and, once again, the orientations~\ref{eq-b1-or} from Lemma~\ref{l-or}, Subsection~\ref{ss-4-2}, whence comes the factorization
\begin{equation}\label{eq-b-ii_2}
\int_{\partial_{A_1,B_1} \mathcal C_{n,p+1+q}^+}\omega^K_{\widetilde \Gamma}=(-1)^{j(|B_1|+1)+1}\left(\int_{\mathcal C_{A_1,B_1}^+}\omega^B_{\Gamma_{A_1,B_1}}\right)\left(\underset{\mathcal C_{[n]\smallsetminus A_1,[p+1+q]\smallsetminus B_1\sqcup\{*\}}^+}\int \omega^K_{\Gamma^{A_1,B_1}}\right),
\end{equation} 
with the same notation as in Identity~\eqref{eq-b-ii_2}.

Once again, because of Lemma~\ref{l-CF}, Subsubsection~\ref{sss-4-3-2}, $iv)$, the subgraph $\Gamma_{A_1,B_1}$ cannot have edges connecting vertices of $\Gamma_{A_1,B_1}$ to vertices on $\mathbb R$ on the right of $p$, hence $\Gamma_{A_1,B_1}$ and $\Gamma^{A_1,B_1}$ are both admissible graphs.

As already noticed for boundary strata of type $ii_1)$, if the subgraph $\Gamma_{A_1,B_1}$ contains at least one loop, the corresponding contribution vanishes, by the very same arguments as above.

Needless to repeat, the sign conventions for the corresponding higher compositions $\bullet$ from Subsection~\ref{ss-2-1} imply that all signs arising in this situation tantamount to the signs (with obvious due modifications) from the proof of Theorem A.7,~\cite{CF}: because of the presence of the form $\omega_{\Gamma_{A_1,B_1}^A}$ in Identity~\eqref{eq-b-ii_1}, the sum in~\eqref{eq-stok-L} over all boundary strata of type $ii_2)$ yields the third term on the left-hand side of Identity~\eqref{eq-L_inf-K}.

Finally, we consider boundary strata of type $ii_3)$: a stratum of this type describes the collapse of points in $\mathbb H$, labelled by $A_1$, and of points on $\mathbb R$, labelled by $B_1$ (which, this time, contains the special point $p+1$), to a single point in $\mathbb R$, which becomes the new special point.
\bigskip
\begin{center}
\resizebox{0.45 \textwidth}{!}{\input{koszul-b-ii_3.pstex_t}}\\
\text{Figure 12- A general configuration of points in a boundary stratum of $\mathcal C^+_{n,p+1+q}$ of type $ii_3)$} \\
\end{center}
\bigskip
We make use of Lemma~\ref{l-CF}, Subsubsection~\ref{sss-4-3-2}, $iii)$, for the restriction of the $4$-colored propagators on the boundary strata $\delta$ and $\varepsilon$ of $\mathcal C_{2,1}$, and of the orientations~\ref{eq-b1-or} from Lemma~\ref{l-or}, Subsection~\ref{ss-4-2} to come to the factorization
\begin{equation}\label{eq-b-ii_3}
\int_{\partial_{A_1,B_1} \mathcal C_{n,p+1+q}^+}\omega^K_{\widetilde \Gamma}=(-1)^{j(|B_1|+1)+1}\left(\int_{\mathcal C_{A_1,B_1}^+}\omega^K_{\Gamma_{A_1,B_1}}\right)\left(\underset{\mathcal C_{[n]\smallsetminus A_1,[p+1+q]\smallsetminus B_1\sqcup\{*\}}^+}\int \omega^K_{\Gamma^{A_1,B_1}}\right),
\end{equation} 
where we have used the same notation as in~\eqref{eq-b-ii_1} and~\eqref{eq-b-ii_3}.

We observe that, in this case, $\Gamma_{A_1,B_1}$ cannot have, once again, edges connecting vertices of $\Gamma_{A_1,B_1}$ to vertices on $\mathbb R$, because of Lemma~\ref{l-CF}, Subsubsection~\ref{sss-4-3-2}, $iv)$; further, the only incoming, resp.\ outgoing, edges of $\Gamma_{A_1,B_1}$ are labelled by propagators of the form $\omega^{+,+}$, resp.\ $\omega^{-,-}$, because of Lemma~\ref{l-CF}, Subsubsection~\ref{sss-4-3-2}, $iii)$. 
In particular, $\Gamma_{A_1,B_1}$, as well as $\Gamma^{A_1,B_1}$, is an admissible graph.

Thanks to the previously discussed sign conventions for the higher compositions $\bullet$ in see Subsection~\ref{ss-2-1}, all signs arising in this situation are the same appearing in the proof of Theorem A.7,~\cite{CF}: because of operators of the form $\omega_{\Gamma_{A_1,B_1}^K}$ in Identity~\eqref{eq-b-ii_1}, the sum in~\eqref{eq-stok-L} over all boundary strata of type $ii_1)$ yields the second term on the left-hand side of Identity~\eqref{eq-L_inf-K}.
\end{proof}

\subsection{The $L_\infty$-morphism $\mathcal U$ is an $L_\infty$-quasi-isomorphism: the Hochschild--Kostant--Rosenberg quasi-isomorphism for $\texttt{Cat}_\infty(A,B,K)$}\label{ss-6-3}
So far, we have only proved that the morphism constructed in Subsection~\ref{ss-6-1} is an $L_\infty$-morphism: it remains to prove that $\mathcal U$ is in fact an $L_\infty$-quasi-isomorphism: equivalently, we have to prove that its first Taylor component $\mathcal U_1$ is a quasi-isomorphism.

We observe now that the $L_\infty$-morphism $\mathcal U$ fits into the following commutative diagram of $L_\infty$-algebras:
\begin{equation}\label{d-keller}
\xymatrix{ & \left(\mathrm C^\bullet(A,A),[\mathrm d_A,\bullet],[\bullet,\bullet]\right) & \\
\left(T_\mathrm{poly}(X),0,[\bullet,\bullet]\right)\ar[ur]^{\mathcal U_A} \ar[dr]_{\mathcal U_B} \ar[rr]^{\mathcal U}& & \left(\mathrm{C}^\bullet(\texttt{Cat}_\infty(A,B,K)),[\mu,\bullet],[\bullet,\bullet]\right)\ar[ul]_{\mathrm p_A} \ar[dl]^{\mathrm p_B}\\
 & \left(\mathrm C^\bullet(B,B),[\mathrm d_B,\bullet],[\bullet,\bullet]\right) &}
\end{equation}

The relative Formality Theorem of~\cite{CF} implies that $\mathcal U_A^1$ and $\mathcal U_B^1$ are $L_\infty$-quasi-isomorphisms.
Hence, if we can prove that the projections $\mathrm p_A$ and $\mathrm p_B$ are quasi-isomorphisms (in particular, $L_\infty$-quasi-isomorphisms), the invertibility property of $L_\infty$-quasi-isomorphisms would imply that also $\mathcal U$ is a quasi-isomorphism.

By Theorem~\ref{t-keller}, Subsection~\ref{ss-3-3}, it suffices to prove that the left derived action $\mathrm L_A$ and the right derived action $\mathrm R_B$ are quasi-isomorphisms.

We will prove that the left derived action $\mathrm L_A$ is a quasi-isomorphism; the proof for $\mathrm R_B$ follows by the same arguments (with due modifications).

\subsubsection{$\mathrm S(Y^*)$ as a (relative) quadratic algebra}\label{sss-6-3-1}
We consider, more generally, a finite-dimensional graded $k$-vector space $Y$ with a fixed direct sum decomposition $Y=X_1\oplus X_2$ into (finite-dimensional) graded subspaces $X_1$, $X_2$.

We further consider the symmetric algebra $\mathrm S(Y^*)$: owing to the decomposition $Y=X_1\oplus X_2$, $\mathrm S(Y^*)\cong \mathrm S(X_1^*)\otimes \mathrm S(X_2^*)$, whence $\mathrm S(Y^*)$ has a structure of left $\mathrm S(X_1^*)$-module.
Conversely, $\mathrm S(X_1^*)$ has a structure of left $\mathrm S(Y^*)$-module, w.r.t.\ the natural projection from $\mathrm S(Y^*)$ onto $\mathrm S(X_1^*)$.

We now set, for the sake of simplicity, $A_0=\mathrm S(X_1^*)$ and $A_1=\mathrm S(X_1^*)\otimes X_2^*$: $A_1$ is a free $A_0$-module in a natural way.
We further have the obvious identification $\mathrm T_{A_0}A_1\cong \mathrm S(X_1^*)\otimes \mathrm T(X_2^*)$, where $\mathrm T_{A_0}(A_1)$ denotes the tensor algebra over $A_0$ of $A_1$, and similarly for the tensor algebra $\mathrm T(X_2^*)$ over $\mathbb C$.
Further, we consider 
\begin{equation*}
R=\left\{1\otimes v_1^*\otimes v_2^*-(-1)^{|v_1^*||v_2^*|}1\otimes v_2^*\otimes v_1^*:\ v_i^*\in X_2^*,\ i=1,2 \right\}\subset \mathrm S(X_1^*)\otimes (X_2^*)^{\otimes 2}\cong A_1\otimes_{A_0} A_1,
\end{equation*}
and, by abuse of notation, we denote by $R$ also the two-sided ideal in $\mathrm T_{A_0}(A_1)$ spanned by $R$.

It is then quite easy to verify that 
\[
\mathrm T_{A_0}(A_1)/R\cong \mathrm S(X_1^*)\otimes \mathrm S(X_2^*)\cong \mathrm S(Y^*),
\]
whence it follows that $A=\mathrm S(Y^*)$ is a quadratic $A_0$-algebra.

The algebra $A^!$, the quadratic dual of $A$, can be also computed explicitly: since $A^!=\mathrm T_{A_0}(A_1^\vee)/R^\perp$, where $A_1^\vee$ is the dual (over $A_0$) of $A_1$ and $R^\perp$ is the (two-sided ideal in $T_{A_0}(A_1^\vee)$ generated by the) annihilator of $R$ in $A_1^\vee\otimes_{A_0} A_1^\vee$, and since $Y$ is finite-dimensional, we have 
\[
A^!\cong \mathrm S(X_1^*)\otimes \Lambda(X_2)\cong \mathrm S(X_1^*)\otimes \mathrm S(X_2[-1])\cong \mathrm S(X_1^*\oplus X_2[-1]),
\]
where the exterior algebra $\Lambda(X_2)$ of $X_2$ is defined by mimicking the standard definition in the category $\texttt{Mod}_k$, and where the second isomorphism is explicitly defined by the so-called {\em d\'ecalage} isomorphism.

Finally, by means of $A$ and $A^!$, we may compute the Koszul complex of $A$: since $\mathrm K_n(A)=A\otimes_{A_0} (A^!_n)^\vee$, where again $(A^!_n)^\vee$ denotes the dual over $A_0$ of $A^!_n$, we obtain
\[
\mathrm K^\bullet (A)\cong \mathrm S(Y^*)\otimes \mathrm S(X_2^*[1])\cong \mathrm S(Y^*\oplus X_2^*[1]),
\]
with the natural formula for the Koszul differential.

\medskip

We can repeat the very same construction with completed tensor and symmetric algebras\footnote{The completion considered here is the $\mathcal I$-adic completion, where $\mathcal I$ is the ideal generated by degree $0$ generators.}, 
yielding $\widehat{A}=\widehat{\rm S}(Y^*)$, $\widehat{A}^!=\widehat{\mathrm S}(X_1^*\oplus X_2[-1])$, and the completed Koszul complex $\widehat{\mathrm K}^\bullet (A)\cong \widehat{\mathrm S}(Y^*\oplus X_2^*[1])$. 

\subsubsection{The Koszul complex of $\mathrm S(Y^*)$}\label{sss-6-3-2}
We now inspect more carefully the Koszul complex $\mathrm K^\bullet(A)$ (viewed as a cohomological complex) of the algebra $A=\mathrm S(Y^*)$.

First of all, we discuss the gradings of $\mathrm K^\bullet(A)$.
The shift by $1$ of the grading of $X_2^*$ induces the {\bf cohomological grading}, which is concentrated in $\mathbb Z_{\leq 0}$.
Alternatively, we may view (the graded vector space of the complex) $\mathrm K^\bullet(A)$ as the (graded vector space of the) relative de Rham complex of $Y$ w.r.t.\ $X_2$, and the cohomological grading is the opposite of the natural grading of the relative de Rham complex as a complex.

Then, $\mathrm K^{-n}(A)$, for $n\geq 0$, is naturally an object of $\texttt{GrMod}_k$, and the corresponding grading is called total grading: furthermore, the total grading can be written as the sum of the cohomological grading and the internal grading.

Exemplarily, $\mathrm K^\bullet(A)$ is generated by $x_i$, $y_j$, $\theta_k$, where $\theta_k$ denotes a basis of $X_2^*[1]$ associated to a basis $y_j$ of $X_2^*$: by, definition, $|\theta_j|=|y_j|-1$.
Thus, a general element $x_{i_1}\cdots x_{i_p}y_{j_1}\cdots y_{j_q}\theta_{k_1}\cdots \theta_{k_r}$ of $\mathrm K^\bullet(A)$ has total, resp.\ cohomological, resp.\ internal, degree 
\[
\sum_{s=1}^p|x_{i_s}|+\sum_{t=1}^q|y_{j_t}|+\sum_{u=1}^r|\theta_{k_u}|-r,\ \text{resp.}\ -r,\ \text{resp.}\ \sum_{s=1}^p|x_{i_s}|+\sum_{t=1}^q|y_{j_t}|+\sum_{u=1}^r|\theta_{k_u}|.
\] 
The Koszul differential $\mathrm d$ is defined w.r.t.\ the previous basis as $\mathrm d=y_j\partial_{\theta_j}$, where $\partial_{\theta_j}$ denotes the derivation w.r.t.\ $\theta_j$ acting from the left with total degree $1$, cohomological degree $1$ and internal degree $0$. 

The Koszul complex $\mathrm K^\bullet(A)$ is endowed with a distinct differential $\mathrm d_\mathrm{dR}=\theta_j\partial_{y_j}$, where the differential $\partial_{y_j}$ acts from the left, with total degree $-1$, cohomological degree $-1$ and internal degree $0$.

The operator $\mathrm L_\mathrm{rel}=[\mathrm d_\mathrm{dR},\mathrm d]$, $[\ ,\ ]$ being the commutator in $\mathrm{End}(\mathrm K^\bullet(A))$ w.r.t.\ internal degree, has total degree $0$, cohomological degree $0$ and internal degree $0$: $\mathrm L_\mathrm{rel}$ is expressed on generators {\em via} $\mathrm L_\mathrm{rel}(x_i)=0$, $\mathrm L_\mathrm{rel}(y_j)=y_j$ and $\mathrm L_\mathrm{rel}(\theta_j)=\theta_j$, and is extended on general elements w.r.t.\ the Leibniz rule.

The homotopy formula $\mathrm L_\mathrm{rel}=[\mathrm d_\mathrm{dR},\mathrm d]$ implies by a direct computation that the Koszul complex of the quadratic algebra $A=\mathrm S(Y^*)$ is a resolution of $A_0=\mathrm S(X_1^*)$ as a left module over $A$, therefore $A$ and $A^!=\mathrm S(X_1^*\oplus X_2[-1])$ are quadratic Koszul algebras over $A_0$.

\medskip

The very same homotopy formula holds for the completed Koszul complex $\widehat{\mathrm K}^\bullet (A)$, which is therefore a resolution of $\widehat{A_0}=\widehat{\mathrm S}(X_1^*)$ as a left $\widehat{A}$-module. 

\subsubsection{Relative Koszul duality}\label{sss-6-3-3}
We consider the category ${}_{A_0}\texttt{Mod}$, for $A_0$ and $A$ as before, of graded left $A_0$-modules, with spaces of morphisms specified {\em via} $\mathrm{Hom}_{A_0-}(V,W)=\bigoplus_{n\in\mathbb Z}\mathrm{Hom}_{A_0-}^n(V,W)$, $n$ referring to the degree: from the arguments of Subsubsections~\ref{sss-6-3-1} and~\ref{sss-6-3-2}, we know that $A$ is a Koszul quadratic algebra, which is additionally commutative (in the graded sense).
We also recall that $A$ is bigraded w.r.t.\ the Koszul grading and w.r.t.\ the internal grading; similarly, the Koszul resolution $\mathrm K^\bullet(A)$ is bigraded w.r.t.\ the Koszul grading and internal grading.
The Koszul differential is compatible with the internal grading (i.e.\ it has internal degree $0$).

On the other hand, $A_0$, $A$ and $\mathrm K^{-n}(A)$ are all objects of the category $\texttt{Mod}_{k}$ w.r.t.\ the internal grading: since the Koszul differential has degree $0$, it makes sense to define
\begin{equation}\label{eq-ext-big}
\mathrm{Ext}_{A-}^n(A_0,A_0)=\bigoplus_{p+q=n}\mathrm{Ext}_{A-}^{(p,q)}(A_0,A_0)=\bigoplus_{p+q=n}\mathrm{ext}_{A-}^p(A_0,A_0[q]),
\end{equation}
where $p$, resp.\ $q$, corresponds to the Koszul, resp.\ internal, grading.
By $\mathrm{ext}_{A-}^\bullet(\bullet,A_0[q])$, we denote the right derived functor of $\mathrm{hom}_{A-}(\bullet,A_0[q])$ in the category ${}_{A}\texttt{grmod}$ of graded left $A_0$-modules, whose spaces of morphisms are defined as $\mathrm{hom}_{A-}(V,W)$, the space of morphisms of degree $0$ from $V$ to $W$.

The right-derived functor $\mathrm{ext}_{A-}^\bullet(\bullet,A_0[q])$ can be computed by means of the Koszul resolution $\mathrm K^\bullet(A)$, hence 
\[
\mathrm{ext}^p_{A-}(A_0,A_0[q])=\mathrm H^p(\mathrm{hom}_{A-}(\mathrm K^{-\bullet}(A),A_0[q]),\mathrm d),
\]
where, by abuse of notation, $\mathrm d$ denotes the differential induced by the Koszul differential $\mathrm d$ by composition on the right.
The Koszul differential $\mathrm d$ acts trivially, whence the cohomology of the previous complex identifies with the complex itself:
\[
\mathrm{ext}^p_{A-}(A_0,A_0[q])=\mathrm{hom}_{A-}(\mathrm K^{-\bullet}(A),A_0[q])\cong (A^!_p)_q,
\]
where the index $q$ on the right hand-side term refers to the internal grading.
The previous chain of isomorphisms follows from the fact that, being $V$ non-negatively graded w.r.t.\ the internal grading, $A^!_p$ and $A_p$ are both free of finite rank over $A_0$, thus, the dual space over $A_0$ of $A^!_p$ is naturally graded w.r.t.\ the internal grading.

$\mathrm{Ext}_{A-}$-groups admit the {\bf Yoneda product}, i.e.\ a pairing of Koszul and internal degree $0$
\[
\mathrm{Ext}^{(m_1,n_1)}_{A-}(A_0,A_0)\otimes \mathrm{Ext}^{(m_2,n_2)}_{A-}(A_0,A_0)\to \mathrm{Ext}^{(m_1+m_2,n_1+n_2)}_{A-}(A_0,A_0).
\]
We consider a representative $\alpha$ of an element of $\mathrm{Ext}^{(m_1,n_1)}_{A-}(A_0,A_0)\cong (A^!_{m_1})_{n_1}$.
More explicitly, $\alpha$ acts by multiplication w.r.t.\ $A_0$ and by derivations on $\mathrm S^{m_1}(X_2^*[1])$, finally setting coordinates on $X_2^*$ to $0$.
Furthermore, $\alpha$ can be lifted to an element $\alpha_n$ of $\mathrm{hom}_{A-}(\mathrm K^{-m_1-n}(A),\mathrm K^{-n}(A)[n_1])$ acting by contraction.

We now consider two elements $\alpha$, $\beta$ of $\mathrm{Ext}^{(m_i,n_i)}_{A-}(A_0,A_0)$, $i=1,2$: the Yoneda product between them is represented by the composition of contractions
\[
\alpha\otimes \beta\mapsto (-1)^{(m_1+n_1)(m_2+n_2)} \beta\circ \alpha_n=\beta\alpha,
\]
viewed as an element of $\mathrm{hom}_{A}(\mathrm K^{-m_1-m_2}(A),A_0[n_1+n_2])\cong (A^!_{m_1+m_2})_{n_1+n_2}$, therefore, the Yoneda product is represented by the opposite product in $A^!$. 

We can finally summarize all arguments so far in the following
\begin{Thm}\label{t-kosz-rel}
For a finite-dimensional graded vector space $Y$, admitting a decomposition $Y=X_1\oplus X_2$, there is an isomorphism 
\[
\mathrm{Ext}_{\mathrm S(Y^*)-}^\bullet(\mathrm S(X_1^*),\mathrm S(X_1^*))^\mathrm{op}\cong \mathrm S(X_1^*\oplus X_2[-1]),
\]
of bigraded algebras w.r.t.\ the Koszul and internal grading.
\end{Thm}
Of course, the arguments above, with due modifications, hold true also when replacing left modules by right modules.

\medskip

We also have an isomorphism 
\[
\widehat{\mathrm{Ext}}_{\widehat{\mathrm S}(Y^*)-}^\bullet(\widehat{\mathrm S}(X_1^*),\widehat{\mathrm S}(X_1^*))^\mathrm{op}\cong \widehat{\mathrm S}(X_1^*\oplus X_2[-1]),
\]
of bigraded algebras; here, it is better to think of $\widehat{\mathrm{Ext}}_{\widehat{\mathrm S}(Y^*)-}^\bullet(\widehat{\mathrm S}(X_1^*),\widehat{\mathrm S}(X_1^*))^\mathrm{op}$ as e.g.\ the cohomology of the DG algebra $\widehat{\mathrm{End}}_{\widehat{\mathrm S}(Y^*)-}^\bullet(\widehat{\mathrm S}(X_1^*))$, with obvious $A_\infty$-structures.

\subsubsection{The proof of Keller's condition}\label{sss-6-3-4}
It is clear that the GAs $A$ and $B$ and the graded vector space $K$ from Subsection~\ref{ss-5-2} fit into the setting of Subsubsection~\ref{sss-6-3-1}; we consider here $K$ as an $A$-$B$-bimodule, where the actions are simply given by multiplication, followed by restriction.

We now recall the $A_\infty$-$A$-$B$-bimodule structure from Subsection~\ref{ss-5-2}.
\begin{Lem}\label{l-A_inf-triv}
For the structure maps~\eqref{eq-A_inf-bimod}, Subsection~\ref{ss-5-2}, hold the triviality conditions 
\[
\mathrm d_K^{0,n}=\mathrm d_K^{m,0}=0,\ \text{if $n,m\geq 2$}. 
\]
Further, $\mathrm d_K^{0,1}$, resp.\ $\mathrm d_K^{1,0}$, endow $K$ with the structure of a right $B$-module, resp.\ left $A$-module, simply given by multiplication followed by restriction: in particular, the $A_\infty$-$A$-$B$-structure on $K$ restricts to the above left $A$- and right $B$-module structures.
\end{Lem}
\begin{proof}
We recall from Subsection~\ref{ss-5-2} the construction of~\eqref{eq-A_inf-bimod}: if e.g.\ we consider the Taylor component
\[
\mathrm d_K^{0,n}(k|b_1|\cdots|b_n)=\sum_{\Gamma\in \mathcal G_{0,1+n}}\mu^K_{1+n}\left(\int_{\mathcal C_{0,1+n}^+}\omega_\Gamma^K(k|b_1|\cdots|b_n)\right).
\]

The discussion on admissible graphs in Subsection~\ref{ss-5-2} imply that a general admissible graph $\Gamma$ in the previous sum has no edges: the corresponding integral is thus non-trivial, only if the dimension of the corresponding configuration space is $0$, which happens exactly when $n=1$.

In such a case, $\mathrm d_K^{0,1}$ is simply given by multiplication followed by restriction on $K$, since there is no integral contribution.
\end{proof}
We observe that Lemma~\ref{l-A_inf-triv} implies that the left $A_\infty$-module structure on $K$, coming by restriction from the $A_\infty$-$A$-$B$-bimodule structure, is the standard one, as well as the right $A_\infty$-module structure; on the other hand, the $A_\infty$-$A$-$B$-bimodule structure is {\bf not} the standard one.
In particular, if we take the bar-cobar construction on $K$, for the left $A$-module structure we get a resolution of $K$, as well as for the right $B$-module structure; however, we do {\bf not} get a resolution of $K$ as an $A$-$B$-bimodule.

Lemma~\ref{l-A_inf-triv} implies, in particular, that the cohomology of $\underline{\mathrm{End}}_{-B}(K)$ coincides with $\mathrm{Ext}^\bullet_{-B}(K,K)$, the latter being the derived functor of $\mathrm{Hom}_{-B}(\bullet,K)$ in the category $\texttt{Mod}_{-B}$.
It is also clear that the graded algebra structure on $\underline{\mathrm{End}}_{-B}(K)$ induces the opposite of the Yoneda product on $\mathrm{Ext}^\bullet_{-B}(K,K)$, see e.g.~\cite{Rine} for a direct computational approach to the Yoneda product.

We know from Subsection~\ref{ss-3-1} that $\mathrm L_A$ is an $A_\infty$-algebra morphism from $A$ to $\underline{\mathrm{End}}_{-B}(K)$: in particular, since the cohomology of the $A_\infty$-algebra $A$ coincides with $A$ itself, $\mathrm L_A$ descends to a morphism of GAs from $A$ to $\mathrm{Ext}^\bullet_{-B}(K,K)^\mathrm{op}$, where the product on $\mathrm{Ext}^\bullet_{-B}(K,K)^\mathrm{op}$ is the opposite of the Yoneda product.
\begin{Prop}\label{p-keller-l}
We consider $A$, $B$ and $K$ as in Subsection~\ref{ss-5-2}, with the corresponding $A_\infty$-algebra structures and $A_\infty$-$A$-$B$-bimodule structure respectively, then the left derived $A$-action $\mathrm L_A$ is a quasi-isomorphism.
\end{Prop}
\begin{proof}
By the previous arguments, $\mathrm L_A$ descends to a morphism of GAs from $A$ to $\mathrm{Ext}^\bullet_{-B}(K,K)$; using the notation from Subsection~\ref{ss-5-2}, the GA $A$ is generated by the commuting variables $\{x_i\}$, for $i$ in $(I_1\cap I_2)\sqcup (I_1\cap I_2^c)$, and the anti-commuting variables $\{\partial_{x_i}\}$, $i$ in $(I_1^c\cap I_2)\sqcup (I_1^c\cap I_2^c)$.

On the other hand, as a corollary of Theorem~\ref{t-kosz-rel}, there is an isomorphism of GAs $\mathrm{Ext}_{-B}^\bullet(K,K)\cong A$: namely, $B=\mathrm S(Y^*)$, for $Y^*=(U\cap V)^*\oplus (U^\perp\cap V)^*\oplus (U\cap V^\perp)[-1]\oplus (U+V)^\perp[-1]$, and we set $X_1=(U\cap V)\oplus \left((U+V)^\perp\right)^*[-1]$, $X_2=(U^\perp\cap V)\oplus (U\cap V^\perp)^*[-1]$.

We will now prove that $\mathrm L_A$ is the identity map of $A$, by evaluating $\mathrm L_A$ on the generators of $A$.

We consider first $x_i$, for $i$ in $I_1\cap I_2$: the Taylor components of $\mathrm L_A^1(x_i)$ are given by 
\[
\mathrm L_A^1(x_i)^m(k|b_1|\cdots|b_n)=\mathrm d_K^{1,n}(x_i|k|b_1|\cdots|b_n)=\sum_{\Gamma\in\mathcal G_{0,1+1+n}}\mu_{1+1+n}^K\left(\int_{\mathcal C_{0,1+1+n}^+}\omega^K_\Gamma(x_i|k|b_1|\cdots|b_n)\right).
\]
An admissible graph $\Gamma$ yielding a non-trivial contribution to the previous expression has at most one edge: since $n=|\mathrm E(\Gamma)|\geq 1$, we have only two possibilities, either $i)$ $\Gamma$ has two vertices of the second type and no edge, or $ii)$ $\Gamma$ has three vertices of the second type and one edge.
Pictorially, 
\bigskip
\begin{center}
\resizebox{0.7 \textwidth}{!}{\input{left-x-adm.pstex_t}}\\
\text{Figure 13 - The only two admissible graphs contributing to $\mathrm L_A^1(x_i)$} \\
\end{center}
\bigskip
In case $ii)$, we get
\[
\mathrm L_A^1(x_i)^1(k|b_1)=\mathrm d^{1,1}_K(x_i|1|b_1)=\left(\int_{\mathcal C_{0,3}^+}\omega^{+,-}\right)(-1)^{|k|}k \left(\iota_{\mathrm d x_i}b_1\right)\vert_K,
\]
and since $b_1$ contains poly-vector fields normal w.r.t.\ $V$, the contraction w.r.t.\ $\mathrm d x_i$ annihilates $b_1$.
Thus, we are left with case $i)$, whence immediately
\[
\mathrm L_A^1(x_i)^0(k)=x_i\, k.
\]
 
We consider then $x_i$, for $i$ in $I_1\cap I_2^c$: again, we have to consider only $\mathrm L_A^1(x_i)^0$ and $\mathrm L_A^1(x_i)^1$.
In the first case, the contribution is trivial, because $\mathrm L_A^1(x_i)_0(k)$ is simply restriction on $K$ of the product $x_i\, k$.
We are left with $\mathrm L_A^1(x_i)^1(k|b_1)$: by construction,
\[
\mathrm L_A^1(x_i)^1(k|b_1)=\left(\int_{\mathcal C_{3,0}^+}\omega^{+,-}\right)(-1)^{|k|}k \left(\iota_{\mathrm d x_i}b_1\right)\vert_K=(-1)^{|k|}k \left(\iota_{\mathrm d x_i}b_1\right)\vert_K,
\]
because the integral can be computed explicitly e.g.\ by choosing a section of $C_{0,3}^+$, which fixes the middle vertex to $0$, and the left-most one to $-1$, and using the explicit formul\ae\ for the $4$-colored propagators, see Subsubsection~\ref{sss-4-3-2}, and is equal to $1$.

We consider $\partial_i=\partial_{x_i}$, for $i$ in $I_1^c\cap I_2$: then,
\[
\mathrm L_A^1(\partial_i)^n(k|b_1|\cdots|b_n)=\mathrm d_K^{1,n}(\partial_i|k|b_1|\cdots|b_n)=\sum_{\Gamma\in\mathcal G_{0,1+1+n}}\mu_{1+1+n}^K\left(\int_{\mathcal C^+_{0,1+1+n}}\omega^K_\Gamma(\partial_i|k|b_1|\cdots|b_n)\right).
\]
The arguments of Subsection~\ref{ss-5-2} imply that the admissible graphs in the previous formula have at most one edge: thus, only two graphs can contribute possibly non-trivially, either $i)$ the only graph with two vertices of the second type and no edge or $ii)$ the only graph with three vertices of the second type and one edge, pictorially
\bigskip
\begin{center}
\resizebox{0.7 \textwidth}{!}{\input{left-partial-adm.pstex_t}}\\
\text{Figure 14 - The only two admissible graphs contributing to $\mathrm L_A^1(\partial_i)$} \\
\end{center}
\bigskip
We consider $|\mathrm E(\Gamma)|=0$: there is only one graph with two vertices of the second type and no edges, whose corresponding contribution vanishes, since we restrict to $K$.
On the other hand, for $|\mathrm E(\Gamma)|=1$, we have only one graph with three vertices of the second type, and one edge, whose contribution is
\[
\mathrm L_A^1(\partial_i)^1(k|b_1)=\left(\int_{\mathcal C_{3,0}^+}\omega^{-,+}\right)k\left(\partial_i(b_1)\right)\vert_K=k \left(\partial_i b_1\right)\vert_K,
\]
where the integral can be computed explicitly e.g.\ by choosing a section of $C_{0,3}^+$, which fixes the middle vertex to $0$, and the left-most one to $-1$.

Finally, we consider $\partial_i$, for $i$ in $I_1^c\cap I_2^c$: by the same arguments as above, we need only consider $\mathrm L_A^1(\partial_i)^0$ and $\mathrm L_A^1(\partial_i)^1$.
We first consider $\mathrm L_A^1(\partial_i)^1$: the computation in the previous case implies that $\mathrm L_A^1(\partial_i)^1(k|b_1)$ vanishes, since $b_1$ does not depend on variables $\{x_i\}$, $i$ in $I_1^c\cap I_2^c$.
Thus, we are left with $\mathrm L_A^1(\partial_i)^0$, which is simply left multiplication by $\partial_i$ by construction.

In the previous computations, $\mathrm L_A^1(\bullet)$ is regarded as an element either of $\mathrm{Hom}(K[1],K[1])$ or of $\mathrm{Hom}(K[1]\otimes B[1],K[1])$: more precisely, we view $\mathrm L_A^1(\bullet)$, in all four cases, as a representative of a cocycle in $\mathrm{Ext}^\bullet_{-B}(K,K)$ w.r.t.\ the bar resolution of $K$ as a right $A$-module.
To identify correctly $\mathrm L_A^1(\bullet)$ with an element of $A$, we still need a chain map from the bar resolution of $K$ to the Koszul resolution of $K$ as a right $B$-module, because of Subsubsection~\ref{sss-6-3-3}: in particular, we need the components from $\mathcal B_0^B(K)=K\otimes B$ to $\mathrm K^0(B)=B$ and from $\mathcal B_1^B(K)=K\otimes B\otimes B$ to $\mathrm K^{-1}(B)$.
(We notice that the abstract existence of such a chain map is guaranteed automatically by standard arguments of homological algebra; the same arguments imply that such a chain map is homotopically invertible.) 

Since $K$ is a subalgebra of $B$, the map $\mathcal B_0^B(K)\to\mathrm K^0(B)$ is obviously given by multiplication; the map $\mathcal B_1^B(K)\to\mathrm K^{-1}(B)$ is a consequence of Poincar\'e Lemma in a linear graded manifold, more explicitly
\[
\mathcal B_1^B(K)\ni (k|b_1|b_2)\mapsto (-1)^{|k|}k \left(\mathrm d y_i\int_0^1 (\partial_{y_i}b_1)(ty)\mathrm dt\right) b_2\in \mathrm K^{-1}(B), 
\] 
where $\{y_i\}$ denotes a set of linear graded coordinates (associated to the chosen coordinates $\{x_i\}$ on $X$) of the graded vector space $X_2$, and where we have hidden linear graded coordinates on $X_1$, because they are left untouched by integration or derivation.
Graded derivations and corresponding contraction operators act from the left to the right.

From the previous computations, we see that $\mathrm L_A^1(x_i)$, $i$ in $I_1\cap I_2^c$, and $\mathrm L_A^1(\partial_i)$, $i$ in $I_1^c\cap I_2$, act non-trivially only on elements of the form $(k|\partial_i|b_2)$, $i$ in $I_1\cap I_2^c$, and $(k|x_i|b_2)$, $i$ in $I_1^c\cap I_2$ respectively: the image of such elements in $\mathcal B_1^B(K)$ w.r.t.\ the previous map is $(-1)^{|k|}k\,\mathrm d y_i\, b_2$, where now $y_i$ is a standard coordinate, if $i$ is in $I_1^c\cap I_2$, or a coordinate of degree $-1$, if $i$ is in $I_1\cap I_2^c$.

Setting then $b_2=1$, the computations in Subsubsection~\ref{sss-6-3-3} imply the desired claim.
\end{proof}
The same arguments, with obvious due modifications, imply that $\mathrm R_B:B\to \underline{\mathrm{End}}_{A-}(K)^\mathrm{op}$ is also a quasi-isomorphism: in fact, the same kind of computations in the proof of Proposition~\ref{p-keller-l}, prove that $\mathrm R_B$ equals the identity map on $B$, identifying the cohomology of $\underline{\mathrm{End}}_{A-}(K)^\mathrm{op}$ with $\mathrm{Ext}_{A-}(K,K)^\mathrm{op}$ in the category of left $A$-modules.
Thus, Keller's condition~\ref{ss-3-3} for the $A_\infty$-algebras $A$ and $\underline{\mathrm{End}}_{-B}(K)^\mathrm{op}$ is verified, from which we can deduce that the projection $\mathrm p_B$ in Diagram~\ref{d-keller} is a quasi-isomorphism in virtue of Theorem~\ref{t-keller}; similarly, the projection $\mathrm p_A$ is also a quasi-isomorphism, whence the commutativity of Diagram~\eqref{d-keller} implies that $\mathcal U$ is a quasi-isomorphism.

Equivalently, the Taylor component $\mathcal U^1$ is a Hochschild--Kostant--Rosenberg-type quasi-isomorphism from $T_\mathrm{poly}(X)$ to the cohomology of the Hochschild cochain complex 
$\left(\widetilde{\mathrm C}^\bullet(\texttt{Cat}_\infty(A,B,K),[\mu,\bullet]\right)$, where $\mu$ is the structure of $A_\infty$-category on $\texttt{Cat}_\infty(A,B,K)$, described in Subsection~\ref{ss-5-2}: 
from the discussion in Subsection~\ref{ss-2-1}, the HKR quasi-isomorphism has three components, $\mathcal U^1_A$, $\mathcal U^1_B$ and $\mathcal U_K^1$.
All three components can be described explicitly in terms of admissible graphs: the components $\mathcal U_A^1$ and $\mathcal U_B^1$ have been already described explicitly in~\cite{CF} in the framework of a formality 
result for graded manifolds.

On the other hand, the third component $\mathcal U^1_K:(T_\mathrm{poly}(X),0)\to \left(\widetilde{\mathrm C}^\bullet(A,B,K),[\mathrm d_K,\bullet]\right)$ is new.
By construction,
\[
\mathcal U^1_K(\gamma)(a_1|\cdots|a_m|k|b_1|\cdots|b_n)=\sum_{\Gamma\in \mathcal G_{1,m+1+n}}\mathcal O_\Gamma^K(\gamma|a_1|\cdots|a_m|k|b_1|\cdots|b_n).
\] 
Since the dimension of the configuration space $\mathcal C_{1,m+1+n}^+$ equals $m+n+1$, only those admissible graphs $\Gamma$ in $\mathcal G_{1,m+1+n}$ with $|\mathrm E(\Gamma)|=m+n+1$ yield possibly non-trivial contributions to the previous sum: such graphs can be of two types, $i)$ HKR-graphs, i.e.\ there are no edges in such graphs between vertices of the second type (hence, all edges connect the only vertex of the first type, corresponding to the multi-vector field $\gamma$, with the vertices of the second type), or $ii)$ HKR-$A_\infty$-graphs, i.e.\ graphs which contain (possibly multiple) edges connecting vertices of the second type, edges connecting the only vertex of the first type to vertices of the second type, and at most $1$ loop at the only vertex of the first type.

We observe that, for an admissible graph $\Gamma$ of type $(1,m+1+n)$ to yield a non-trivial contribution to the previous expression, the only vertex of the first type must be at least bivalent (i.e.\ there are at least two edges departing of incoming from this vertex).
Because of similar reasons, there is no $0$-valent vertex of the second type (i.e.\ a vertex of the second type, which is the initial or the final point of no edges).

Pictorially, the component $\mathcal U^1_K$ of the HKR-type quasi-isomorphism $\mathcal U$ of Theorem~\ref{t-form-cat} is a sum of the following two types of graphs:
\bigskip
\begin{center}
\resizebox{0.7 \textwidth}{!}{\input{HKR-gr.pstex_t}}\\
\text{Figure 15 - Two possible admissible graphs of type $(1,7)$ contributing to $\mathcal U_K^1$} \\
\end{center}
\bigskip

\section{Maurer--Cartan elements, deformed $A_\infty$-structures and Koszul algebras}\label{s-7}

Here an below we consider the case when $K$ has finitely many non-trivial taylor components. 

In Section~\ref{s-6}, we have constructed an $L_\infty$-quasi-isomorphism $\mathcal U$ from $T_\mathrm{poly}(X)$ to $\mathrm{C}^\bullet(\texttt{Cat}_\infty(A,B,K)$.

We consider a formal parameter $\hbar$: the ring $k_\hbar=k[\![\hbar]\!]$ is a complete topological ring, w.r.t.\ the $\hbar$-adic topology.
Accordingly, we denote by $T^\hbar_\mathrm{poly}(X)$ the trivial deformation $T_\mathrm{poly}(X)[\![\hbar]\!]$, where the Schouten--Nijenhuis bracket is extended to $T^\hbar_\mathrm{poly}(X)$ $k_\hbar$-linearly, and by $A_\hbar$, $B_\hbar$ and $K_\hbar$ the trivial $k_\hbar$-deformations of $A$, $B$ and $K$ respectively as in Subsection~\ref{ss-5-2}, where the GA-structures on $A$ and $B$ and the $A_\infty$-$A$-$B$-bimodule structure is extended $k_\hbar$-linearly to the respective algebras and modules.

In this framework, a $\hbar$-dependent MCE of $T^\hbar_\mathrm{poly}(X)$ is defined to be a $\hbar$-dependent polynomial bivector $\pi_\hbar$, which satisfies the Maurer--Cartan equation $[\pi_\hbar,\pi_\hbar]=0$.
The $\hbar$-formal Poisson bivector $\pi_\hbar$ is assumed to be of the form $\pi_\hbar=\hbar\pi_1+\mathcal O(\hbar^2)$: in particular, the Maurer--Cartan equation translated into a (possibly) infinite set of equations for the components $\pi_n$, $n\geq 1$, e.g.\ $\pi_1$ is a standard Poisson bivector on $X$.

Since $\mathcal U$ is an $L_\infty$-morphism, the image of $\pi_\hbar$ w.r.t.\ (the $k_\hbar$-linear extension of) $\mathcal U$ is also a MCE of $\mathrm{C}^\bullet(\texttt{Cat}_\infty(A,B,K))$, i.e.
\[
\mathcal U(\pi_\hbar)=\sum_{n\geq 1}\frac{1}{n!}\mathcal U^n(\underset{\text{$n$-times}}{\underbrace{\pi_\hbar|\cdots|\pi_\hbar}}).
\]
Again, the MCE $\mathcal U(\pi_\hbar)$ splits into three components, which we denote by $\mathcal U_A(\pi_\hbar)$, $\mathcal U_B(\pi_\hbar)$ and $\mathcal U_K(\pi_\hbar)$, viewed as elements of $\mathrm C^1(A_\hbar,A_\hbar)$, $\mathrm C^1(B_\hbar,B_\hbar)$ and $\mathrm C^1(A_\hbar,B_\hbar,K_\hbar)$.

\subsection{Deformation quantization of quadratic Koszul algebras}\label{ss-7-1}
We assume that we are in the framework of Subsection~\ref{ss-5-2}, where now $U=\{0\}$ and $V=X$, whence $A=\mathrm \wedge(X)$, $B=\mathrm S(X^*)$ and $K=k$: $A$ and $B$ are once again regarded as GAs, and $K$ is endowed with the (non-trivial) $A_\infty$-$A$-$B$-bimodule structure described in Subsection~\ref{ss-5-2}.

Furthermore, Theorem~\ref{t-kosz-rel}, Subsubsection~\ref{sss-6-3-3}, yields the well-known Koszul duality between $A$ and $B$, i.e.\
\[
\mathrm{Ext}^\bullet_{A-}(K,K)=B,\ \mathrm{Ext}^\bullet_{-B}(K,K)=A,
\]
where $K$ is viewed as a left $A$-module and right $B$-module respectively, as a consequence of Lemma~\ref{l-A_inf-triv}, Subsubsection~\ref{sss-6-3-4}. 

The Koszul complex of $A$ in the category ${}_A\texttt{GrMod}$ identifies with the deRham complex of $X$, with differential given by contraction w.r.t.\ the Euler field of $X$, as can be readily verified by repeating the arguments of Subsubsection~\ref{sss-6-3-1} in the present situation: in particular, the Koszul complex is acyclic, whence $A$ and $B$ are Koszul algebras over $k$.

We recall that the property of a non-negatively graded algebra $A$ over a field $K=A_0$ (more generally, over a semisimple ring $K=A_0$) of being Koszul, is equivalent to the existence of a (projective or free) resolution of $K$ in the category of graded right $A$-modules, whose component of cohomological degree $p$ is concentrated in internal degree $p$ (``internal'' refers to the grading in the category $\widehat{\texttt{GrMod}_k}$). 

For our purposes, we are interested in another {\em criterion} for a non-negatively graded algebra of being Koszul: namely, $A$ is a Koszul algebra, if and only the $\mathrm{Ext}_{A-}^\bullet(K,K)$-groups are concentrated in bidegree $(i,-i)$, $i\geq 0$.
We observe that the Koszul property implies that $A$ is quadratic, see e.g.~\cite{BGS, Sh} for details.

For a very detailed discussion of Koszul algebras, we refer to~\cite{BGS}; still, for a better understanding of the upcoming computations, we develop the above {\em criterion} in some details.

The graded Bar resolution of $K$ in the category of graded left $A$-modules, denoted by $\mathcal B^{A,+}_\bullet(K)$, is defined {\em via}
\[
\mathcal B_p^{A,+}(K)=A\otimes A_+^{\otimes p}\otimes K,
\]
where $A_+=\bigoplus_{n\geq 1}A_n$ and the tensor products have to be understood over the ground field $k$; the differential is a slight modification of the standard bar-differential.

By the definition of the category ${}_A\texttt{GrMod}$, we have
\begin{equation}\label{eq-bi-gr}
\mathrm{Hom}_{A-}(\mathcal B_p^{A,+}(K),K)=\bigoplus_{q\in \mathbb Z}\mathrm{hom}_{A-}(\mathcal B_p^{A,+}(K),K[q]).
\end{equation}

The differential on $\mathcal B_\bullet^{A,+}(K)$ has homological degree $1$ and Koszul degree $0$, where now the Koszul grading refers to the non-negative degree on $\mathcal B_\bullet^{A,+}(K)$ coming from the grading of $A$; by duality, $\mathrm{Hom}_{A-}(\mathcal B_p^{A,+}(K),K)$ has a differential of bidegree $(1,0)$, where the first, resp.\ second, grading is the cohomological, resp.\ Koszul, one.

Hence, we have a natural bigrading on $\mathrm{Ext}_{A-}^\bullet(K,K)$, inherited from Identity~\eqref{eq-bi-gr}. 
Further, since $K$ is concentrated in Koszul degree $0$, $K[q]$ is concentrated in degree $-q$. 
Since by construction $\mathcal B_p^{A,+}(K)$ is concentrated in Koszul degree bigger or equal than $p\geq 0$, it follows immediately that in general $-q\geq p$, i.e. $\mathrm{Ext}_{A-}^\bullet(K,K)=\bigoplus_{p+q\leq 0} \mathrm{ext}_{A-}^p(K,K[q])$.

In particular, the same arguments leading to the bigrading of $\mathrm{Ext}_{A-}^\bullet(K,K)$ yield, assuming $A$ is a Koszul algebra, the following condition on the bigrading:
\begin{equation}\label{eq-kosz-bi-gr}
\mathrm{Ext}_{A-}^\bullet(K,K)=\bigoplus_{p+q=0}\mathrm{Ext}_{A-}^{p,q}(K,K)\overset{!}=\bigoplus_{p+q=0}\mathrm{ext}^p_{A-}(K,K[q]).
\end{equation}
For a proof of the converse statement, we refer again to~\cite{BGS}.

Lemma~\ref{l-A_inf-triv} implies that the cohomology of $\underline{\mathrm{End}}_{A-}(K)$ identifies with $\mathrm{Ext}^\bullet_{A-}(K,K)$ in the category ${}_A\texttt{Mod}$, and, similarly, the cohomology of $\underline{\mathrm{End}}_{-B}(K)$ identifies with $\mathrm{Ext}^\bullet_{-B}(K,K)$ in the category $\texttt{GrMod}_{B}$.

For computational reasons, we choose a set of linear coordinates $\{x_i\}$, $i=1,\dots,d$, on $X$: thus, $A$ is generated by $\{x_i\}$ and $B$ is generated by $\{\partial_{x_i}=\partial_i\}$, $i=1,\dots,d$. 

The chain map from $\mathcal B_\bullet^B(K)$ to $\mathrm K^{-\bullet}(B)$ used in the proof of Proposition~\ref{p-keller-l}, Subsubsection~\ref{sss-6-3-4}, simplifies considerably: in particular, the image of $(1|b_1|b_2)$ in $\mathcal B_1^B(K)$ equals $-\left(\int_0^1\left(\partial_i b_1\right)\!(t x)b_2(x)\mathrm dt\right)\mathrm d x_i$ in $\mathrm K^{-1}(B)$.
\begin{Prop}\label{p-keller-sym}
The left derived action $\mathrm L_A$ descends to an isomorphism from $A$ to $\bigoplus_{p\geq 0}\mathrm{Ext}^{p,-p}_{-B}(K,K)$.
\end{Prop}
\begin{proof}
Adapting to the present situation the arguments of the proof of Proposition~\ref{p-keller-l}, Subsubsection~\ref{sss-6-3-4}, we find
\begin{equation}\label{eq-left-gen}
\mathrm L_A^1(\partial_i)^n(1|b_1|\cdots|b_n)=\begin{cases}
(\partial_{x_i}b_1)(0),& n=1,\\
0,& \text{otherwise}.
\end{cases}
\end{equation}
Viewing $1\otimes x_i\otimes 1$ as an element of Koszul degree $1$ in $\mathcal B_1^{B,+}(K)$, and recalling the previous discussion on the bigrading on $\mathrm{Ext}_{-B}^\bullet(K,K)$, the previous computation implies, in particular, that the image of $\mathrm L_A$ is contained in $\mathrm{Ext}_{-B}^{1,-1}(K,K)$.
Since $\mathrm L_A$ is an algebra morphism, and by the above {\em criterion} for Koszulness, $B$ is a Koszul algebra.

Finally, using the chain map from the bar resolution to the Koszul resolution of $K$ in ${}_B\texttt{GrMod}$, $\mathrm L_A^1(\partial_i)_1=-\iota_{\partial_i}$, where the expression on the right-hand side is viewed as a $B$-linear morphism from $\mathrm K_1(B)$ to $K$.
\end{proof}
We observe that, repeating these arguments {\em verbatim}, we may prove that $\mathrm R_B$ is an algebra isomorphism from $B$ to $\mathrm{Ext}^\bullet_{A-}(K,K)\cong A$, and that the image of $Y^*=A_1$ w.r.t.\ $\mathrm R_B$ is contained in the piece of bidegree $(1,-1)$ of $\mathrm{Ext}_{A-}^\bullet(K,K)$.

We now consider a $\hbar$-formal quadratic Poisson bivector on $X$, and the corresponding MCE $\mathcal U(\pi_\hbar)$ with components $\mathcal U_A(\pi_\hbar)$, $\mathcal U_B(\pi_\hbar)$ and $\mathcal U_K(\pi_\hbar)$.

It is easy to verify that $\mathcal U_A(\pi_\hbar)$ and $\mathcal U_B(\pi_\hbar)$ define associative products on $A_\hbar$ and $B_\hbar$ respectively: namely, e.g.\ $\mathcal U_A(\pi_\hbar)$ can be written explicitly as
\[
\mathcal U_A(\pi_\hbar)^m(\underset{\text{$m$-times}}{\underbrace{\bullet|\cdots|\bullet}})=\sum_{n\geq 1}\frac{1}{n!}\sum_{\Gamma\in \mathcal G_{n,m}}\mathcal O_\Gamma^A(\underset{\text{$n$-times}}{\underbrace{\pi_\hbar|\cdots|\pi_\hbar}}|\underset{\text{$m$-times}}{\underbrace{\bullet|\cdots|\bullet}}),
\]
borrowing notation from Subsection~\ref{ss-6-1}.

For a general admissible graph $\Gamma$ of type $(n,m)$, we have
\[
\mathcal O_\Gamma^A(\underset{\text{$n$-times}}{\underbrace{\pi_\hbar|\cdots|\pi_\hbar}}|\underset{\text{$m$-times}}{\underbrace{\bullet|\cdots|\bullet}})=\mu_{n+m}^A\!\left(\int_{\mathcal C_{n,m}^+}\omega_\Gamma^A(\underset{\text{$n$-times}}{\underbrace{\pi_\hbar|\cdots|\pi_\hbar}}|\underset{\text{$m$-times}}{\underbrace{\bullet|\cdots|\bullet}})\right).
\]
The integral in the previous expression on the right-hand side is non-trivial, only if the degree of the integrand equals $2n+m-2$, which is the dimension of $\mathcal C^+_{n,m}$.

The degree of the integrand equals $2n$, since we restrict to $A$ by means of the multiplication operator $\mu^A_{n+m}$, and since no edge in this situation can depart from vertices of the second type, if the corresponding contribution to the previous expression is non-trivial: this forces $m=2$.

We may thus consider $\mu_A+\mathcal U_A(\pi_\hbar)$, where $\mu_A$ is the $K_\hbar$-linear extension of the product on $A$ to $A_\hbar$: it is easy to verify that it defines an associative product $\star_A$ on $A_\hbar$, which, for $\hbar=0$, reduces to the standard product on $A$. 
Similar arguments imply that $\mathcal U_B(\pi_\hbar)$ defines an associative product $\star_B$ on $B_\hbar$, which reduces, for $\hbar=0$, to the standard product on $B$.

Furthermore, the expressions $\mathrm d_{K_\hbar}^{m,n}=\mathrm d_K^{n,m}+\mathcal U_K(\pi_\hbar)^{m,n}$, for non-negative integers $m$, $n$, define an $A_\infty$-$A_\hbar$-$B_\hbar$-bimodule structure on $K_\hbar$, which reduces, for $\hbar=0$, to the $A_\infty$-$A$-$B$-bimodule structure on $K$ described in Subsection~\ref{ss-5-2}.
\begin{Lem}\label{l-taylor-0-m}
The Taylor components $\mathrm d_{K_\hbar}^{m,n}$ satisfy the following triviality conditions:
\[
\mathrm d_{K_\hbar}^{m,0}=\mathrm d_{K_\hbar}^{0,n}=0,\ \text{if either $m=n=0$ or $m,n\geq 2$}.
\]
\end{Lem}
\begin{proof}
We consider exemplarily a Taylor component $\mathrm d_{K_\hbar}^{0,n}$, for $n\geq 0$: more explicitly,
\[
\mathrm d_{K_\hbar}^{0,n}(1|b_1|\cdots|b_n)=\sum_{l\geq 0}\frac{1}{l!}\sum_{\Gamma\in\mathcal G_{l,1+n}}\mathcal O_\Gamma^K(\underset{\text{$l$-times}}{\underbrace{\pi_\hbar|\cdots|\pi_\hbar}}|1|b_1|\cdots|b_n),
\]
with the same notation as above.

For a general admissible graph $\Gamma$ of type $(l,1+n)$,  
\[
\mathcal O_\Gamma^K=\mu_{l+1+n}^K\!\left(\int_{\mathcal C^+_{l,1+n}}\omega_{\Gamma}^K\right).
\]
Such an operator gives a non-trivial contribution to $\mathrm d_{K_\hbar}^{0,n}$, only if $|\mathrm E(\Gamma)|=2l+n-1$, where $2n+l$ is the dimension of $\mathcal C^+_{l,1+n}$.
Since a general vertex of the first type of $\Gamma$ has at most two outgoing edges, and a general vertex of the second type has no outgoing edges, and since we restrict to $K$, it follows that $|\mathrm E(\Gamma)|=2l$, whence $n=1$.
Similar arguments imply the claim for $\mathrm d_{K_\hbar}^{m,0}$, when $m\geq 2$ or $m=0$.
\end{proof}
We now discuss the grading on the deformed algebras $A_\hbar$, $B_\hbar$: we recall that the corresponding undeformed algebras possess a natural grading.
\begin{Lem}\label{l-grad-def}
The natural grading of $A$ and $B$ is preserved by the associative products $\star_A$ and $\star_B$ respectively.
\end{Lem}
\begin{proof}
Exemplarily, we consider a general non-trivial summand in
\[
\mathcal U_A(\pi_\hbar)^2(a_1|a_2)=\sum_{n\geq 1}\frac{1}{n!}\sum_{\Gamma\in\mathcal G_{n,2}}\mathcal O_\Gamma^A(\underset{\text{$m$-times}}{\underbrace{\bullet|\cdots|\bullet}}|a_1|a_2),
\]
associated to an admissible graph $\Gamma$ of type $(n,2)$.

By the same arguments used in the proof of Lemma~\ref{l-taylor-0-m}, such a graph has the property $|\mathrm E(\Gamma)|=2n$, which, by construction of the operator $\mathcal O_\Gamma^A$, implies that $\mathcal O_\Gamma^A$ contains exactly $2n$-derivations.
Since the polynomial degree of the element $(\underset{\text{$m$-times}}{\underbrace{\bullet|\cdots|\bullet}}|a_1|a_2)$ equals $2n+\mathrm{deg}(a_1)+\mathrm{deg}(a_2)$, the claim follows directly, where $\mathrm{deg}(\bullet)$ denotes the polynomial degree, and we recall that $\pi_\hbar$ is a quadratic bivector. 

Similar arguments imply the claim for $B_\hbar$.
\end{proof}
As a consequence of Lemma~\ref{l-taylor-0-m}, $K_\hbar$ has a structure of left $A_\hbar$- and right $B_\hbar$-module and the degree-$0$-component of both $B_\hbar$ and $A_\hbar$ identifies with $K_\hbar$, where the degree is specified by Lemma~\ref{l-grad-def}: hence, the cohomology of $\underline{\mathrm{End}}_{A_\hbar-}(K_\hbar)^\mathrm{op}$ identifies with $\mathrm{Ext}^\bullet_{A_\hbar-}(K_\hbar,K_\hbar)^\mathrm{op}$ in the category ${}_{A_\hbar}\texttt{GrMod}$, and the product on $\underline{\mathrm{End}}_{A_\hbar-}(K_\hbar)^\mathrm{op}$ descends, by a direct computation, to the opposite of the Yoneda product on $\mathrm{Ext}^\bullet_{A_\hbar-}(K_\hbar,K_\hbar)$, where $(A_\hbar,\star_A)$ and $(B_\hbar,\star_B)$ are GAs in view of Lemma~\ref{l-grad-def}.
Similarly, the cohomology of $\underline{\mathrm{End}}_{-B_\hbar}(K_\hbar)$ identifies with $\mathrm{Ext}^\bullet_{-B_\hbar}(K_\hbar,K_\hbar)$ in $\texttt{GrMod}_{B_\hbar}$, and composition descends, again, to the opposite of the Yoneda product.
\begin{Rem}\label{r-def_end}
We have been very sketchy in the definition of e.g.\ $\underline{\mathrm{End}}_{-B_\hbar}(K_\hbar)$: in fact, we define it, as a graded vector space, as the direct sum of the homogeneous components of the $\hbar$-trivial deformation $\underline{\mathrm{End}}_{-B}(K)_\hbar$, where the product is extended $\hbar$-linearly and continuously w.r.t.\ $\hbar$-adic topology, but whose differential is now the graded commutator with the deformed differential $\mathrm d_{K_\hbar,B_\hbar}$.
Thus, in the previous identification, we should also write $\mathrm{Ext}^\bullet_{A-}(K,K)_\hbar^\mathrm{op}$, where the product is also extended $\hbar$-linearly and continuously w.r.t.\ the $\hbar$-adic topology, but we still keep the previous notation.
\end{Rem}
The Taylor components $\mathrm d_{K_\hbar}^{m,n}$, $m$, $n$ non-negative integers, define, by the arguments of Subsection~\ref{ss-3-1}, the left derived action $\mathrm L_{A_\hbar}$ of $A_\infty$-algebras from $A_\hbar$ to $\underline{\mathrm{End}}_{-B_\hbar}(K_\hbar)$, and similarly for the right derived action $\mathrm R_{B_\hbar}$: $\mathrm L_{A_\hbar}$ descends to an algebra morphism from $A_\hbar$ to $\mathrm{Ext}^\bullet_{-B_\hbar}(K_\hbar,K_\hbar)$.

Lemma~\ref{l-grad-def} yields a bigrading on $\mathrm{Ext}^\bullet_{-B_\hbar}(K_\hbar,K_\hbar)$ and on $\mathrm{Ext}^\bullet_{-B_\hbar}(K_\hbar,K_\hbar)$ in the respective categories by the previous arguments. 
\begin{Lem}\label{l-keller-def}
The left derived action $\mathrm L_{A_\hbar}$ maps $A_\hbar$ to $\bigoplus_{p\geq 0}\mathrm{Ext}^{p,-p}_{-B_\hbar}(K_\hbar,K_\hbar)$.
\end{Lem}
\begin{proof}
First of all, we consider 
\[
\begin{aligned}
&\mathrm L_{A_\hbar}^1(\partial_i)^n(1|b_1|\dots|b_n)=\mathrm d_{K_\hbar}^{1,n}(\partial_i|1|b_1|\cdots|b_n)=\sum_{l\geq 0}\frac{1}{l!}\sum_{\Gamma\in \mathcal G_{l,1+n+1}}\mathcal O_\Gamma^K(\underset{\text{$l$-times}}{\underbrace{\pi_\hbar|\cdots|\pi_\hbar}}|\partial_i|1|b_1|\cdots|b_n),
\end{aligned}
\]
$b_j$ in $B_\hbar$, $j=1,\dots,n$, for $n\geq 1$, using the same notation as above.

We consider a general admissible graph $\Gamma$ of type $(l,1+n+1)$, $l\geq 0$, $n\geq 0$: its contribution is
\[
\mathcal O_\Gamma^K(\underset{\text{$l$-times}}{\underbrace{\pi_\hbar|\cdots|\pi_\hbar}}|\partial_i|1|b_1|\cdots|b_n)=\mu_{n+2}^K\left(\int_{\mathcal C^+_{l,1+1+n}}\omega_\Gamma^K(\underset{\text{$l$-times}}{\underbrace{\pi_\hbar|\cdots|\pi_\hbar}}|\partial_i|1|b_1|\cdots|b_n)\right).
\]
The degree of the integrand equals $|\mathrm E(\Gamma)|$, which, by all previous discussions, equals $2l+1$; since the dimension of $\mathcal C^+_{l,n+2}$ is $2l+n$, the previous integral is non-trivial, only of if $n=1$.

Thus, it remains to consider only 
\[
\mathrm L_{A_\hbar}^1(\partial_i)^1(1|b_1)=\sum_{l\geq 0}\frac{1}{l!}\sum_{\Gamma\in \mathcal G_{l,1+1+1}}\mathcal O_\Gamma^K(\underset{\text{$l$-times}}{\underbrace{\pi_\hbar|\cdots|\pi_\hbar}}|\partial_i|1|b_1),\ b_1\in A_\hbar.
\]
For a general admissible graph $\Gamma$ in $\mathcal G_{l,1+1+1}$, we consider the element $\mathcal O_\Gamma^K(\underset{\text{$l$-times}}{\underbrace{\pi_\hbar|\cdots|\pi_\hbar}}|\partial_i|1|b_1)$ of $K_\hbar$: by construction, it is non-vanishing, only if its polynomial degree w.r.t.\ $\{x_j\}$ is $0$.
The arguments of the proof of Lemma~\ref{l-grad-def} imply that its degree in the symmetric part is $\mathrm{deg}(b_1)-1$, which is equal to $0$, only if $\mathrm{deg}(b_1)=1$, i.e.\ $a_1$ is a monomial of degree $1$.  

The claim follows.
\end{proof}
Further, it is follows immediately from previous discussions 
\[
\mathrm L_{A_\hbar}\vert_{\hbar=0}=\mathrm L_A,\ \mathrm{Ext}^\bullet_{-B_\hbar}(K_\hbar,K_\hbar)\vert_{\hbar=0}=\mathrm{Ext}^\bullet_{-B}(K,K).
\]
We also observe that all deformed structures are obviously $\hbar$-linear, in particular, the differential on $\underline{\mathrm{End}}_{-B_\hbar}(K_\hbar)$ is $\hbar$-linear.

Summarizing the previous results, we have an $\hbar$-linear morphism $\mathrm L_{A_\hbar}$ of DG algebras from $(A_\hbar,0,\star_A)$ to $\left(\underline{\mathrm{End}}_{-B_\hbar}(K_\hbar),[\mathrm d_{K_\hbar,B_\hbar},\bullet],\circ\right)$, which restricts to a quasi-isomorphism, when $\hbar=0$: then, a standard perturbative argument w.r.t.\ $\hbar$ implies that $\mathrm L_{A_\hbar}$ is also a quasi-isomorphism, i.e.\ Keller's condition is verified for $\mathrm L_{A_\hbar}$.

In virtue of Lemmata~\ref{l-grad-def} and~\ref{l-keller-def}, Keller's condition implies that $\mathrm{Ext}^\bullet_{-B_\hbar}(K_\hbar,K_\hbar)$ is concentrated in bidegrees $(p,-p)$, $p\geq 1$, whence it follows that $B_\hbar$ is a Koszul algebra over $K_\hbar$.

On the other hand, the same arguments imply the validity of Keller's condition for $\mathrm R_{B_\hbar}$: this, in turn, implies that $A_\hbar$ is a Koszul algebra.
\begin{Thm}\label{t-keller-def}
We consider the $d$-dimensional vector space $X=k^d$, and a $\hbar$-formal quadratic Poisson bivector $\pi_\hbar=\hbar\pi_1+\mathcal O(\hbar^2)$ on $X$; further, we set $A=\wedge(X)$, $B=\mathrm S(X^*)$ and $K=k$, with the $A_\infty$-structures discussed in Subsection~\ref{ss-5-2}.

Then, the MCE $\pi_\hbar$ defines, by means of the $L_\infty$-morphism $\mathcal U$ of Theorem~\ref{t-form-cat}, Subsection~\ref{ss-6-2}, GA-algebra structures on $A_\hbar$ and $B_\hbar$, and an $A_\infty$-$A_\hbar$-$B_\hbar$-bimodule structure on $K_\hbar$, which deforms $A$ and $B$ to Koszul algebras $A_\hbar$ and $B_\hbar$, which are again Koszul dual to each other, i.e.\
\[
\mathrm{Ext}^\bullet_{A_\hbar-}(K_\hbar,K_\hbar)^\mathrm{op}\cong B_\hbar,\ \mathrm{Ext}^\bullet_{-B_\hbar}(K_\hbar,K_\hbar)^\mathrm{op}\cong A_\hbar,
\] 
in the respective categories.
\end{Thm}
\begin{Rem}
We observe that Theorem~\ref{t-keller-def} is an alternative proof of the main result of~\cite{Sh}: the main differences lie in the fact that $i)$ we make use of Kontsevich's formality result in the framework examined in~\cite{CF}, and $ii)$ instead of deforming Koszul's complex of $A$ and $B$ to a resolution of $A_\hbar$ and $B_\hbar$, we consider already at the classical level (i.e.\ when $\hbar=0$) a non-trivial $A_\infty$-$A$-$B$-bimodule structure on $K=k$, which we later deform by means of a quadratic MCE in $T_\mathrm{poly}^\hbar(X)$.
\end{Rem}

\end{document}